\documentclass[12pt,reqno,twoside]{amsart}
\title[Distribution of lattice orbits]{Distribution of lattice orbits
on homogeneous varieties}  
\author{Alex Gorodnik}
\address{University of Michigan, Ann Arbor, MI 48109 {\tt
gorodnik@umich.edu}} 
\author{Barak Weiss}
\address{Ben Gurion University, Be'er Sheva, Israel 84105
{\tt barakw@math.bgu.ac.il}}
\usepackage[all]{xy}

\usepackage{amsmath}
\usepackage{amssymb}

\newcommand{\goth}[1]{{\mathfrak{#1}}}
\newcommand{\Q}{{\mathbb {Q}}}

\newcommand{\R}{{\mathbb{R}}}
\newcommand{\Z}{{\mathbb{Z}}}
\newcommand{\C}{{\mathbb{C}}}

\newcommand{\N}{{\mathbb{N}}}

\newcommand{\cl}{\overline}

\newcommand{\Ad}{{\operatorname{Ad}}}
\newcommand{\GL}{\operatorname{GL}}
\newcommand{\SL}{\operatorname{SL}}

\newcommand{\Id}{\operatorname{Id}}

\newcommand{\spa}{{\rm span}}

\newcommand{\interior}{{\rm int}}

\newcommand{\df}{{\, \stackrel{\mathrm{def}}{=}\, }}
\newcommand{\HG}{{H \backslash G}}

\newcommand{\smallcirc}{{\circ}}
\newcommand {\equ}[1]     {\eqref{#1}}

\newcommand{\til}{\widetilde}
\newcommand{\supp}{{\rm supp}}

\newcommand{\sm}{\smallsetminus}
\newcommand{\vre}{\varepsilon}
\newcommand{\Mat}{\mathrm{Mat}}

\newcommand {\ignore}[1]  {}

\font\comment = cmbx10 scaled \magstep0

\newtheorem{thm}{Theorem}[section]
\newtheorem{lemma}[thm]{Lemma}
\newtheorem{prop}[thm]{Proposition}
\newtheorem{cor}[thm]{Corollary}

\newtheorem{claim}{Claim}
\newtheorem{remark}[thm]{Remark}

\newtheorem{definition}[thm]{Definition}
\newtheorem{example}[thm]{Example}

\begin{document}
\date{\today}

\begin{abstract}
Given a lattice $\Gamma$ in a locally compact group $G$ and a closed
subgroup $H$ of $G$, one has a
natural action 
of $\Gamma$ on the homogeneous space $V = \HG$. For an increasing family
of finite subsets $\{\Gamma_T: T>0\},$ a dense orbit $v
\cdot \Gamma, \ v \in V$ and compactly supported function $\varphi$
on $V$, we consider the sums $S_{\varphi, v}(T) =
\sum_{\gamma 
\in \Gamma_T} \varphi (v \gamma).$ 
Understanding the asymptotic behavior of
$S_{\varphi, v}(T)$ is a delicate problem which has only been
considered for certain very special 
choices of $H$, $G$ and $\{\Gamma_T\}$. We develop a general abstract
approach to 
the problem, and apply it
to the case when $G$ is a Lie group and either $H$ or $G$ is
semisimple. When $G$ is a group of matrices equipped with a norm, we
have $S_{\varphi, v}(T) \sim \int_{G_T} \varphi(vg) \, dg, $ where
$G_T = \{g \in G: \| g\|<T\}$ and $\Gamma_T = G_T \cap \Gamma.$ We
also show that the asymptotics of $S_{\varphi, v}(T)$ is governed by 
$\int_V \varphi \, d\nu,$ where $\nu$ is an explicit limiting density
depending on the choice of $v$ and $\| \cdot \|$.  
\end{abstract}

\maketitle

\setcounter{tocdepth}{1}
\tableofcontents

\section{Introduction}
Let $V$ be a manifold equipped with a transitive right action of a Lie
group $G$, so 
that $V$ is identified with $\HG$ for some subgroup $H$ of $G$
stabilizing a point of $V$. Let $\Gamma$ be a lattice in $G$, that is,
a discrete subgroup of finite co-volume. In this
paper we study the asymptotic distribution of $\Gamma$-orbits. More
precisely, we fix a proper function $D: G \to [0,\infty)$, set
$G_T = \{g \in G : D(g)<T\}$, and study, for an
arbitrary $\varphi \in C_c(V)$ and $v \in V$, the asymptotic behavior
of the sum 
$$S_{\varphi, v}(T) = \sum_{\gamma \in \Gamma \cap G_T} \varphi(v \cdot
\gamma)$$ as $T \to \infty.$

We have not assumed that there is an invariant measure on $V$, and
certainly not a finite one. Thus this problem does not belong to the
classical framework of ergodic theory, and has many surprising
features. 

\subsection{A simple case}
Perhaps the simplest nontrivial case is when 
$$V = \R^2 \sm \{0\}, \ \  G=\SL(2,\R), \ \ D(g) = \|g\|,$$
 where $$\left \|\left( \begin{matrix}a &
b \\ c & d \end{matrix} \right ) \right \| = \left(a^p + b^p + c^p +
d^p\right)^{1/p}, \, 1 \leq p< \infty.$$  
Then $H$ is (up to conjugation) 
the one-parameter group of upper-triangular unipotent matrices. 
The distribution of $\SL(2,\Z)$-orbits in $V$ was considered by Nogueira
\cite{Nogueira}, and the case of a general lattice by Ledrappier
\cite{Ledrappier}. Ledrappier showed that for $\varphi \in C_c(V)$ and
$v \in V$ such that $v \cdot \Gamma$ is dense in $V$, one
has:

$$S_{\varphi, v}(T) \sim \left( c_{\Gamma} \int
\frac{\varphi(Y)}{\| v \| \cdot \|Y\|} dY \right ) \, T,
$$
where $dY$ denotes Lebesgue measure on $\R^2$,
$\| \cdot \|$ is the $p$-norm on $\R^2$, and $c_{\Gamma}$ is an
explicit constant depending on the lattice $\Gamma$ (here and
throughout the paper, the notation $A(T) \sim B(T)$ means that
$\frac{A(T)}{B(T)} \to_{T \to \infty} 1$).

In this example, one already notices that the limiting measure
$\frac{1}{\|v\|} \, \frac{dY}{\|Y\|}$ is not $\Gamma$-invariant and
depends essentially on the choice of $v$ and the choice of norm used
on $G$. A computation with the haar measure of $G$ shows that in this
case the asymptotics of $S_{\varphi, v}(T)$ is the same as that of the
corresponding average along the (unique) $G$-orbit, namely 
$$S_{\varphi, v}(T) \sim \til{S}_{\varphi, v}(T), \mathrm{\ \ \ where \
\ } \til{S}_{\varphi, v}(T) = \int_{G_T} \varphi(v \cdot g)\, dm(g),
$$
and $m$ is haar measure on $G$.

We reconsider and generalize this example in \S 12.4.

\subsection{The general problem}
It is natural to try to extend this result to a general
case. Namely, given an arbitrary 
locally compact second countable topological group $G$, a lattice
$\Gamma$ in $G$, a transitive $G$-space $V$ and a proper function $D:
G \to [0, \infty)$, one would like to:

 \begin{itemize}
\item
Obtain an explicit expression for
the asymptotics of $S_{\varphi, v}(T)$.

\item
Show that 
$S_{\varphi, v}(T) \sim \til{S}_{\varphi, v}(T).$
\end{itemize}

 A number of special
cases of this program were carried out in 
work of the first-named author, all involving $G=\SL(n,\R), \,
D (g) = \sqrt{\mathrm{tr} \ ({}^t g \, g)}$ and a general lattice 
$\Gamma$: in \cite{gorodnik1}, the
case in which $V$ is the 
Furstenberg boundary of $G$, and in \cite{gorodnik2}, the case in which
$V$ is the space of $k$-frames in $\R^n, \, 1 \leq k < n.$ Additional
cases were studied in \cite{lp, Maucourant, GM, GO}.

Our goal in this paper is twofold. First we abstract some ideas of the
above-mentioned papers and develop a general axiomatic framework for
studying the problem. We specify certain explicit
conditions on $G, H, D, \Gamma$ under which the above problems can be
solved. These  
conditions are nontrivial to check, and in the second part of the
paper we study the conditions in two important cases: 
when $H$ is semisimple, or when $G$ is semisimple and the
$\Gamma$-action admits a finite invariant measure.

We use two strategies to analyze
$\Gamma$-orbits. First, when 
there is a finite $\Gamma$-invariant measure on $V$, the standard 
construction of an induced action gives a $G$-action with the same
asymptotic behavior as the given $\Gamma$-action. Secondly, one may
use `duality' to replace the question about $\Gamma$-orbits on $\HG$
with a problem about $H$-orbits on $G/\Gamma$ (a similar approach was
introduced in \cite{drs} and employed in \cite{EMS - Annals} to study the
{\em discrete} orbits). Although more general, 
the second strategy is more difficult to implement because the
question about $H$-orbits involves averaging along 
certain `skew-balls' with respect to the function $D$. 

Following either strategy, one is led
to a question about the distribution of an orbit of
a connected Lie group on a
finite measure homogeneous space. Such questions have been thoroughly
studied in recent years, building on Ratner's classification of measures
invariant under unipotent flows and subsequent work of many
authors. The precise result we require follows from work of Nimish Shah
\cite{Nimish: proc ind}. We refer the reader to \cite{KSS handbook}
for a recent 
and detailed account of this theory. 

Our strategy also requires a precise understanding of the asymptotics
of volumes of certain `skew balls' in $H$, with respect to a norm in a
linear representation,  or with respect to the natural $G$-invariant
metric on the symmetric space of $G$.  For example we determine
the asymptotics of the volume of the sets $\{g \in G: \| g \| < T\},$
where $G$ is a semisimple Lie group realized as a matrix group in
$\Mat_d(\R)$ and $\| \cdot \|$ is an arbitrary linear norm on
$\Mat_d(\R)$. Even this very natural question has only been settled in
the literature in certain special cases.

Precise formulation of our results appears in the next section. To
conclude this introduction, we list some applications which illustrate
the scope of our methods.

\subsection{Applications}
\subsubsection{Values of indefinite irrational quadratic forms on frames}
Recall that the Oppenheim conjecture proved by Margulis \cite{m89}
states that for a nondegenerate indefinite quadratic form $Q$ of dimension $d\ge 3$
that is not a scalar multiple of a rational form,
the set $Q(\mathbb{Z}^d)$ of its values at integer points is dense in $\mathbb{R}$.
A strengthening of this result was obtained by Dani and Margulis and by Borel and Prasad. 
Denote by $\mathcal{F}_d$ the space of unimodular frames in
$\mathbb{R}^d$, i.e. 
$$
\mathcal{F}_d=\left\{(f_1,\ldots, f_d): f_i\in\R^d, \,
\hbox{Vol}\left( \R^d/(\Z f_1 \oplus \cdots \oplus \Z f_d) \right)=1\right\}.
$$
Let $\mathcal{F}_d(\Z)$ be the space of integer unimodular frames.
For a frame $f\in \mathcal{F}_d$, denote by $\bar{Q}(f)$ the corresponding
Gram matrix: 
$$
\bar{Q}(f)=\left(Q(f_i,f_j)\right)_{i,j=1,\ldots,d} \in \Mat_d(\R).
$$
Then for a nondegenerate indefinite quadratic form $Q$ of dimension $d\ge 3$
that is not a scalar multiple of a rational form, it was observed in  \cite{dm89} and \cite{bp92} that the set
$\bar{Q}(\mathcal{F}_d(\Z))$ 
is dense in $\bar{Q}(\mathcal{F}_d)$. Our results
imply a quantitative 
strengthening  of this fact. Let
$$G = \{g \in \GL(d,\R): \det g = \pm 1\}.$$
Recall that the map 
\begin{equation}
\label{eq: identify}
g \in G \mapsto g \, \mathbf{e} \in \mathcal{F}_d,
\end{equation}
where $\mathbf{e}$ is the frame consisting of the standard basis
vectors, is a diffeomorphism. This endows $\mathcal{F}_d$ with a
natural measure $\mu$ coming from haar measure $m$ on $G$, which we
normalize by 
requiring that $m \left(G/G(\Z) \right)=1.$ Note that
$\bar{Q}(\mathcal{F}_d)$ is a submanifold of $\Mat_d(\R)$, consisting
of symmetric matrices of a given signature.  We have: 

\begin{cor}
\label{th_oppenheim}
Let $Q$ be as above. 
Fix any norm $\|\cdot\|$ on $\R^d$ and let $A \subset \bar{{Q}}(\mathcal{F}_d)$
be a bounded set whose boundary has measure zero. 
Then
\[
\# \, \left\{f\in\mathcal{F}_d(\Z):\|f\|<T, \, \bar{Q}(f)\in
A \right\}
\sim
\mu \left(\left\{f\in\mathcal{F}_d:\|f\|<T, \, \bar{Q}(f)\in
A \right\}\right)
\]
(here $\|(f_1, \ldots, f_d)\| = \max_{i=1, \ldots,
d}\|f_i\|$). Explicitly, these quantities are asymptotic to 
\[ \begin{array}{cc}
C T^{p(q-1)} &  \mathrm{if \ } p<q \\
C \log T \, T^{p(p-1)} &   \mathrm{if \ } p=q
\end{array}
\]
where $(p,q), \, 1 \leq p \leq q$ is the signature of $Q$ and
$C$ is an explicitly computable constant depending on $A$, $Q$ and $\|
\cdot \|$.
\end{cor}

Note that in contrast with the quantitative version of the Oppenheim
conjecture, proved by 
Eskin, Margulis and Mozes \cite{emm}, in this result there are no
difficulties arising in signatures $(2,1)$ and $(2,2)$.





\subsubsection{Dense projections of lattices}
Suppose $G$ is a semisimple non-simple group, and $\Gamma$ is an
irreducible lattice, i.e. the projection
of $\Gamma$ onto each non-cocompact factor of $G$ is dense. Our results imply a
quantitative version of this fact. Here we formulate two special cases. 

Let $G$ be a Lie group which is a direct product $G=SH$ of
its closed subgroups $S$ and $H$. Suppose $\Gamma$ is a lattice in $G$
such that $H\Gamma$ is dense in 
$G$. Denote by $P: G\to S$ the projection map. Suppose $H$ is simple,
and $\Psi_H : H \to \GL(V_H), \, \Psi_S: S \to \GL(V_S)$ are two
irreducible representations with compact kernels. 

Let $m$ be a haar measure on $G$, normalized by requiring that the covolume
of $\Gamma$ in $G$ is equal to 1. Let $\lambda$ and $\nu$ haar
measures on $H$ and $S$ respectively, with $m=\lambda \otimes \nu$.

\begin{cor}
\label{cor: dense projections}
 Let the notation be as above. 
\begin{itemize}
\item[(i)]
Let $\Psi : G \to \GL(V_H \oplus V_S)$ be the direct sum
representation of $G$, and for a linear norm $\| \cdot \|$ on $V_H
\oplus V_S$ let 
\begin{equation*}
H_T =\{h \in H: \|\Psi(h)\| <T\}, \ \ \ \Gamma_T = \{\gamma \in \Gamma
: \|\Psi(\gamma)\|<T\}.
\end{equation*}
For every $\varphi\in C_c(S)$ and every $s_0 \in S$,
$$
\frac{1}{\lambda(H_T)}\sum_{\gamma\in\Gamma_T} \varphi(s_0P(\gamma))
\mathop{\longrightarrow}_{T\to\infty}  \int_S \varphi \, d\nu.
$$
\item[(ii)]
Let $\Psi : G \to \GL(V_H \otimes V_S)$ be the tensor product
representation of $G$. Choose bases $\{v_1, \ldots, v_k\}, \, \{u_1,
\ldots, u_{\ell} \}$ of $V_H$ and $V_S$, and for $1 \leq p <\infty$
let $\| \cdot \|$ denote both the $p$-norm on 
$V_S$ associated with the basis $\{u_1, \ldots, u_{\ell} \}$ and the $p$-norm on
$V_H \otimes V_S$ associated with
the basis $\{v_i \otimes u_j\}$.  
Then  for every $\varphi\in
C_c(S)$ and every $s_0 \in S$,
$$
\frac{1}{\lambda(H_T)}\sum_{\gamma\in\Gamma_T}
\varphi(s_0P(\gamma))\mathop{\longrightarrow}_{T\to\infty} c\int_S\varphi(s)\,
\frac{d\nu(s)}{\|\Psi_S(s_0s)\|^m}, 
$$
where $m$ and $c$ are explicitly given positive constants. 
\end{itemize}
\end{cor}

\subsubsection{Lattice actions by translations}

Let $L$ be a noncompact simple connected Lie group, and $\Lambda$ and
$\Delta$ lattices in $L$. 
The group $\Lambda$ acts on $L/\Delta$ by left translation: if $\pi: L
\to L/\Delta$ is the quotient map, then the action is given by 
$$ \lambda \, \pi(g) = \pi(\lambda g).
$$
It is a simple application of Ratner's theorems that for any $x=\pi(g)
\in 
L/\Delta$, $\Lambda x$ is either finite (in case $\Lambda$ and $g
\Delta g^{-1}$ are commensurable) or dense. Note that a special
case of this observation was used recently by Vatsal \cite{va} in his study of
Heegner points. 
We obtain an equidistribution result for the dense $\Lambda$-orbits in
$L/\Delta$, generalizing a recent result of \cite{Hee}.

Let $m$ be the Haar measure on $L$, normalized so that $m \left(L/\Delta
\right)=1$. We denote the induced finite measures on $L/\Delta$ and
$L/\Lambda$ also by $m$.  
We assume that $L$ is a closed subgroup of $\GL(d,\R)$, fix a
norm $\|\cdot\|$ on $\Mat_d(\R)$, and set 
$$L_T = \{\ell \in L: \| \ell \| < T \}.
$$

\begin{cor}\label{cor_equi_appl2}
Let $\pi(g_0) \in L/\Delta$ be such that $\overline{\Lambda \pi(g_0)}=L/\Delta$.
Then for $\varphi\in C_c(L/\Delta)$,
\[
\frac{1}{m(L_T)}\sum_{\lambda\in\Lambda \cap L_T} \varphi(\lambda^{-1} \pi(g_0))
\mathop{\longrightarrow}_{T \to  \infty}
\frac{1}{m(L/\Lambda)}\int_{L/\Delta}\varphi \, dm.
\]

\end{cor}

\subsubsection{Lattice actions by automorphisms}
Our analysis also yields results for actions of lattices on manifolds
which do not admit a transitive action of the ambient Lie group. For
example, for $d \geq 2$, let $\Gamma = \SL(d, \Z)$ act by linear
isomorphisms on 
the torus $\mathbb{T} = \R^d / \Z^d$. 
It is well-known that an orbit $\Gamma x$ for this action is either
finite (if $x \in \mathbb{T}(\Q)$) or dense. Our results imply a
quantitative strengthening of this fact. 

Let $G=\SL(d,\R)$, let $m$ denote the haar measure on $G$, normalized
so that $m \left( G/\Gamma \right) =1,$ and let $\mu$ be the haar
probability measure on $\mathbb{T}$. 
We have:
\begin{cor}
\label{cor: application torus}
Fix a norm $\|\cdot\|$ on $\Mat_d(\R)$. Let $\varphi \in
C_c(\mathbb{T})$.
Then for any $x \in \mathbb{T} \sm \mathbb{T}(\Q),$
$$
\frac{\sum_{\| \gamma \|<T} \varphi (\gamma^{-1} x_0)}{\# \left\{ \gamma \in \Gamma
: \| \gamma \|<T \right\} }\mathop{\longrightarrow}_{T \to
\infty} \int_{\mathbb{T}} \varphi \, d\mu.
$$
\end{cor}

There is nothing special about the torus in this statement; the result
follows from the more general Theorem \ref{cor_auto_equi} which 
encompasses other actions, e.g. the action of $\Gamma$ on certain 
nilmanifolds. 


 \subsection{Acknowledgements}
We thank Amos Nevo and Ralf Spatzier for many useful and insightful 
conversations. Thanks also to Hee Oh, Dave Morris and Michael Levine. We thank the Center for 
Advanced Studies in Mathematics at Ben-Gurion University for funding
Gorodnik's visit to Be'er Sheva, when this work was conceived, and the
Technion's hospitality, which made additional
progress possible. The first author is partially supported by NSF grant
0400631.

\section{Notation and Statement of Results}
\subsection{The general setup}
\label{subsection: statements}

Let $G$ be a second countable locally compact noncompact topological
group, let $\Gamma$ be a 
lattice in $G$, and let $\pi: G
\to G/\Gamma$ be the natural quotient map.

There is a natural left-action of $G$ and any of its subgroups on
$G/\Gamma$ defined 
by $g_1\, \pi(g_2)=\pi(g_1g_2).$ 

Since it admits a lattice,
$G$ is unimodular, i.e. the haar measure $m$ is 
invariant under both left and right multiplication.
Let $m'$ denote the $G$-invariant measure on
$G/\Gamma$ induced by $m$, that is, 
$$m'(X) \df m\left(\Omega \cap \pi^{-1}(X)\right),$$ 
for some Borel fundamental domain $\Omega$ for the right-action of
$\Gamma$ on $G$.  
It follows from the invariance of $m$ that $m'$ is independent of the
choice of $\Omega$. Normalize $m$ so that $m'(G/\Gamma)=1$. 

By a {\em distance function} on $G$ we mean a function $G \to [0,
\infty)$ which is continuous and proper. Let $D$ be a
distance function on $G$. 
The words {\bf the general setup holds} mean that $G, \,
\Gamma, \, D,\, m, \, m'$ are as above. 

Fixing $D$, for a subset $L \subset G$ we write 
$$L_T \df \{g \in L: D(g)<T\}.$$

We list some hypotheses about this setup.

\begin{itemize}

\item[{\bf UC}]
{\bf Right uniform continuity of $D$.}
For any $\vre>0$ there is a neighborhood $\mathcal{U}$ of identity in
$G$ such that for all $g \in G, \, u \in \mathcal{U},$
\begin{equation}
\label{eq: rt uniform cont}
D(gu) < (1+\vre)D(g).
\end{equation}

\item[{\bf I1}]
{\bf Moderate volume growth for balls in $G$.} 
For any $\vre>0$ there are $\delta>0$ and $T_0$ such that for all
$T>T_0$:
\begin{equation}
\label{eq: I1}
m \left(G_{(1+\delta)T} \right) \leq (1+\vre) m \left(G_T
 \right). 
\end{equation}


\item[{\bf I2}]
{\bf $\Gamma$-points equidistributed in $G$ w.r.t. $D$.}
$$\# \Gamma_T \sim m(G_T).$$

\end{itemize}


\subsection{Inducing the action}
\label{subsection: inducing}
Suppose $X$ is a space on which $\Gamma$ acts on the right, preserving
a finite invariant measure. There
is a standard construction (see \S \ref{section: induction}) of a left
$G$-space $Y$ with a $G$-map $\pi_{G/\Gamma}: Y \to G/\Gamma$, such
that the fiber over $\pi(e)=[\Gamma] \in G/\Gamma$ is isomorphic to
$X$. In case $X =
\HG$, it
is simply the action of $G$ on the product $\HG \times
G/\Gamma$, given by
\begin{equation}
\label{eq: G action, induction}
g \cdot \left(\tau(g_2), \pi(g_1) \right) = \left(
\tau(g_2g^{-1}), \pi(gg_1)\right).
\end{equation}

The following reduces the study of asymptotic behavior of
$\Gamma$-orbits on $X$ to that of $G$-orbits on $Y.$

\begin{prop}
\label{prop: induction, special}
Suppose hypotheses I1, I2 and UC hold. Let $x_0 \in X$, and let $y_0
\in Y$ be the corresponding point in the fiber $\pi_{G/\Gamma}^{-1}
\left([\Gamma] \right)$. Suppose that 
there is a measure 
$\mu$ on $X$ such that:

\indent{$\left( * \right)$ \ \ For any $F \in C_c(Y)$,
\[
\frac{1}{m(G_T)} \int_{G_T} F(g^{-1} \cdot y_0)\, dm(g)
\mathop{\longrightarrow}_{T \to \infty} \int_{Y} F\, d\nu
,
\]}
where $d\nu =  d\mu \, dm'$.

Then for any $\varphi \in C_c(X),$
\begin{equation*}
\frac{1}{\# \Gamma_T} \sum_{\gamma
\in \Gamma_T} \varphi (x_0 \, \gamma) \mathop{\longrightarrow}_{T \to \infty} \int_X
\varphi \, d\mu.
\end{equation*}
\end{prop}

\subsection{The general setup, continued}
Suppose the general setup holds, let $H$ be a
closed subgroup of $G$, and let $\tau: G \to H \backslash G$ be the
natural quotient map. Let $\lambda$ be the left haar measure on
$H$. There is a natural right-action of $G$ on $H \backslash G$
defined by  $\tau(g_1)\, g_2= \tau(g_1g_2)$. We have the following
diagram.

\[
\xymatrix{
& G \ar[dl]_{\tau} \ar[dr]^{\pi}  \\
H \backslash G & & G/\Gamma}
\]

If $H$ has been defined, the words {\bf the general setup holds}
mean in addition that $H, \, \lambda$ and $\tau$ are as above. 

If $g_1, \, g_2, \, g \in G$ and $L \subset G$, we write
\begin{equation}
\label{eq: def skew balls}
L_T[g_1, g_2] 
 \df \{\ell \in L: \, D(g_1\ell g_2) <T \}, \ \ \ \
L_T[g] \df L_T[e,g].
\end{equation}

We sometimes call these `skew balls'. 
We list some additional hypotheses.

\begin{itemize}
\item[{\bf S}]
{\bf Locally continuous section.}
For any $x \in \HG$ there is a Borel map $\sigma : \HG \to G$
which is continuous in a neighborhood of $x$ and satisfies $\tau
\smallcirc \sigma = \Id_{\HG}.$

\item[{\bf D1}]
{\bf Uniform volume growth for skew-balls in $H$.} 
For any bounded
$B \subset G$ and any $\vre>0$ there are $T_0$ and $\delta>0$ such
that for all $T>T_0$ and all 
$g_1,g_2 \in B$ we have: 
\begin{equation}
\label{eq: A1}
\lambda
\left(H_{(1+\delta)T} [g_1,g_2] \right) \leq (1+\vre)\lambda \left(H_T
[g_1,g_2] \right). 
\end{equation}

\item[{\bf D2}]
{\bf Limit volume ratios.}
For any $g_1,\, g_2\in G$, 
the limit

\begin{equation}\label{eq_dual_as}
\alpha(g_1, g_2) \df \lim_{T\to\infty}
\frac{\lambda(H_T[g_1^{-1},g_2])}{\lambda(H_T)} 
\end{equation}
exists and is positive and finite.


\end{itemize}


\subsection{Duality.}
The following result establishes a link between the asymptotic
behavior of a $\Gamma$-orbit on $\HG$ and the behavior of a
corresponding $H$-orbit on $G/\Gamma$.

\begin{thm}
\label{thm: duality}
Suppose the general setup holds, and conditions S, D1, UC are
satisfied. Let $g_0 \in G$ and assume the following:

\indent{$\left( ** \right)$ \ \ For any 
$b \in G$ and any $F \in C_c(G/\Gamma)$, 
$$\frac{1}{\lambda \left(H_T [g_0, b]
\right)}\int_{H_T
[g_0, b]} F(h^{-1}\pi(g_0)) \, d\lambda(h)
  \mathop{\longrightarrow}_{T \to \infty} 
\int_{G/\Gamma} F(x)\, dm'(x) .
$$
}

Then for every non-negative compactly supported $\varphi$ on $\HG$,
\begin{equation}
\label{eq: to demonstrate}
\sum_{\gamma \in \Gamma_T} \varphi(\tau(g_0)\gamma)
\sim \int_{G_T} \varphi(\tau(g_0)g) \, dm(g).
\end{equation}

\end{thm}

\subsection{The limiting density}\label{subsection 2.3}
The conclusion of Theorem \ref{thm: duality} describes the asymptotics
of a $\Gamma$-orbit in terms of the asymptotics of a $G$-orbit.
In order to calculate the latter we need some more
terminology, and an additional
assumption. Note that we are calculating the asymptotics of a $G$-orbit in a space
on which $G$ acts transitively. Nevertheless the computation is not   
trivial.

We have not assumed that $H$ is unimodular, and therefore a
$G$-invariant measure on $\HG$ need not exist. To remedy this, we need
to discuss measures on $\HG$. We suppose that condition D2 is satisfied,
and define $\alpha(g_1, g_2)$ by (\ref{eq_dual_as}).

Let $Y$ be a lift of $\HG$ to $G$. That is, $Y$ is a Borel subset of
$G$ such that the product map 
$$
H\times Y\to G:(h,y) \mapsto hy
$$
is a Borel isomorphism.
Since the measures $m$ and $\lambda$ are left $H$-invariant, 
\begin{equation}
\label{eq: measure on fibers}
dm(hy) = d\lambda(h) \, d \nu_Y(y)
\end{equation}
for some Borel measure $\nu_Y$ on $Y$ 
. Moreover (\ref{eq: measure on
fibers}) determines $\nu_Y$ uniquely. When $H$ is
unimodular, $\nu_Y$ is identified with a
measure on $\HG$ which we denote by $\nu_{\HG}$. It is independent of
$Y$ and is the unique (up
to scaling) $G$-invariant measure on $\HG$.





For $x_0=\tau(g_0) \in \HG$, define a measure $\nu_{x_0}$ on $Y$ by 
\begin{equation}\label{eq_nu_g}
d \nu_{x_0}(y) \df \alpha(g_0,y) \, d\nu_Y(y).
\end{equation}
This is easily seen to be well-defined (independent of the choice of
$g_0$). 
Pushing forward the measures $\nu_{x_0}$ and $\nu_Y$ via the map 
$\tau|_Y$ defines measures on $H\backslash G$,  
which we denote by the same letters.
Although $\nu_Y$ (as a measure on $\HG$) depends on the choice
of $Y$, the measure $\nu_{x_0}$ 
depends only on $x_0$. 
See 
Proposition \ref{prop: measure well defined} below for more details.

\begin{thm}
\label{thm: identifying the limit}
Assume the general setup holds, and conditions UC, D1, 
D2 are satisfied.
Let $x_0 \in \HG$
compactly supported continuous function $\varphi$ on 
$\HG$, 
\begin{equation}\label{eq_th_dual}
\frac{1}{\lambda(H_T)}
\int_{G_T} \varphi(x_0\, g)\, dm(g)
\mathop{\longrightarrow}_{T\rightarrow\infty}
\int_{H\backslash G} \varphi \,
d\nu_{x_0}. 
\end{equation}
\end{thm}
Note that the conclusion of Theorem \ref{thm: identifying the limit} may fail if condition D2 does
not hold, see \S \ref{section: examples}.2.

As a corollary we obtain an equidistribution result for certain
$\Gamma$-orbits. To state it, for $A \subset \HG$ and $x_0 \in \HG$, let
$$
N_T(A,x_0) \df \# \{\gamma \in \Gamma_T: x_0 \, \gamma \in A\}.
$$
\begin{cor}
\label{cor: limit for indicator sets}
Assume the general setup holds, and conditions S, UC, D1, D2 are satisfied.
Let $g_0 \in G, \, x_0 = \tau(g_0)$, and assume $\left( ** \right)$ holds.
Then 
for any $\varphi \in C_c(\HG)$ we have: 
\begin{equation}\label{eq_last}
\frac{1}{\lambda(H_T)}
\sum_{\gamma \in \Gamma_T} \varphi(x_0 \, \gamma) \mathop{\longrightarrow}_{T\rightarrow\infty}
\int_{\HG} \varphi \, d\nu_{x_0},
\end{equation}
and for any bounded $A \subset \HG$ with
$\nu_{x_0}(\partial 
A)=0$ we have:
\begin{equation}\label{eq: indicator sets}
\frac{N_T(A, x_0)}{\lambda(H_T)}
\mathop{\longrightarrow}_{T\rightarrow\infty} 
\nu_{x_0}(A)
\end{equation}
\end{cor}

\subsection{Distance functions}
In order to apply our general results we must show that
the conditions listed above hold in some
cases of interest. Showing this 
can be quite complicated, and requires specific methods for specific
cases. We verify the conditions under several hypotheses on $D$, $G$
and $H$. In all cases we consider in this paper, $G$ will be a Lie
subgroup and $H$ a Lie subgroup. In this setup condition S holds by a
standard application of the implicit function theorem, 
see
e.g. \cite[Thm. 3.58]{Warner}.
Note that if $G$ is not assumed to be a Lie group, condition S might not
hold (see \S \ref{section: examples}.1).
In this paper we will always assume furthermore that at least one of $H$ and $G$ is
semisimple and connected. 

\subsubsection{Linear norms}
Suppose $G$ is realized, via a linear representation $\Psi :
G \to \GL(d, \R)$, as a
closed subgroup of 
$\Mat_d(\R)$, and $\| \cdot \|$ is a norm on $\Mat_d(\R)$
(considered as a vector space). Assume that $\ker \Psi$ is compact. In
this setting we call 
$$D(g) \df \max \left\{ 1, \| \Psi(g) \| \right \}$$
a {\em matrix norm distance function}. 
Although $D$ depends on both $\Psi$ and $\| \cdot \|$, to simplify
notation we omit this dependence from the notation.

\subsubsection{Symmetric space distance functions}
\label{subsection: symmetric spaces}
 Suppose $G$ is a connected, semisimple Lie group, $K$ a maximal
compact subgroup of $G$, and $H$ a semisimple
subgroup. Thus $K \backslash G$ (resp., $G/K$) is the right
(resp. left) symmetric space of $G$. It is equipped with a natural
$G$-invariant Riemannian metric induced by the Cartan--Killing form on
the Lie algebra $\goth{g}$ of $G$ (see \cite[Ch. IV]{Helgason} or
\cite[Ch. 2]{Eberlein}). Let $P: G    
\to K \backslash G$ (resp., $P': G \to G/K$) be the natural projection. 
  Let $d$ (resp. $d'$) 
be the corresponding metric on $K \backslash G$ (resp., $G/K$). In
this setting we 
call 
$$D(g) \df \exp \left(d\left(P(e), P(g)\right)\right)$$ 
(respectively $\exp(d'(P'(e), P'(g)))$) 
a {\em symmetric space distance function.} Note that $D$ is proper
since $K$ is compact and $x \mapsto d(x_0,x)$ is proper on
$K\backslash G$.

\begin{definition}
Suppose the general setup holds. We say that $G, D$ are {\em standard}
if $G$ is semisimple and connected, and $D$ is either a matrix norm
distance function or a symmetric space distance function. 
We say that $G, H, D$ are {\em
standard} in either of the
following two cases:
\begin{itemize}
\item
$H$ is
semisimple and connected and $D$ is a matrix norm distance function. 
\item
Both $G$ and $H$ are semisimple and connected and $D$ is a symmetric
space distance 
function. 
\end{itemize}
\end{definition}

\begin{prop}
\label{prop: verifying axioms 1}
If $G, D$ are standard and $G$ is balanced then UC, I1 and I2 are
satisfied. 
If $G, H, D$ are standard then conditions UC, D1, and D2 are 
satisfied.

\end{prop}

We deduce the proposition from the following two results, which are of
independent interest, regarding the  asymptotics of
volumes of `balls' with respect to a standard distance function. 

\begin{thm}
\label{thm: volume asymptotics1}
Let $H$ be connected semisimple Lie group and let $\Psi: H \to \GL(V)$ be a representation
with compact kernel. Then  for any matrix norm distance function we have 
$$\lambda(H_T) \sim C (\log T)^{\ell} T^m,$$
where $\ell \in \Z_+$ and $m>0$
are explicitly computable constants depending
only on $\Psi$, and $C>0$ is an explicitly computable constant
depending continuously on the choice of norm on $\Mat_{d}(\R)$. 
\end{thm}

In \S \ref{section: verifying the axioms 1} we introduce a technical
condition, 
called condition {\bf G}, about $\Psi$ (see Definition \ref{def:
condition G}). The condition holds for `most' choices of $\Psi$ (cf. Remark
\ref{remark: condition generic}). We prove Theorem \ref{thm: volume
asymptotics}, which is a precise version of Theorem \ref{thm: volume
asymptotics1}, under the assumption
that condition {\bf G} holds, and assuming that 
$\Psi$ is 
irreducible. The proof of Theorem \ref{thm: volume 
asymptotics1} without
making these assumptions is a rather lengthy computation which will
appear in a separate paper.

In \S \ref{section: verifying the axioms 2} we treat the case of a
symmetric space distance function, and prove:

\begin{thm}
\label{thm: symmetric volume asymptotics 1}
Let $G$ be a connected semisimple Lie group, $H$ a connected semisimple
subgroup, and let $D: G \to [0, \infty)$  be a symmetric space
distance function. There are non-negative
constants $\ell$ and $m$, and for any
$g_1, \, g_2 \in G$,  a positive constant $C = C(g_1, g_2)$ such that 
\begin{equation}
\label{eq: formula involving C()}
\lambda\left(H_T[g_1,g_2]\right) \sim C (\log T)^{\ell}
T^{m}.
\end{equation}
The constants $\ell\geq 0, \, m>0, \, C(g_1, g_2)>0$ are given explicitly, 
$C(g_1, g_2)$ is continuous in $g_1, g_2$, and the
convergence in \equ{eq: formula involving C()} is uniform
for $g_1, g_2$ 
in compact subsets of $G$.
\end{thm}

\subsection{Balanced semisimple groups}
We will see below that {\em when $H$ is simple} and $G, H, D$ are
standard, condition $\left(**
\right)$ holds whenever $H\pi(g_0)$ is dense in $G/\Gamma$.  
However, in the case of a general semisimple group, we need to make an
additional assumption. 

Suppose that $H$ is a semisimple non-simple Lie
group with finite center, and thus a nontrivial almost direct product
of its simple 
factors. It may happen that some of
these factors do not contribute to the volume growth of balls in $H$,
and therefore should be ignored when computing the asymptotics of
$H$-orbits. To make this precise, we make the following definition.

Let $H = H_1 \cdots H_t$ be the decomposition of $H$ into 
an almost direct product of its  simple factors.
Fix measurable sections $\sigma_i: H \to H_i, \, i=1,
\ldots, t$; that is, $\sigma_i(h) \in H_i$ and $h=\sigma_1(h) \cdots
\sigma_t(h)$ for each $h \in H$. 

\begin{definition}
\label{def: balanced} 
We say that $H$ is {\em balanced} if for every $j \in \{1, \ldots,
t\}$, every $g_1, g_2 \in G$ and every compact $L \subset H_j$,
\begin{equation}
\label{eq: balanced}
\lim_{T \to \infty} \frac{\lambda\left(\{h \in H_T[g_1, g_2]:
\sigma_j(h) \in L 
\right\})}{\lambda(H_T[g_1, g_2])}  
=0.
\end{equation}
\end{definition}

Note that this definition  does not depend on the choice of the
sections because for any two sections $\sigma_i,
\, \sigma'_i$ we have $\sigma'_i(h) \in Z\, \sigma_i(h)$ for every $h \in
H$, where $Z$ is the (finite) center of $H$. Note also that if $H$ has compact factors then $H$ is not balanced by our
definition.

We remark that if $D$ is a matrix norm distance function,
corresponding to a 
representation $\Psi: G \to \GL_d(\R) \subset \Mat_{d}(\R)$ with
compact kernel, and a norm
$\| \cdot \|$ on 
$\Mat_d(\R)$, then the condition that $H$ is balanced depends on
$\Psi$ but not on $\| \cdot \|$ (see Proposition \ref{prop:
contribution to volume}). 

The importance of the assumption that $H$ is balanced lies in the
following:

\begin{thm}
\label{thm: using Ratner}
Suppose the general setup holds, $G, H, D$ are standard and $H$ is 
balanced. Suppose $g_0
\in G$ satisfies $\cl{H\pi(g_0)} = G/\Gamma$. Then $\left(**\right)$
holds. 
\end{thm}

Note that the conclusion of theorem 2.10 may fail if $H$ is not
balanced, see \S \ref{section: examples}.3.

\subsection{Equidistribution results}
Collecting the  results stated previously, we obtain the following
results concerning the distribution of lattice orbits:

\begin{cor}
\label{cor: putting together, induction}
Suppose the general setup holds, $G, D$ are standard, and $G$ is
balanced. Suppose there is a $\Gamma$-invariant probability measure
$\mu$ on $\HG$, and $x_0 \in \HG$ satisfies $\HG = \cl{x_0 \,
\Gamma}.$

Then for any $\varphi \in C_c(\HG),$
\begin{equation*}
\frac{1}{\# \Gamma_T} \sum_{\gamma
\in \Gamma_T} \varphi (x_0 \,  \gamma) \mathop{\longrightarrow}_{T \to
\infty} \int_{\HG} 
\varphi \, d\mu.
\end{equation*}
\end{cor}

\begin{cor}
\label{cor: putting together, duality}
Suppose the general setup holds, $G, H, D$ are standard, and $H$ is
balanced. Suppose
$x \in \HG$ with $\HG= \cl{x \,  \Gamma}$. Then for 
any $\varphi \in C_c(V),$
\begin{equation}
\label{eq: asymptotics like G}
S_{\varphi, x}(T) \sim \til S_{\varphi, x}(T)
\end{equation}
and 
\begin{equation}
\label{eq: precise asymptotics}
\frac{S_{\varphi, x}(T)}{\lambda(H_T)}
\mathop{\longrightarrow}_{T\rightarrow\infty} 
\int_{\HG} \varphi \, d\nu_{x}.
\end{equation}
Also, for any bounded $A \subset \HG$ with $\nu_{\HG} \left(\partial A
\right)=0$, 
$$\lim_{T \to  \infty} \, \frac{N_T(A,x)}{\lambda(H_T)} = \nu_{x}(A).
$$
\end{cor}


\section{Induction}
\label{section: induction}
Let $\Gamma$ act continuously on a locally compact Hausdorff space $X$
from the right. In this section
we recall the definition of the {\em induced action} of $G$, and
relate the asymptotics of a $\Gamma$-orbit with that of a
corresponding $G$-orbit.

On  $\til Y \df X \times G$, there is a left action of $G$ given by
$$
g_1\cdot (x,g)=(x, g_1g), 
$$
and a right action of $\Gamma$ given by
$$
(x,g)\cdot\gamma=(x \cdot \gamma, g\gamma). 
$$
It is clear that these actions commute, so that $G$ acts on the
quotient $Y = \til Y/\Gamma,$ and the quotient map $\pi_{G/\Gamma} : Y
\to G/\Gamma$ is a continuous $G$-map. Since $\Gamma$ is discrete in
$G$, the 
projection is a covering map, so that $Y$ inherits (locally) the
topological properties of $\til Y$. As a $G$-space, $Y$ is a fiber
bundle with base $G/\Gamma$ and fiber isomorphic to $X$. In
particular, if $\mu$ is a $\Gamma$-invariant probability measure on
$X$ then $d\mu \, dm'$ defines a finite $G$-invariant measure on $Y$.  

Assume that $\sigma:
G/\Gamma \to G$ is a Borel section. Then 
$$\gamma(g)=g^{-1} \, \sigma (\pi (g)) \in \Gamma$$
 and the map 
$$\pi_{X}=\pi^{\sigma}_X: Y \to X, \ \ \pi_X: y=[(x,g)] \mapsto x
\cdot \gamma(g)$$ 
is well-defined (does not depend on the representative $(x,g) \in
\til Y$ of
the $\Gamma$-orbit $[(x,g)] \in Y$). 

Note that when $X = \HG$ is as in \S \ref{subsection: inducing}, the
space $Y$ is isomorphic to $\HG \times G/\Gamma$, with the action
given by \equ{eq: G action, induction} where the isomorphism is defined by the map
$$
(x,g)\mapsto (xg^{-1},g),\quad x\in X,\;g\in G.
$$

The following reduces the study of asymptotic behavior of
$\Gamma$-orbits on $X$ to that of $G$-orbits on $Y$. It immediately
implies Proposition \ref{prop: induction, special}.

\begin{prop}
\label{prop: induction}
Suppose hypotheses I1 and UC hold. Let $x_0 \in X$, let $\til y_0
= (x_0,e ) \in \til Y$ and let $y_0 = [\til y_0] \in Y$. Suppose that
there is a measure 
$\mu$ on $X$ such that
for any $F \in C_c(Y)$,
\begin{equation}
\label{eq: hypothesis ***}
\frac{1}{m(G_T)} \int_{G_T} F(g^{-1} \, y_0)\, dm(g)
\mathop{\longrightarrow}_{T \to \infty} \int_Y F\, d\nu
,
\end{equation}
where $\nu =  d\mu \, dm'$.

Then for any $\varphi \in C_c(X),$
\begin{equation}
\label{eq: first consequence}
\frac{1}{m(G_T)} \sum_{\gamma \in \Gamma_T} \varphi (x_0 \cdot \gamma)
\mathop{\longrightarrow}_{T \to \infty} \int_X \varphi \, d\mu.
\end{equation}

If in addition I2 holds then 
\begin{equation}
\label{eq: second consequence}
\frac{1}{\# \Gamma_T} \sum_{\gamma
\in \Gamma_T} \varphi (x_0 \cdot \gamma) \mathop{\longrightarrow}_{T \to \infty} \int_X
\varphi \, d\mu.
\end{equation}
\end{prop}

\begin{proof}[Proof of Proposition \ref{prop: induction}]
There is no loss of generality in assuming that $\varphi \geq 0.$ 

Let $\vre>0.$ By condition I1, there are positive $\delta, \, T_0$
such that for all $T>T_0$, 
\begin{equation}
\label{eq: comparable balls}
\frac{m(G_{(1+\delta)T})}{m(G_T)} <1+\vre.
\end{equation}
By condition UC, there exists a symmetric neighborhood
$\mathcal{O}$ of the identity in $G$ such that for
every $T>0$, 
\begin{equation}
\label{eq_norm_simp_est}
 G_T \, \mathcal{O} \subset G_{(1+\delta)T}.
\end{equation}

Let $\chi \in C_c(G)$ be a non-negative function such that
$\supp \, \chi \subset \mathcal{O}$ and 
$\int_{G}\chi\, dm=1$. 
Define functions $\til F: \til{Y} \to \R, \, F: Y \to \R$ by
$$
\til F(\til y)= \til F(x,g)=\chi(g) \, \varphi(x), \ \ \
\ F(y) = \sum_{\gamma \in \Gamma} \til{F} (\til y \gamma),
$$
where $\til y \in \til Y, \, y=[\til y] \in Y.$ Since $\til F$ is
compactly supported, the sum in the 
definition of $F$ is actually finite, $F \in C_c(Y)$, and 
\[
\begin{split}
 \int_{Y} F \, d \nu 
& = \int_{G/\Gamma} \int_{X} F(y) \, d\mu(\pi_X(y)) \,
dm'(\pi_{G/\Gamma}(y))\\  
& =\sum_{\gamma \in \Gamma} \int_{\sigma(G/\Gamma)} \int_X \chi(g \gamma)
\varphi(x \cdot \gamma)\, d \mu(x) \, dm(g) \\ 
& = \int_{X} \varphi\, d\mu \, \left(\sum_{\gamma \in \Gamma}
\int_{\sigma(G/\Gamma)} \chi(g \gamma) \, dm(g) \right)\\
& = \int_{X} \varphi \, d\mu \, \int_G \chi \, dm = \int_{X} \varphi \, d\mu.
\end{split}
\]

Applying \equ{eq: hypothesis ***} we obtain:
\begin{equation}
\label{eq: limit of F}
\frac{1}{m(G_T)}\int_{G_T} F(g^{-1} \, y_0) \, d
m(g) \mathop{\longrightarrow}_{T \to \infty} \int_{X} \varphi \, d\mu.
\end{equation}

Now let 
\begin{equation}
\label{eq: def IT}
\begin{split}
I_T(x_0) & \df \sum_{\gamma \in \Gamma} \left(\int_{G_T} \chi(g^{-1}  \gamma)
\, dm(g) \right) \varphi(x_0 \cdot \gamma) \\
& = \sum_{\gamma \in \Gamma}\int_{G_T} \til F \left((x_0, g^{-1}
) \cdot \gamma \right )\, d m(g)\\
& = \int_{G_T} F(g^{-1} \, y_0)\, d m(g) .
\end{split}
\end{equation}
We  claim that
\begin{equation}\label{eq_I_N}
I_{T/(1+\delta)} (x_0)\le \sum_{\gamma \in\Gamma_T}\varphi(x_0 \cdot
\gamma)\le I_{(1+\delta)T}(x_0).  
\end{equation}
Let $g \in G \sm G_{(1+\delta)T}$. Then 
by \equ{eq_norm_simp_est}, $g \notin G_T \, \mathcal{O}$, and
hence $g^{-1} \gamma \notin \mathcal{O}$ for all $\gamma \in
\Gamma_T$. This implies that $\chi(g^{-1} \gamma)=0$, that is,

$$\gamma \in \Gamma_T \ \ \Longrightarrow \ \
\int_{G_{(1+\delta)T}}\chi (g^{-1}\gamma)\, dm(g) =  \int_{G}\chi
(g^{-1}\gamma)\, dm(g) = 1.$$
The right hand inequality in 
(\ref{eq_I_N}) follows. 

Now if $\gamma \in \Gamma \sm \Gamma_T$ then for all $g \in
G_{T/(1+\delta)}$ we have $g^{-1} \gamma \notin \mathcal{O}$, hence
$\chi(g^{-1} \gamma )=0$ and so 
$$\gamma \notin \Gamma_T \ \ \Longrightarrow \ \ 
\int_{G_{T/(1+\delta)}} \chi(g^{-1}  \gamma)
\, dm(g)=0
.$$ 
This implies the second inequality in \equ{eq_I_N}.

\medskip

We obtain, for all large enough $T$:
\[
\begin{split}
 \sum_{\gamma \in\Gamma_T}\varphi(x_0 \cdot
\gamma) &  \stackrel{\equ{eq_I_N}}{\leq} I_{T(1+\delta)}(x_0) \\
&\stackrel{\equ{eq: limit of F}}{\leq} (1+\vre) m(G_{(1+\delta)T})
\int_X \varphi \, d\nu \\
& \stackrel{\equ{eq: comparable balls}}{\leq} (1+\vre)^2 m(G_T)
\int_X \varphi \, d\nu.
\end{split}
\]

Since $\vre>0$ was arbitrary, 
$$\limsup_{T \to \infty} \frac{1}{m(G_T)} \sum_{\gamma
\in\Gamma_T}\varphi(x_0 \cdot \gamma)\leq \int_X \varphi \, d\nu.$$
The opposite inequality for $\liminf$ is similarly established. This
proves \equ{eq: first consequence} and 
\equ{eq: second consequence} immediately follows in light of condition
I2.
\end{proof}

\section{Duality}

In this section we prove Theorem \ref{thm: duality}.

We first record the following useful consequence of condition UC.

\begin{prop}
\label{prop: A2 alternative form}
Suppose the general setup holds and condition UC is satisfied. Then
for any bounded $B \subset G$ and any $\eta>0$ there is a neighborhood 
$\mathcal{O}$ of identity in $G$ such that for all $g \in G, \, h \in
\mathcal{O}, \, b \in B,$ we have 
\begin{equation}
\label{eq: A2}
\left | \frac{D(ghb)}{D(gb)} - 1 \right | < \eta.
\end{equation}

\end{prop}

\begin{proof}
Let $\vre = \eta$ and let $\mathcal{U}$ be as in the formulation of
condition UC. Since $B$ is bounded, there is a small enough
neighborhood $\mathcal{O}$ of identity in $G$ such that for all $h\in
\mathcal{O}, \, b \in B$ we have $b^{-1}hb \in \mathcal{U} \cap
\mathcal{U}^{-1}$. This 
implies that for $g\in G, \, h \in \mathcal{O}, \, b \in B,$ 
$$
D(ghb) = D(gb \, b^{-1}hb) <(1+\eta)D(gb)$$
and
$$D(gb) = D(ghb \, b^{-1}h^{-1}b) \leq (1+\eta) D(ghb).
$$
Putting these together we obtain
$$\frac{D(gb)}{1+\eta} < D(ghb) <(1+\eta)D(gb),
$$
and (\ref{eq: A2}) follows.
\end{proof}

We now make some reductions:
\begin{claim}
\label{claim: 1}
There is no loss of generality in assuming:
\begin{enumerate}
\item[(a)]
 $g_0 = e$.
\item[(b)]
There is an open set $\mathcal{U} \subset \HG$ and a Borel section
 $\sigma: \HG  \to G$ such that $\sigma|_{\mathcal{U}}$ is continuous
 and $\supp \, \varphi \subset  \mathcal{U}.$ 
\end{enumerate}
\end{claim}

In order to reduce to (a), let $H' \df g_0^{-1}Hg_0, \,
\lambda'$ a 
left haar measure on $H', \, \tau': G \to H'
\backslash G$ the natural quotient map and $\varphi'(\tau'(g)) \df
\varphi(\tau(g_0g))$ and note that the hypotheses and conclusion of
the theorem for $H'$, $\lambda'$, $\tau'$, $\varphi', \, e$ are 
equivalent to those for $H$, $\lambda$, $\tau$, $\varphi, \, g_0$. 

Using condition S, for each $x \in \supp \, \varphi$ there is a
section $\sigma_x: \HG \to G$ which is continuous in a neighborhood
$\mathcal{U}_x$ of $x$. Using the compactness of $\supp \, \varphi$
we may choose a finite subset of $\{\mathcal{U}_x: x \in \supp \, 
\varphi\}$ covering $\supp \, \varphi$, and using a standard partition
of unity argument for this 
cover we may assume that for some $x \in \HG$, $\supp \, \varphi
\subset \mathcal{U}_x$ and take $\sigma = \sigma_x$. This reduces our
problem to the case (b).
The claim is proved.

\medskip

 Let 
$$B \df \sigma(\supp \, \varphi),$$
a compact subset of $G$, in view of Claim \ref{claim: 1}.
Also let 
$$h_g \df \til \sigma (g)g^{-1} \in H, \ \ \ \mathrm{where} \ \ \til
\sigma \df \sigma \circ \tau: G \to G.$$


Given $c>1$ we will find
$T_0$ such that for all $T>T_0$, 
\begin{equation}
\label{eq: bound from above}
\sum_{\gamma \in \Gamma_T} \varphi(\tau(\gamma)) <
c\int_{G_T} \varphi(g) \, dm(g). 
\end{equation}

Let 
\begin{equation}
\label{eq: size of vre}
0<\vre < \sqrt{c}-1,
\end{equation}
and let $\delta>0$ be as in hypothesis D1. Let
$\eta>0$ be small enough so that 
\begin{equation}
\label{eq: eta}
(1+\eta)^4 < 1+\delta. 
\end{equation}
Using hypothesis UC and Proposition \ref{prop: A2 alternative form}, 
let $\mathcal{O}$ be a small enough 
neighborhood of the 
identity in $H$ so that for all $g \in G, \, h \in \mathcal{O}, \, b
\in B,$ 
$$D(ghb)<(1+\eta)D(gb).$$ 

Let $\psi: H \to \R$ be a non-negative continuous function such that
$$\supp \, \psi \subset \mathcal{O}\ \ \mathrm{and} \ \ \int_{H}
\psi(h)\, d\lambda(h)=1.$$ 

Define a function $f: G \to \R$ by 
$$f(g) = \psi(h_g)\varphi(\tau(g)).$$

Note that $f$ is non-negative, continuous, of compact support, and
\begin{equation}
\label{eq: using spare}
\tau(g) \notin \supp \, \varphi \ \ \Longrightarrow \ \ \forall h \in
H, \, f(h^{-1}g)=0.
\end{equation}

We now show:

\medskip

\begin{claim}
\label{claim: first}
For any $T>0$,
\begin{equation}
\label{eq: A}
\sum_{\gamma \in \Gamma_T} \varphi(\tau(\gamma))  \leq
\sum_{\gamma \in \Gamma}
\int_{H_{(1+\eta)T}[\til \sigma(\gamma)]} f(h^{-1}\gamma) \, d
\lambda(h)
\end{equation}
and 
\begin{equation}
\label{eq: B}
\int_{G_T} \varphi(\tau(g))\, dm(g)  \geq \int_G
\int_{H_{T/(1+\eta)}[\til \sigma (g)]}
f(h^{-1}g) d \lambda(h)\, dm(g). 
\end{equation}
\end{claim}

Let $g \in G_T$, that is, $D(h_g^{-1}\til \sigma (g))=D(g)<T$. Suppose
$\tau(g) \in \supp \, \varphi$. 
Then for $h \in \mathcal{O}$ one has 
$$D(h_g^{-1}h\til \sigma (g)) <(1+\eta)T,
$$
that is, $h_g^{-1}h \in H_{(1+\eta)T}[\til \sigma (g)].$
Since $\supp \, \psi \subset \mathcal{O},$ this shows that 
\begin{equation*}
\label{eq: supp contained}
\supp \, \psi \subset 
h_{g}H_{(1+\eta)T}[\til \sigma(g)].
\end{equation*}
So, using left-invariance of $\lambda$ and the fact that $h_gh_0 =
h_{h_0^{-1}g}$:
\[
\begin{split}
\varphi(\tau(g)) & = \varphi(\tau(g))
\int_{h_{g}H_{(1+\eta)T}[\til \sigma(g)]} \psi(h)\,
d\lambda(h)\\ 
&= \varphi(\tau(g)) \int_{H_{(1+\eta)T}[\til \sigma(g)]} \psi(h_{g}h) \,
d\lambda(h)\\
&= \int_{H_{(1+\eta)T}[\til \sigma(g)]} \psi(h_{h^{-1}g})
\varphi(\tau(h^{-1}g)) \, d\lambda(h)\\
& =  \int_{H_{(1+\eta)T}[\til \sigma(g)]} f(h^{-1}g)\, d
\lambda(h).
\end{split}
\]

Specializing to $g = \gamma \in \Gamma_T$ one obtains (\ref{eq: A}). 

\medskip

Let $\til{G} \df \tau^{-1}(\supp \, \varphi)$.
Similarly, it is easy to check that for $g \in \til{G} \sm G_T,$
$$\supp \, \psi \cap
h_gH_{T/(1+\eta)}[\til \sigma (g)] = \varnothing.$$

Therefore, arguing as in the previous computation:
\[
\begin{split}
& \int_G \int_{H_{T/(1+\eta)}[\til \sigma(g)]} f(h^{-1}g)d \lambda(h)
\, dm(g) \\
& \stackrel{(\ref{eq: using spare})}{=} \int_{\til{G}}
\int_{H_{T/(1+\eta)}[\til \sigma(g)]} f(h^{-1}g)d \lambda(h) \, dm(g) \\
& \leq \left[\int_{G_T}+ \int_{\til{G} \sm G_T} \right] \left (\varphi(\tau(g))
\int_{h_gH_{T/(1+\eta)}[\til \sigma(g)]} \psi(h)d\lambda(h) \right) dm(g)\\ 
& = \int_{G_T} \varphi(\tau(g))
\int_{h_gH_{T/(1+\eta)}[\til \sigma(g)]} \psi(h)d\lambda(h)\, dm(g) \\
& \leq \int_{G_T} \varphi(\tau(g))\, dm(g). 
 \end{split}
\]

This proves (\ref{eq: B}).

\medskip

We now claim:

\medskip

\begin{claim}
\label{claim: second}
There are $b_1, \ldots, b_N \in B$ and a symmetric neighborhood
$\mathcal{O}$ of the identity in $G$, such that for $x \in
\mathcal{O}b_i, \, i=1, \ldots, N$,
\begin{equation*}
\label{eq: in claim}
H_{(1+\eta)T}[x] \subset H_{(1+\eta)^2T}[b_i]
\end{equation*}
and such that 
\begin{equation}
\label{eq: cover in claim}
B \subset \bigcup_{i=1}^N \mathcal{O}b_i.
\end{equation}
\end{claim}

To see this, using hypothesis UC and Proposition \ref{prop: A2
alternative form} we let $\mathcal{O}$ be a small
enough symmetric neighborhood of identity in $G$ such that for all $g
\in H, \, h \in \mathcal{O}$ and $b \in B$, we have 
$(1-\eta/2)D(gb)<D(ghb)$. This implies that $H_{(1+\eta)T}[x] \subset
H_{(1+\eta)^2T}[b]$ for all $x \in \mathcal{O}b.$ Taking a finite
subcover $\{\mathcal{O}b_i: 
i=1, \ldots, N\}$ of the cover $\{\mathcal{O}b : b \in B\}$ proves the
claim.


\medskip

Using a partition of unity subordinate to the cover (\ref{eq: cover in
claim}), there is
no loss of generality in assuming that for some $b \in B$,
\begin{equation}
\label{eq: for s}
\supp \, \varphi \subset \tau(\mathcal{O}b).
\end{equation}

Thus for some $b \in B$, and
for all $\gamma \in \Gamma,$

\begin{equation}
\label{eq: next claim}
\int_{H_{(1+\eta)T}[\til \sigma(\gamma)]} f(h^{-1}\gamma )\, d
\lambda(h) \leq \int_{H_{(1+\eta)^2T}[b]}
f(h^{-1}\gamma )\, d 
\lambda(h).
\end{equation}

Now defining 
$$F (\pi(g))  \df \sum_{\gamma \in \Gamma} f(g \gamma) $$
(actually a finite sum for each $g$) one obtains via the monotone
convergence theorem:


\begin{equation}
\label{eq: nextest}
\begin{split}
& \sum_{\gamma \in \Gamma} \int_{H_{(1+\eta)T}[\til \sigma (\gamma)]}
f(h^{-1}\gamma) \, d \lambda(h) \\
& \stackrel{(\ref{eq: next claim})}{\leq} \, \sum_{\gamma \in \Gamma
} \int_{H_{(1+\eta)^2T}[b]} f(h^{-1}
\gamma)\, d\lambda(h) \\ 
& = \int_{H_{(1+\eta)^2T}[b]} F(h^{-1}\pi(e))\, d\lambda(h).
\end{split}
\end{equation}

Using $\left( ** \right)$, there is $T_1$ such that for all
$T>T_1$: 

$$
\left |\frac{\int_{H_{(1+\eta)^2T}
[b]} F(h^{-1}\pi(e))\, d\lambda(h)}{\lambda \left(H_{(1+\eta)^2T} [b]
\right)}- \int_{G/\Gamma} F(x)\, dm'(x)\right |<\vre,
$$

hence 

\begin{equation}
\label{eq: follows from hypothesis}
\int_{H_{(1+\eta)^2T} [b]} F(h^{-1}\pi(e))\, d\lambda(h) \leq (1+\vre)
\lambda \left(H_{(1+\eta)^2T} [b] \right) \int_{G/\Gamma}
F(x)\, dm'(x). 
\end{equation}

Reversing the argument with $G$ in place of $\Gamma$ yields:
\begin{equation}
\label{eq: converse for G}
\begin{split}
& \lambda \left(H_{(1+\eta)^2T} [b] \right) \int_{G/\Gamma}
F(x)\, dm'(x) \\
& \stackrel{(\ref{eq: A1}), \, (\ref{eq: eta})}{\leq} (1+\vre) \lambda 
\left(H_{T/(1+\eta)^2} [b] \right) \int_{G} f(g)\, dm(g) \\  
& = (1+\vre) \int_{H_{T/(1+\eta)^2} [b]} \int_G f(h^{-1} g)\, dm(g) \,
d\lambda(h) \\ 
& = (1+\vre) \int_{G} \int_{H_{T/(1+\eta)^2} [b]} f(h^{-1}g)\,
d\lambda(h) \, dm(g) \\
& \stackrel{(\ref{eq: using spare}),\,(\ref{eq: for s})}{\leq}
(1+\vre) \int_{G} \int_{H_{T/(1+\eta)} [\til \sigma (g)]} f(h^{-1}g)\,
d\lambda(h) \, dm(g) \\
  & \stackrel{(\ref{eq: B})}{\leq} (1+\vre)
\int_{G_{T}} \varphi(\tau(g))\, dm(g).
\end{split}
\end{equation}

Now putting together (\ref{eq: A}), (\ref{eq: nextest}), (\ref{eq:
follows from hypothesis}), (\ref{eq: converse for G}) and
(\ref{eq: size of vre}) proves  
(\ref{eq: bound from above}). 

\medskip

The proof of the second inequality 
\begin{equation}
\label{eq: bound from below}
\sum_{\gamma \in \Gamma_T} \varphi(\tau(\gamma)) >
\frac1c\int_{G_T} \varphi(\tau(g)) dm(g)
\end{equation}
is similar and is left to the reader. The two inequalities (\ref{eq:
bound from above}) and (\ref{eq: bound from below}), along with the
assumption that $g_0=e$, imply
(\ref{eq: to demonstrate}). \qed

\section{The limiting density}
\label{section: describing the
limit}
Let the notation be as in \S \ref{subsection 2.3}. In particular assume
condition D2. In this section we
discuss the measures $\{\nu_{x} : x \in \HG\}$ and prove Theorem
\ref{thm: identifying the limit}. We first introduce some additional
notation.

Let $N_G(H)$ denote the normalizer of $H$ in $G$. There is a
homomorphism 
$\Delta: N_G(H) \to \R_+$ such that for any Borel subset $A$ of $H$
and any $n \in N_G(H)$, 
\begin{equation}
\label{eq: defn modular}
\lambda(nAn^{-1}) = \Delta(n) \, \lambda(A).
\end{equation}
Note that $\Delta|_H$ is the modular function of $H$.

\begin{prop}
\label{prop: measure well defined}
Let $g_0 \in G, \, x_0 = \tau(g_0)$. Suppose that condition D2 is
satisfied. Then $\alpha(g_1, g_2)$
and the measure 
$\nu_{x_0}$ have the following properties:
\begin{itemize}
\item[(i)]
For $h \in H, \, n \in N_G(H)$ and $g_1, g_2 \in G$, 
$$\alpha(nhg_1, ng_2) = \Delta(n) \, \alpha(g_1, g_2).$$
In particular, for $h_1, h_2 \in H, \, g_1, g_2 \in G, $
$$\alpha(h_1g_1, h_2g_2) = \Delta(h_2)\, \alpha(g_1, g_2)
.$$

\item[(ii)]
$\nu_{x_0}$ is well-defined (does not depend on $g_0$) and does not
depend on the section $Y$.
\item[(iii)]
If $G, \, H$ are Lie groups then $\nu_{x_0}$ is
absolutely continuous with respect to the smooth measure class 
on $\HG$. 
\item[(iv)]
If $H$ is unimodular then 
$$\alpha_{x_0}: \HG \to \R_+, \ \ \ \alpha_{x_0}(\tau(g)) =
\alpha(g_0, 
g)$$ 
is well-defined and 
$d \nu_{x_0}(x) = \alpha_{x_0} (x) \,
d\nu_{\HG}(x)$, where $\nu_{\HG}$ is the $G$-invariant measure on $\HG$. 

\end{itemize}
\end{prop}

\begin{proof}
For $h \in H, \,
n \in N_G(H), \, g_1, g_2 \in G$ we have
\[
\begin{split}
H_T[(nhg_1)^{-1}, ng_2] & = \{h' \in H: D(g_1^{-1} (nh)^{-1} h' n g_2)<T \} \\
& = \{h' \in H: D(g_1^{-1} (nh)^{-1} n \, n^{-1} h' n \,
g_2)<T \} \\ 
& = n\, \{h'' \in H: D(g_1^{-1} h^{-1} h'' g_2)<T \} \,
n^{-1} \\ 
& = n h \, \{h''' \in H: D(g_1^{-1} h''' g_2)<T \}\, n^{-1} \\
& = n h \, H_T[g_1^{-1}, g_2] \, n^{-1},
\end{split}
\]
hence, using (\ref{eq: defn modular}) and the left-invariance of
$\lambda$, 
$$\lambda(H_T[(nhg_1)^{-1}, ng_2]) = \Delta(n) \, \lambda(H_T[g_1,
g_2])$$ 
and (i) follows.

It follows from (i) that $\alpha(hg_0, y)= \alpha(g_0, y)$ for all
$h \in H$, so $\nu_{x_0}$ is well-defined. 
Suppose $Y, \, Z$ are two different sections for $\tau$. There are
maps $\bar{h}:  Z \to H, \, \bar{y}: Z \to Y$ defined by the formula
$z=\bar{h}(z) \bar{y}(z)$ for all $z \in Z$. Let $\bar{y}_{*} \nu_Z$
denote the measure on $Y$ obtained by pushing forward $\nu_Z$ via
$\bar{y}$. Using (\ref{eq: measure on 
fibers}) we have for each $h \in H,$
\[
\begin{split}
 dm(hz) & = dm(h\bar{h}(z)\bar{y}(z)) = d\lambda(h\bar{h}(z)) \,
d\nu_Y(\bar{y}(z)) \\
& = \Delta(\bar{h}(z)) \, d\lambda(h) \, d\bar{y}_*\nu_Z(y).
\end{split}
\]

By the uniqueness of the decomposition (\ref{eq: measure on 
fibers}) we have $\Delta(\bar{h}(z)) \, d\bar{y}_*\nu_Z(y) =
d\nu_Y(y).$  In particular, setting $\bar{y}(z)=y, \, \bar{h}(z)= h =
zy^{-1}$ we have  
\[
\begin{split}
\alpha(g_0,y) d\nu_Y(y) & = \alpha(g_0, y z^{-1} \, z)
\Delta(zy^{-1}) \, d\bar{y}_* \nu_Z(y) \\
& = \Delta (yz^{-1}) \, \alpha(g_0, z) \Delta(zy^{-1})
d\nu_Z(z)\\
& = \alpha(g_0,z) d\nu_Z(z).
\end{split}
\]
proving (ii).

If $H$ is a Lie subgroup of a Lie group $G$ then there is a manifold
structure on $\HG$ which is obtained locally by finding, for each $x_0
\in \HG$, a submanifold $\mathcal{U} \subset G$ such that $x_0 \in
\tau (\mathcal{U})$ and $\tau|_{\mathcal{U}}$ is one to one. Now,
using (ii), we can assume that the lift $Y$ contains $\mathcal{U}$. It
is clear from (\ref{eq: measure on fibers}) and $G$-invariance of $m$
that $\nu_Y |_{\mathcal{U}}$ is smooth 
and (iii) follows.

The last assertion follows immediately from the preceding ones.
\end{proof}

\medskip

\begin{proof}[Proof of Theorem \ref{thm: identifying the limit}.]
By separating $\varphi$ into its negative and positive parts, we may
assume without loss of generality that $\varphi\ge 0$. Assume also
that $\varphi\ne 0$ (otherwise there is nothing to prove). 
Let $Y=\sigma (\HG)$ where $\sigma: \HG \to G$ is a Borel section, and let
$B = \sigma(\supp \, \varphi).$ 
We may assume that $B$ is bounded, e.g. by applying condition S and a
partition of unity argument. Let $\nu_Y$ be the measure on $Y$ satisfying (\ref{eq: measure on fibers}).
In view of Proposition \ref{prop:
measure well defined}, there is no loss of generality in defining
$\nu_{x_0}$ by (\ref{eq_nu_g}). We have:
\begin{equation}
\label{eq: newer one}
\begin{split}
\int_{G_T}\varphi(\tau(g_0)g) \, dm(g) & = \int_{\{g:D(g_0^{-1}g) < T\}}
\varphi(\tau(g)) \, dm(g)\\    
&=\int_Y \int_{\{h:D(g_0^{-1}hy)<T\}}\varphi(\tau(y)) \, d\lambda(h)
\, d\nu_Y(y)\\
&= \int_Y\varphi(\tau(y))\lambda(H_T[g_0^{-1},y]) \, d\nu_Y(y).
\end{split}
\end{equation}

Let $T_0$ be as in condition
D1. We claim that there is $C$ such that for all $T>T_0$ and all $y
\in B$, 
$$
\frac{\lambda(H_T[g_0^{-1},y])}{\lambda(H_T)}<C.$$

Let $C_1>1$. By Proposition \ref{prop: A2 alternative form}, there exists a neighborhood $\mathcal{O}$ of identity in $G$
such that for all $h\in\mathcal{O}$ and $b\in B$,
$D(g_0^{-1}hb)< C_1\, D(g_0^{-1}b)$.
Choose a finite cover $B\subset \cup_{i=1}^N \mathcal{O}b_i$ for some $b_i\in B$. Then
for every $y\in B$,
$$
H_T[g_0^{-1},y]\subset \bigcup_{i=1}^N H_{C_1T}[g_0^{-1},b_i]
$$
and
$$
\lambda(H_T[g_0^{-1},y])\le N\max_i \lambda(H_{C_1T}[g_0^{-1},b_i]).
$$
Using condition D1 a finite number of
times (depending on $C_1$) we find a constant
$C_2$ such that $\lambda(H_{C_1T}[g_0^{-1},b_i]) \leq C_2
\lambda(H_T[g_0^{-1},b_i])$. Since the limit defining $\alpha(g_0, b_i)$
exists, there is a constant $C_3>0$ 
such that 
$$\frac{\lambda(H_T[g_0^{-1},b_i])}{\lambda(H_T)} < C_3, \ \ \mathrm{for
\ all \ } T>T_0 \ \mathrm{and} \ i=1,\ldots, N.$$
Therefore 
$$
\frac{\lambda(H_T[g_0^{-1},y])}{\lambda(H_T)} 
< N\, C_2\, C_3 \df C.$$

Hence, by (\ref{eq_dual_as}), \equ{eq: newer one}, and Lebesgue's
dominated convergence theorem,
\begin{eqnarray*}
\lim_{T\to\infty} \frac{1}{\lambda(H_T)} \int_{G_T}
\varphi(\tau(g_0)g) \, dm(g) &=&
\int_B \varphi(\tau(y))\, \lim_{T\to \infty}
\frac{\lambda(H_T[g_0^{-1},y])}{\lambda(H_T)} \,  
d\nu_Y(y)
\\
& = & \int_{B}\varphi(\tau(y))\alpha (g_0,y) \, d\nu_Y(y)\\
&=&\int_Y \varphi \, d\nu_{x_0}.
\end{eqnarray*}
\end{proof}

\medskip

\begin{proof}[Proof of Corollary \ref{cor: limit for indicator sets}.]
Statement (\ref{eq_last}) is immediate from (\ref{eq: to demonstrate}) and (\ref{eq_th_dual}), and (\ref{eq: indicator sets}) follows by a
standard argument for approximating $1_A$ from above and below by
continuous functions, which we omit. 
\end{proof}

\section{Basics}
\label{section: preliminaries}
Our next goal will be to verify the conditions we have listed. Most of
the verification is reduced to the computation of asymptotics of
volumes of certain sets. We make the reduction in this section, and
also list some simple relations between our conditions.

We first show that I2 follows from I1, UC, and a very special case of
$\left(**\right)$. Namely, we have:  

\begin{prop}
\label{verifying I2}
Suppose that I1 and UC hold, and for any $F\in C_c(G/\Gamma)$,
\begin{equation}\label{eq_add1}
\frac{1}{m(G_T)}\int_{G_T} F(g^{-1}\pi(e)) \, dm(g)\mathop{\longrightarrow}_{T
\to \infty} \int_{G/\Gamma}
F \, dm'.
\end{equation}
Then I2 holds, that is, $\#\Gamma_T\sim m(G_T)$ as $T\to\infty$.
\end{prop}

\begin{proof}
Let $\vre>0$. By UC, there exists a symmetric neighborhood
$\mathcal{O}$ of $e$ 
in $G$ such that 
\begin{equation}\label{eq_add0}
 G_T\mathcal{O}\subset G_{(1+\vre)T}
\end{equation}
for every $T>0$. Let $f\in C_c(G)$ be such that
\begin{equation}\label{eq_add3}
f\ge 0,\quad \hbox{supp}(f)\subset \mathcal{O}, \quad \int_G fdm=1,
\end{equation}
and let $F \in C_c(G/\Gamma)$ be defined by
$F(\pi(g))=\sum_{\gamma\in\Gamma} f(g\gamma).$
We have
\begin{align*}
\int_{G_T} F(\pi(g^{-1}) )\, dm(g)&= \sum_{\gamma\in\Gamma} \int_{G_T}
f(g^{-1}\gamma)\, dm(g)  \\
& \stackrel{(\ref{eq_add0}),(\ref{eq_add3})}{=} \sum_{\gamma\in
\Gamma_{(1+\vre)T}} \int_{G_T} f(g^{-1}\gamma) \, dm(g)\\ 
&= \sum_{\gamma\in \Gamma_{(1+\vre)T}} \int_{G_T^{-1}\gamma}
f \, dm\stackrel{(\ref{eq_add3})}{\le} \#\Gamma_{(1+\vre)T}, 
\end{align*}
and
\begin{align*}
\int_{G_T} F(\pi(g^{-1}) )\, dm(g)&=  \sum_{\gamma\in\Gamma} \int_{G_T} f(g^{-1}\gamma)\, dm(g) \\
&\stackrel{(\ref{eq_add0}),(\ref{eq_add3})}{\ge} \sum_{\gamma\in
\Gamma_{T/(1+\vre)}} \int_{G_T} f \, dm =\#\Gamma_{T/(1+\vre)}.
\end{align*}
This implies that 
\begin{eqnarray*}
\limsup_{T\to\infty}\, \frac{\#\Gamma_T}{m(G_T)} \le \limsup_{T\to\infty}
\frac{1}{m(G_T)}\int_{G_{(1+\vre)T}} F(\pi(g^{-1})) \, dm(g)\\ 
\stackrel{(\ref{eq_add1})}{\le} \limsup_{T\to\infty}  \frac{m(G_{(1+\vre)T})}{m(G_T)}
\int_{G/\Gamma} F\, dm'=\limsup_{T\to\infty}
\frac{m(G_{(1+\vre)T})}{m(G_T)}
\end{eqnarray*}
for every $\vre>0$. Thus, by I1,
$$
\limsup_{T\to\infty} \frac{\#\Gamma_T}{m(G_T)}\le 1.
$$
The lower estimate for $\liminf$ is proved similarly.
\end{proof}

\begin{prop}
\label{prop: matrix distance functions have UC}
If $D$ is a matrix norm distance function then $D$ satisfies condition
UC.
\end{prop}

\begin{proof}
The set 
$$\Mat_d^1(\R) \df \{z \in \Mat_d(\R) : \, \|z \|=1\}$$
is compact. Hence, given $\vre>0$, there is a small enough
neighborhood $\mathcal{V}$ of identity in $\Mat_{d}(\R)$ such that for
all $x \in \mathcal{V}$ and all $z \in \Mat_d^1(\R)$,
$$ \|zx\| < (1+\vre).$$ 
Let $\mathcal{U} = \Psi^{-1}(\mathcal{V})$. There is a bounded set $B
\subset G$ such that for all $g \in G \sm B$ and all $g' \in
\mathcal{U}$ we have $\|\Psi(g)\| >1$ and $\|\Psi(gg')\|>1$ so that 
$$D(g)=\|\Psi(g)\| \ \ \ \mathrm{and} \ \ \ D(gg') = \|\Psi(gg')\|.$$
 For any $g \in G \sm B$,
$z=\frac{\Psi(g)}{\|\Psi(g)\|} \in \Mat_d^1(\R)$, hence for any $g'
\in \mathcal{U} $ we have
$$\frac{D(gg')}{D(g)} = \frac{\|\Psi(gg')\|}{\|\Psi(g)\|} = \|z \Psi(g')
\| < 1+\vre.
$$
By making $\mathcal{U}$ smaller if necessary we can also ensure that 
$$\frac{D(gg')}{D(g)} < 1 +\vre
$$
for all $g \in B$ and $g' \in \mathcal{U}$. This implies UC.
\end{proof}

\begin{prop}
\label{prop: right invariant metric have UC}
Suppose $(X,d)$ is a metric space equipped with a right  
(respectively, left)
action of $G$ which is continuous 
and isometric. Then for any $x_0 \in X$, the function
$$D(g) = \exp( d(x_0, x_0g) )\mathrm{ \ \ \ \
  (resp., \, } \exp(d(x_0, gx_0)) \mathrm{)}
$$
satisfies UC.

\end{prop}

\begin{proof}
Suppose the action is on the right. Given $\vre>0$, let $$\mathcal{U} \df
\{g \in G: d(x_0, x_0g) < \log (1+\vre)\}.$$ Then for every $g\in G$ and $u\in\mathcal{U}$,
\begin{align*}
d(x_0,x_0gu)-d(x_0,x_0g)=d(x_0u^{-1},x_0g)-d(x_0,x_0g)\\
\le d(x_0u^{-1},x_0)< \log (1+\vre).
\end{align*}
This implies UC. The proof for a left action is similar.
\end{proof}

\ignore{
\begin{prop}
\label{prop: right invariant metric have UC}
Suppose $d(\cdot, \cdot)$ is a continuous proper right-invariant or
left-invariant 
metric on $G$. Then $D(g) = \max \left \{ 1, d(g,e) \right \}$ satisfies UC.

\end{prop}

\begin{proof}
Suppose first that $d$ is right-invariant. The function $D: G \to [1,
\infty)$ is clearly continuous and proper. Given $\vre$, let 
$K \df \{g \in G: d(g,e) \leq 1\}$, a compact set. Since $D$ and $d$
are continuous, there is a small enough symmetric neighborhood
$\mathcal{U}$ of $e$ in $G$
such that 
\begin{equation}
\label{eq: first}
\forall g_0 \in \mathcal{U}, \ \ d(e,g_0)<\vre
\end{equation}
and 
\begin{equation}
\label{eq: second}
\forall g \in K\mathcal{U}, \, g_0 \in \mathcal{U}, \ \ \left|
\frac{D(gg_0)}{D(g)} -1 \right| < \vre.
\end{equation}

Thus \equ{eq: second} implies  \equ{eq: rt uniform cont} for $g \in K
\mathcal{U}$. In case $g \notin K\mathcal{U}$ we have 

\[
\begin{split}
D(gg_0) & = d(gg_0,e) =d(g, g_0^{-1}) \\
& \leq d(g,e)+d(e, g_0^{-1}) \stackrel{\equ{eq: first}}{\leq}
d(g,e)+\vre =d(g,e)\left(1+\frac{\vre}{d(g,e)}\right) \\
& \leq (1+\vre)d(g,e) \leq
(1+\vre)D(g).
\end{split}
\]
This proves \equ{eq: rt uniform cont}.

Now if $d$ is left-invariant then $\til d (x,y) \df d(x^{-1}, y^{-1})$
is right-invariant and by left-invariance, 
$$d(g,e)=d(e,g^{-1}) = \til d(e,g) = \til 
d (g,e),$$
so this case follows from the previous one. 
\end{proof}
}

\begin{prop}
\label{prop: D1 D2 from asymptotics}
Suppose $F: [0, \infty) \to [0, \infty)$ and $C: G \times G \to (0,
\infty)$ are 
continuous functions such that for all $\delta>0$,
$$\frac{F\left((1+\delta)T\right)}{F\left(T \right)}
\mathop{\longrightarrow}_{T \to \infty} 1.
$$
If the general setup holds and $m(G_T) \sim F(T)$ then I1 holds.

If the general setup holds and 
$$ \lambda \left(H_T[g_1, g_2] \right) \sim C(g_1, g_2) \, F(T)
$$
as $T \to \infty$, uniformly for $g_1, g_2$ in compact subsets of $G$,
then conditions D1 and D2 are satisfied.

\end{prop}

\begin{proof}
Given a bounded $B \subset G$ and $\vre>0$, there is $T_1$ such that
for all $T>T_1$ and all $g_1, g_2 \in B$, 
$$\frac{C(g_1, g_2) \, F(T)}{(1+\vre)^{1/3}} < \lambda  \left(H_T[g_1,
g_2] \right) < (1+\vre)^{1/3} \,  C(g_1, g_2) \, F(T).
$$
Also there is $T_2>0$ such that for all $T>T_2$, 
$$\frac{F\left((1+\delta)T\right)}{F\left(T \right)}
< (1+ \vre)^{1/3}.
$$
Then setting $T_0 = \max\{T_1, T_2\}$ we obtain for $T>T_0$,
\[
\begin{split}
\lambda  \left(H_{(1+\delta)T}[g_1, g_2] \right) & <  (1+\vre)^{1/3} \,
C(g_1, g_2) \, F((1+\delta)T) \\
& < (1+ \vre)^{2/3}   C(g_1, g_2) \,
F(T) \\
& < (1+\vre) \lambda \left(H_{T}[g_1, g_2] \right).
\end{split}
\]
This proves D1. 

Now for any $g_1, g_2 \in G$ we have
\[
\lim_{T \to \infty} \frac{\lambda(H_T[g_1^{-1},g_2])}{\lambda(H_T)} 
= \lim_{T \to \infty} \frac{\lambda(H_T[g_1^{-1},g_2])}{F(T)} \,
\frac{F(T)}{\lambda(H_T)} 
=\frac{C(g_1^{-1}, g_2)}{C(e,e)},
\]
and we obtain D2. The proof of I1 is similar and is omitted.
\end{proof}

\begin{proof}[Proof of Proposition \ref{prop: verifying
axioms 1}.]
Condition UC follows from Propositions \ref{prop: matrix distance functions
have UC}, \ref{prop: right invariant 
metric have UC}. 
If $G, D$ are standard and $G$ is balanced, then condition I1 follows
from Theorems  
\ref{thm: volume asymptotics1} and \ref{thm: symmetric volume
asymptotics 1} and  Proposition \ref{prop: D1 D2 
from asymptotics}, and condition I2 follows from Proposition
\ref{verifying I2} 
and Theorem \ref{thm: using Ratner}.

If $G,H,D$ are standard, then conditions D1, D2 follow from 
Theorems 
\ref{thm: volume asymptotics1} and \ref{thm: symmetric volume
asymptotics 1} and  Proposition \ref{prop: D1 D2 from asymptotics}. 
\end{proof}

\section{Matrix norms and volume computations}
\label{section: verifying the axioms 1}
In this section we consider the asymptotics of volumes of `balls' $H_T$
with respect to a 
matrix norm distance function. We first 
recall some standard details on Cartan (or polar) decomposition and
haar measure.

Let $H$ be a connected semisimple Lie group. 
Let $K$ be a maximal compact subgroup of $H$ and let $A$ be the
associated {\em split
Cartan} subgroup of $H$; by this we mean that $A$ is a maximal
connected group which is invariant under a Cartan involution
associated to $K$, such that $\Ad(A)$ is diagonalizable over
$\R$. Denote
by $\goth{h}$ and $\goth{a}$ the Lie algebras 
of $H$ and $A$ respectively. 
One can write
$$
\goth{h} = \goth{h}_0 \oplus \bigoplus_{\alpha \in \Phi}\goth{h}_{\alpha}
$$
for $\Phi\subset \goth{a}^*$, the dual space of $\goth{a}$, where 
$$
\goth{h}_{\alpha} = \{X \in \goth{h} : \, \forall Y \in \goth{a}, \,
\Ad(\exp(Y))X = e^{\alpha(Y)} X\},
$$
and $\Phi$ is the restricted root system of $H$ relative to $A$. Let $\Phi^+$ be the set of
positive roots 
with respect to some ordering and
$$
\goth{a}^+=\{Y\in\goth{a}: \, \forall \alpha\in\Phi^+,\, \alpha(Y)
\geq 0\}
$$
the corresponding (closed) Weyl chamber.
Let $\rho \in \goth{a}^*$ be defined by
\begin{equation}
\label{eq: defn of rho}
\rho = \frac12 \sum_{\alpha \in \Phi} m_{\alpha} \alpha.
\end{equation}
Let $\Delta = \{\alpha_1,
\ldots, \alpha_r\}\subset \Phi^+$ be the set of simple roots that corresponds to $\mathfrak{a}^+$
and $\{ \til \beta_1,\ldots, \til \beta_r\}$ the dual basis of $\goth{a}$, that is, 
$$\alpha_i(\til \beta_j) = \delta_{ij}.$$
It will be convenient to rescale each $\til \beta_i$ according to
$\rho$. Namely, let 
$$\beta_i = \frac{\til \beta_i}{2 \rho(\til \beta_i)},$$
so that 
\begin{equation}
\label{eq: normalization of beta}
2 \rho \left(\sum_{i=1}^r t_i \beta_i \right) = \sum_{i=1}^r t_i \ \ \
\ \mathrm{and} \ \ i \neq j \ \ \Longrightarrow \ \ \alpha_i(\beta_j)=0.
\end{equation}

Every $h \in H$ can be written
as $h=k_1 a k_2,$ where $k_1,k_2 \in K$ and $a \in A^+ \df
\exp(\goth{a}^+)$. Note that this decomposition is not
unique. Nevertheless, $\lambda$ can be expressed in terms of this
decomposition of $G$ (see
\cite[p. 186]{Helgasson}). Namely, letting $dk$ denote the probability
Haar measure on $K$ we have for $f\in 
C_c(H)$:
\begin{equation}\label{eq_l_H}
\int_H f(h)\, d\lambda(h)=\int_K \int_{\goth{a}^+}\int_K
f(k_1\exp(Y)k_2)\xi(Y)\, dk_1 \, dY \, dk_2,
\end{equation}
where $dY$ is a scalar multiple of the Lebesgue measure on $\goth{a}$,
and 
\begin{equation}
\label{eq: defn xi}
\xi(Y)=\prod_{\alpha \in \Phi^+} \sinh^{m_{\alpha}}(\alpha(Y)), \
\ \ m_{\alpha} = \dim \goth{h}_{\alpha}.
\end{equation}


Now suppose 
$\Psi: H \to \GL(V)$ is an irreducible 
representation. 
Let $d=\dim V$ and
and $V=V_1 \bigoplus \cdots \bigoplus V_s$ be a direct sum decomposition
into the weight-spaces of $\mathfrak{a}$ with weights $\lambda_1, \ldots, \lambda_s \in \goth{a}^*$, that is, for
all $Y \in \goth{a}, \, v \in V_i$,
\begin{equation}\label{eq_last2}
\Psi(\exp(Y)) v = e^{\lambda_i(Y)} v.
\end{equation}
One of the $\lambda_i$'s is a highest weight --
indeed, this is well-known over $\C$ and holds also over $\R$, see
e.g. \cite[Chap. IV]{Guivarch book}. By re-ordering the basis of $V$
assume with no loss of 
generality that the highest weight is $\lambda_1$, i.e.
\begin{equation}
\label{eq: lambda maximal weight}
\forall Y \in \goth{a}^+, \ \forall j \in \{2, \ldots, s\}, \
\ \ \lambda_1(Y) \geq \lambda_j(Y).
\end{equation}

Now let 
\begin{equation}
\label{eq: expression for m1}
m_1 = \min_{i=1, \ldots, r} \lambda_1 (\beta_i).
\end{equation}
By re-ordering $\beta_1, \ldots, \beta_r$ assume with no loss of
generality that 
$$m_1 = \lambda_1(\beta_1).$$

We now make the following assumption.

\begin{definition}
\label{def: condition G}
Say that $\Psi$ {\em satisfies condition {\bf G}} if  
\begin{equation}\label{eq_last3}
m_2 = \min_{j =2, \ldots, r} \lambda_1(\beta_j)
>m_1=\lambda_1(\beta_1).
\end{equation}
\end{definition}

\begin{remark}
\label{remark: condition generic}
{\rm
Irreducible representations of $H$ are usually described in terms of
the corresponding dominant weights.
Note that if condition {\bf G} fails, then for some $i \neq j$, the dominant
weight $\lambda$ of $\Psi$ satisfies 
$$\rho(\til \beta_j) \lambda(\til \beta_i)= \rho(\til \beta_i)
\lambda(\til \beta_j). 
$$
Thus the condition fails only for a subset of $\Psi$ whose
dominant weight is contained in a finite union of proper linear
subspaces of $\goth{a}^*$. 
}
\end{remark}

Let 
$$\mathcal{J} = \{ j \in \{1, \ldots, s\} : \lambda_j (\beta_1) =
m_1\}.$$

Note that $1 \in \mathcal{J}$, and by \equ{eq: lambda maximal weight},
we have for $j \in \mathcal{J}$, 
$$\lambda_j(\beta_1) = \max_{t =1, \ldots,  s} \lambda_t(\beta_1).$$

We fix a basis $V$ such that $\Psi(A)$ is a diagonal subgroup.
For $j=1, \ldots, d$ let $E_j \in \Mat_{d}(\R)$ denote the matrix
whose $j,j$-th entry 
is 1 and all other entries are zero, and for $\tau = (t_2, \ldots,
t_r) \in \R^{r-1}$, define
\begin{equation}\label{eq_last4}
\bar{\tau} = \sum_{i=2}^r t_i \beta_i, \ \ \ \ \ \ \ E_{\tau} = \sum_{(j,k):\,
k\in \mathcal{J},E_j\in V_k} e^{ 
\lambda_k(\bar{\tau})} E_j.
\end{equation}

Also, let 
$$\hat{\Phi} = \{ \alpha \in \Phi: \alpha(\beta_1)=0\} 
,$$
that is, $\hat{\Phi}$ contains those roots whose expression as a
linear combination of simple roots does not involve $\alpha_1$, and let
$$\hat{\xi} (\tau) = \left(\frac12 \right)^{\sum_{\alpha \notin \hat{\Phi}}
m_{\alpha}} \, \prod_{\alpha \in \hat{\Phi}} 
\left(\frac{1}{2} - \frac{1}{2e^{2\alpha(\bar{\tau})}}
\right)^{m_{\alpha}} \, e^{\sum_{i=2}^rt_i}.$$

We say that a collection $\mathcal{N}$ of norms on $\Mat_d(\R)$ is {\bf
bounded} if 
there is $c>1$
such that for any two norms $\| \cdot \|_1, \, \| \cdot \|_2$ in
$\mathcal{N}$, and 
any nonzero $A \in \Mat_d(\R), \ \frac{1}{c}< \frac{\|A\|_1}{\|A\|_2} <c$.

\begin{prop}
\label{prop: integral for D converges}
For any norm $\| \cdot \|$ on $\Mat_{d}(\R)$, the indefinite integral
\begin{equation} 
\label{eq: expression for D}
D = \int_{\tau \in [0,\infty)^{r-1}}
\frac{\hat{\xi}(\tau)}{\|E_{\tau} \|^{1/m_1}} \, d\tau 
\end{equation}
converges, and the convergence is uniform for $\| \cdot \|$ in a
bounded collection of norms.
\end{prop}

\begin{thm}
\label{thm: volume asymptotics}
Let $H$ be connected and semisimple, and let $\Psi: H \to \GL(V)$ be an
irreducible representation satisfying condition {\bf G}. Choosing a
basis of 
$V$, identify $\Psi(H)$ with a subset of $\Mat_d(\R)$, where $d = \dim
V$. Then  for any 
linear norm $\| \cdot \|$ on $\Mat_d(\R)$ we have
$$\lambda(H_T) \sim C T^m,$$
where $m=\frac1{m_1}, \ m_1$ is given by \equ{eq: expression for m1}, 
and
$$C = \int_K \int_K \,  
 \int_{\tau \in [0,\infty)^{r-1}}
\frac{\hat{\xi}(\tau)}{\|\Psi(k_1)E_{\tau} \Psi(k_2)\|^m} \, d\tau \, 
dk_1 \, dk_2.$$ 
\end{thm}

\begin{remark}
\label{remark: about constants}
{\rm
\begin{enumerate}
\item
In the above expression, $m$ depends only on $\Psi$, and
$C$ depends continuously on $\| \cdot \|$; this means that for
any 
$\vre>0$  there is $\delta>0$ such that for a norm $|\cdot |$,
$$\forall v \neq 0, \, \left| 1- \frac{\|v\|}{|v|} \right | < \delta \
\ \ \Longrightarrow \ \ \ \left|C(\| \cdot \|) - C(|\cdot |)\right|
<\vre.$$ 
\item
Since for any $k_1, k_2 \in G$, the map $A \mapsto \|\Psi(g_1) A
\Psi(g_2)\|$ is a 
linear norm on $\Mat_d(\R)$, the integral for $C$ converges in light
of Proposition \ref{prop: integral for D converges}. 
\item
In the general case, that is when $\Psi$ is reducible or does not
satisfy condition 
{\bf G}, there exist $k \in \Z_+, \, m >0$ and $C>0$, where $k$ and
$m$ depend only on 
$\Psi$ and $C$ depends continuously on $\| \cdot \|$, such that 
$$\lambda(H_T) \sim C (\log T)^k T^m.
$$
Details will appear elsewhere.
\end{enumerate}
}
\end{remark}

\begin{prop}
\label{prop: general norm volume}
Suppose $\mathcal{N}$ is a bounded collection of norms on
$\Mat_d(\R)$. Then for any $\delta>0$ there is $T_0$ such that for all
$T \geq T_0$ and all $\| \cdot \| \in \mathcal{N}$ we have:
\begin{equation}
\label{eq: uniform in norms}
\left | \int_{\goth{a}^+(T, \| \cdot \|)} \xi(Y) dY - DT^{m} \right | <
\delta \, T^{m}, 
\end{equation}
where
$$\goth{a}^+(T, \| \cdot \|) = \{Y \in \goth{a}^+: \| \Psi
(\exp(Y)) \| < T\},
$$
$D$ is as in \equ{eq: expression for D}, and $m = 1/m_1$.
\end{prop}

\begin{proof}[Proof of Theorem \ref{thm: volume asymptotics} assuming Proposition \ref{prop: general norm volume}]
 Note that the collection of norms  
$$\{x \mapsto \|\Psi(k_1)x\Psi(k_2)\|: k_1, k_2 \in K\}$$
 is bounded. Thus the result follows using \equ{eq_l_H}.
\end{proof}

\begin{proof}[Proof of Proposition \ref{prop: integral for D converges}]
Let $\mathcal{N}$ be a bounded collection of norms on
$\Mat_d(\R)$. For the purpose of this proof, the notation $X \ll Y$
will mean that $X$ and $Y$ are quantities depending on various
parameters, and there is $C$, depending only on $\mathcal{N}$ and
independent of the other parameters, such that $X \leq CY.$ 

First note that for a fixed basis $\mathcal{B}$ of $\Mat_d(\R)$, for
every $\| \cdot \| \in \mathcal{N}$, and every $A=\sum_{E \in
\mathcal{B}} a_E E  \in
\Mat_d(\R)$, 

$$\| A \| \ll \max_{E \in \mathcal{B}} |a_E|  \ll \|A\|.
$$

Since $E_1, \ldots, E_d$ can be completed to a basis of
$\Mat_{d}(\R)$, this implies 
that $\|E_{\tau}\|$ is bounded below by
a positive constant independent of $\tau\in [0, \infty)^{r-1}$ and of $\| \cdot \| \in
\mathcal{N}$. Thus we need only consider the behavior of the integrand
as $\tau \to \infty$. For $E_j\in V_1$, we have 
\begin{equation}
\label{eq: estimate for denominator}
e^{m_2 \sum_{i=2}^r t_i} \stackrel{\equ{eq_last3}}{\ll}
e^{\sum_{i=2}^r t_i \lambda_1(\beta_i)
}\|E_j\| =
e^{\lambda_1(\bar{\tau})} \| E_j \|
\stackrel{\equ{eq_last4}}{\ll}
\| E_{\tau} \|.
\end{equation}

On the other hand it is clear that 
\begin{equation}
\label{eq: bound via def rho}
\hat{\xi}(\tau) \ll  e^{\sum_{i=2}^r t_i}.
\end{equation}

Putting together \equ{eq: estimate for denominator} and \equ{eq: bound
via def rho}, and using  condition {\bf G} we obtain
\[
\begin{split}
D & =  \int_{\tau \in [0,\infty)^{r-1}}
\frac{\hat{\xi}(\tau)}{\|E_{\tau} \|^{1/m_1}} \, d\tau \\
& \stackrel{\equ{eq: estimate for denominator}}{\ll} \int_{\tau \in [0,\infty)^{r-1}}
\frac{e^{\sum_{i=2}^r t_i}}{\left(e^{m_2 \sum_{i=2}^r
t_i}\right)^{1/m_1}} d\tau \\
&  = 
\int_{\tau \in [0,\infty)^{r-1}} e^{(1-m_2/m_1) \sum_{i=2}^r t_i }d\tau \stackrel{\equ{eq_last3}}{\ll}
\infty. 
\end{split}
\]
\end{proof}

\begin{proof}[Proof of Proposition \ref{prop: 
general norm volume}]
We first show that for all sufficiently large $M$, for all $T>0$, we have 
\begin{equation}
\label{eq: small outside region}
\int_{\goth{b}^+(M,T, \| \cdot \|)} \xi(Y) dY  <
\left(\frac{\delta}{3}\right) \, T^{m}, 
\end{equation}
where
\begin{equation}
\label{eq: defn of goth b}
\goth{b}^+(M,T, \| \cdot \|) = \left \{ \sum_{i=1}^r t_i \beta_i \in
\goth{a}^+(T, \| \cdot \|) : \exists j \in \{2, \ldots ,r\}, \,  t_j \ge M
\right \}. 
\end{equation}

It suffices to show that for any $j \in \{2, \ldots, r\}$ there is $M_0$
such that for $M>M_0$ and any $T>0$  
\begin{equation}
\label{eq: small outside region1}
\int_{\goth{b}} \xi(Y) dY  <
\left(\frac{\delta}{3r}\right) \, T^{m}, 
\end{equation}
where 
\begin{equation*}
\goth{b} = \left \{ \sum_{i=1}^r t_i \beta_i \in
\goth{a}^+(T, \| \cdot \|) : t_j \ge M \right \}.
\end{equation*}

To prove \equ{eq: small outside region1} we will need the following
easy lemma, which is proved 
by induction on $k$: 

\begin{lemma}
\label{lem: integral of exps}
Given $k \in \N$ and $0<m_1 < \cdots < m_k$, there is $c>0$ such that
for all $S>0$,  
$$\int_{\Delta} e^{\sum_{i=1}^kt_i} dt_1 \cdots dt_k \leq c e^{S/m_1},$$
where 
$$\Delta=\left \{(t_1, \ldots, t_k): \forall i,\, t_i \geq 0, \,
\sum_1^k m_i t_i \leq S \right\}.
$$
\end{lemma}

By \equ{eq: lambda maximal weight}, and by comparing with the
max-norm, we find that there is
$\til C$ such that for all $\| \cdot \| \in \mathcal{N}$ 
and all $T>0$, 
\begin{equation}
\begin{split}
\ \goth{a}^+(T, \|
\cdot \|) 
\subset
 \left \{\sum_{i=1}^r s_i \beta_i : \forall i, \, s_i \geq 0,\,
\sum_{i=1}^r s_i 
\lambda_1(\beta_i) \leq \log T + \til C \right \}. 
\end{split}
\end{equation}

Suppose $ j \in \{2, \ldots, r\}$ and $t_j\ge M$. To simplify notation
suppose $j=r$. Then we obtain that 
\begin{align*}
\goth{b} \subset  \left \{ \sum_{i=1}^r s_i
\beta_i: \forall i, \, s_i \geq 0, \, s_r\ge M,\, \sum_{i=1}^{r} s_i \lambda_1(\beta_i) < \log T + \til C \right \}
\end{align*}

Now using condition {\bf G},  choose $\til m_1 < \til m_2 < \cdots <
\til m_r$ such that 
$$\til m_1 = m_1 \ \  \mathrm{and} \ \ \til m_j \leq
\lambda_1(\beta_j), \ \ j\in \{2, \ldots, r\}.$$
Then it follows that $\goth{b} \subset M\beta_r+\goth{c}$ (that is, the translation of $\goth{c}$ by the vector $M\beta_r$) where
$$
\goth{c} =  \left \{
\sum_{i=1}^r s_i 
\beta_i: \forall i, \, s_i \geq 0, \, \sum_{i=1}^{r} s_i \til m_i <
\log T + \til C - 
m_2 M  \right \}.
$$

Applying Lemma \ref{lem: integral of exps} we find that
\[
\begin{split}
 \int_{\goth{b}} \xi(Y) dY  & \leq  \int_{\goth{b}} e^{2 \rho \left(
\sum_{1}^r t_i \beta_i
\right)} dt_1 \cdots dt_r 
\leq
\int_{M\beta_r+\goth{c} } e^{\sum_1^r t_i} dt_1 \cdots dt_r \\
& = e^M
\int_{\goth{c} } e^{\sum_1^r t_i} dt_1 \cdots dt_r 
\leq c e^{M+( \log T + \til C -
m_2 M)/m_1 } \\
& \le ce^{\til C} e^{ (1 -
m_2 /m_1)M} \, T^{1/m_1}.
\end{split}
\]
Since $m_2>m_1$, this implies \equ{eq: small outside
region1}, and we have proved  \equ{eq: small outside
region}.

\medskip

Now set
\begin{equation}
\label{eq: expression for C tau}
C_{\tau} = \frac{\hat{\xi}(\tau)}{\| E_{\tau} \|^m},
\end{equation}
so that 
\begin{equation}
\label{eq: D is integral of C}
D = \int_{[0,\infty)^{r-1}} C_{\tau} d\tau.
\end{equation}

Using Proposition \ref{prop: integral for D converges}, we may assume
by enlarging $M$ that 
\begin{equation}
\label{eq: D almost}
\left|\int_{\R^{r-1} \sm [0,M]^{r-1}}
C_{\tau}  \, d\tau 
\right| < \delta/3.
\end{equation}

\medskip

We will show below that
\begin{claim}
\label{claim: one}
There is $T_0$ such that for all $\| \cdot \| \in \mathcal{N}$, all
$T>T_0$ and all 
$\tau =(t_2, \ldots, t_r) \in [0,M]^{r-1}$, 
\begin{equation}
\label{eq: convergence for fixed tau}
\left | \frac{\int_{\goth{d}(\tau, T)}\xi\left (\sum_{i=1}^r t_i
\beta_i \right )dt_1}{C_{\tau}T^m}-1 
\right | < 
\frac{\delta}{3M^{r-1} \, \max_{\tau \in [0,M]^{r-1}} C_{\tau}}, 
\end{equation}
where 
$$\goth{d}(\tau, T) = \left\{t_1 : (t_1, t_2, \ldots, t_r) \in
\goth{a}^+(T, \| \cdot 
\|)\right\}.
$$
\end{claim}

Assuming the validity of Claim \ref{claim: one}, and writing 
$$\goth{a} = \goth{a}^+(T, \| \cdot \|), \ \ \ \goth{b}= \goth{b}^+(M,
T, \| \cdot 
\|) ,
$$
we have for all $T>T_0$:

\[
\begin{split}
 & \left | \int_{\goth{a}} \xi(Y) dY - DT^{m} \right | \\
 & \stackrel{\equ{eq: D is integral of C}}{<}  \left| \int_{\goth{a}
\sm 
\goth{b} } \xi(Y) dY - T^m \int_{[0,M]^{r-1}} C_{\tau} d\tau
\right | \\
& + \left
|\int_{\goth{b}} \xi(Y)dY \right | + T^m \left| \int_{\R^{r-1} \sm
[0,M]^{r-1}} C_{\tau} d\tau 
\right | \\
& \stackrel{\equ{eq: small outside
region}, \equ{eq: D almost}}{<} 
\left|\int_{[0,M]^{r-1}} \, \left[ \int_{\goth{d}(\tau, T)} \xi(t_1\beta_1
+ \bar{\tau})
\, dt_1  - C_{\tau} T^m \right ]\, d\tau
\right| +  \frac{2\delta}{3} T^m \\ 
&  \stackrel{\equ{eq: convergence for fixed
tau}}{< } \delta T^m. 
\end{split}
\]

In order to prove the claim, we will first show that for any $\eta>0$
there is $T_0$ such that for all $\| \cdot \| \in \mathcal{N}$ and all
$T>T_0$,  
\begin{equation}
\label{eq: estimate on goth d}
\left [0, \frac{\log T - \log \|E_{\tau}\| - \eta}{m_1} \right ] \subset
\goth{d}(\tau, T) \subset \left [0,\frac{\log T - \log \|E_{\tau}
\| + \eta}{m_1}\right ]. 
\end{equation}

Indeed, suppose that $t_1 \in \goth{d}(\tau, T)$, that is,
$$\left \|\Psi\left(\exp(\sum_{i=1}^r t_i \beta_i) \right) \right
\|<T.$$ 
Write $u = \max_{s \notin \mathcal{J}} \lambda_s(\beta_1)<m_1.$ Then we
have: 
\[
\begin{split}
\left \|\Psi\left(\exp(\sum_{i=1}^r t_i \beta_i) \right) \right
\| & = \left \| \sum_{(j,k):\,
E_j\in V_k} e^{\lambda_k\left(\sum_{i=1}^r t_i
\beta_i  \right)} E_j \right \| \\
& = \left \| \sum_{(j,k):\,
E_j\in V_k} e^{t_1 \lambda_k(\beta_1)} \,
e^{\sum_{i=2}^r t_i \lambda_k (\beta_i)}E_j \right  \| \\
& = \left \|e^{m_1 t_1} E_{\tau} + \sum_{(j,k): k \notin \mathcal{J},E_j\in V_k} e^{t_1
\lambda_k(\beta_1)} \, 
e^{\sum_{i=2}^r t_i \lambda_k (\beta_i)}E_j  \right \| \\
& \geq e^{m_1t_1} \|E_{\tau}\|- Ce^{ut_1}
\end{split}
\]
(here $C$ is a constant depending only on $\mathcal{N}$ and $M$). 
Adjusting $C$ if necessary we obtain:
\begin{equation}
\label{eq: adjust C}
e^{m_1t_1} \|E_{\tau} \| (1 - Ce^{(u-m_1)t_1}) < T.
\end{equation}
Taking $t_1$ and $T$ large enough we can ensure that 
$$-\log \left (1 -Ce^{(u-m_1)t_1} \right) < \eta,
$$
and by plugging this in \equ{eq: adjust C} and taking logs we find:
$$m_1 t_1 < \log T - \log \|E_{\tau} \| + \eta.$$
This proves the right hand inclusion in \equ{eq: estimate on goth
d}. The proof of the left hand 
inclusion is similar. 

\medskip

For any $\alpha \in \Phi \sm \hat{\Phi}$, $\alpha(\sum_{i=1}^r t_i
\beta_i) \to_{t_1 \to \infty} \infty$ and hence 
$$\frac{ \sinh \left(\alpha(\sum_{i=1}^r t_i \beta_i)
\right)}{e^{\alpha(\sum_{i=1}^r t_i \beta_i )}} \mathop{\longrightarrow}_{t_1
\to \infty} \frac{1}{2}.
$$
Hence for fixed $\tau$, 
\[
\begin{split}
\frac{\xi\left(\sum_{i=1}^r t_i \beta_i
\right)}{e^{t_1}} & = \frac{\xi(\sum_{i=1}^r t_i \beta_i) \,
}{e^{\sum_{i=1}^r t_i }} \, e^{\sum_{i=2}^r t_i} \\
& = \frac{\prod_{\alpha \in \Phi} \sinh\left(\alpha(\sum_{i=1}^r t_i
\beta_i)\right)^{m_{\alpha}}}{e^{2
\rho(\sum_{i=1}^r t_i \beta_i)}} \, e^{\sum_{i=2}^r t_i}\\ 
& = 
\prod_{\alpha \in \Phi} \left( \frac{ \sinh\left(
\alpha(\sum_{i=1}^r t_i \beta_i) \right)}{e^{\alpha(\sum_{i=1}^r t_i
\beta_i )}}\right)^{m_{\alpha}} \, e^{\sum_{i=2}^r t_i}
\\
& \mathop{\longrightarrow}_{t_1 \to \infty} \left(\frac12 \right)^{\sum_{\alpha
\notin \hat{\Phi}} 
m_{\alpha}} \prod_{\alpha \in \hat{\Phi}} \left( \frac{ \sinh\left( 
\alpha(\sum_{i=1}^r t_i \beta_i) \right)}{e^{\alpha(\sum_{i=1}^r t_i
\beta_i )}}\right)^{m_{\alpha}} \, e^{\sum_{i=2}^r t_i}
\\ 
& =\hat{\xi} \left (
\tau \right ).
\end{split}
\]
Note also that the convergence is uniform for $\tau$ in compact sets. It follows that 
\begin{equation}
\label{eq: uniform convergence}
 \frac{ \hat{\xi}(\tau) \int_{\goth{d}(\tau, T)} e^{t_1}
dt_1}{\int_{\goth{d}(\tau, T)}\xi\left(\sum_{i=1}^r t_i \beta_i \right)
dt_1} \mathop{\longrightarrow}_{T \to \infty} 1,
\end{equation}
and the convergence is uniform for $\tau \in [0,M]^{r-1}$.

\medskip

We now have:
\[
\begin{split}
\lim_{T \to \infty} \frac{\int_{\goth{d}(\tau, T)}\xi\left (\sum_{i=1}^r t_1
\beta_i \right )dt_1}{C_{\tau}T^m} & \stackrel{\equ{eq: uniform
convergence}}{=} \lim_{T \to \infty} \frac{\hat{\xi}(\tau) \int_{\goth{d}(\tau, T)} e^{t_1}
dt_1 }{C_{\tau}T^m} \\
& \stackrel{\equ{eq: estimate on goth d}}{=} \lim_{T \to \infty}
\frac{\hat{\xi}(\tau) \int_0^{(\log T - \log
\|E_{\tau}\|)/m_1}e^{t_1} dt_1}{C_{\tau}T^m}  \\
& = \lim_{T \to \infty} \frac{\hat{\xi}(\tau)}{C_{\tau}T^m}
\left(\frac{T}{\|E_{\tau} 
\|}\right) ^{1/m_1} \stackrel{\equ{eq: expression for C tau}}{=}1.
\end{split}
\]
The convergence in the above expression is uniform in $\tau$,
therefore Claim \ref{claim: one} is valid. This completes the proof of the
Proposition.  
\end{proof}

\section{Balanced semisimple groups}
As a corollary of the computations in the previous sections, we obtain
some information about balanced semisimple groups (cf. Definition
\ref{def: balanced}). 
First we have:

\begin{prop}
\label{prop: contribution to volume}
If $D$ is a matrix norm distance function, corresponding to a 
representation $\Psi: G \to \GL_d(\R) \subset \Mat_{d}(\R)$ with
compact kernel and a norm
$\| \cdot \|$ on 
$\Mat_d(\R)$, then the condition that $H$ is balanced depends on
$\Psi$ but not on $\| \cdot \|$.
\end{prop}

\begin{proof}
Let  $\| \cdot \|, \, \| \cdot \|'$ be two
norms on $\Mat_d(\R)$, and suppose that $H$ is balanced with respect
to $\| \cdot \|$. 
Then for some $C>1$ and every $x\in \Mat_d(\R)$,
\begin{equation}\label{eq_norms_eqv}
\frac{\|x\|}{C} \le\|x\|'\le C\|x\|.
\end{equation}

Let $H_i, \, i=1, \ldots, t$ be the simple factors of $H$ and let
$\sigma_i : H \to H_i$ be measurable sections. 
Given $j \in \{1, \ldots, t\}$, 
$g_1, g_2 \in G$ and a compact $L \subset H_j$ let
\[
\begin{split}
S_T & = S_{T}[g_1, g_2]  = \{h \in H_T[g_1, g_2] :
\sigma_j(h) \in L \}, \\
H'_T & =H'_T[g_1, g_2] = \{h \in H: \|g_1 h g_2\|'<T\}, \\
S'_T & =S'_{T}[g_1, g_2] = \{h \in H'_T[g_1, g_2] :
\sigma_j(h) \in 
L\}. 
\end{split}
\]

It follows from (\ref{eq_norms_eqv}) that
$$
S'_{T} \subset S_{CT} \quad\hbox{and}\quad 
H_{T/C} \subset H'_T.
$$

It is a consequence of Theorem \ref{thm: volume asymptotics1} that
$$\limsup_{T \to \infty} \frac{\lambda(H_{CT})}{\lambda(H_{T/C})} < \infty.
$$
Therefore 
\[
\frac{\lambda(S'_T)}{\lambda(H'_T)} 
 \le \frac{\lambda(S_{CT}
)}{\lambda(H_{T/C})} 
 = \frac{\lambda(S_{CT})}{\lambda(H_{CT})} \,
\frac{\lambda(H_{CT})}{\lambda(H_{T/C})} \mathop{\longrightarrow}_{T \to
\infty} 0. 
\]
 This shows that $H$ is balanced with respect to the norm $\|\cdot\|'$.
\end{proof}

We also have the following result, showing that balanced
representations of semisimple nonsimple groups are rather atypical.

\begin{prop}
\label{prop: G and balanced}
Let $H$ be a semisimple nonsimple Lie group, realized as a matrix
group via an  irreducible representation $\Psi: H \to \GL(V).$
Suppose condition {\bf G} is satisfied. Then $H$ is not balanced.

\end{prop}

\begin{proof}
Let $H=H_1 \cdots H_t, \, t\geq 2$ be a representation of $H$ as an
almost direct
product,  let $K$ be a maximal compact
subgroups of $H$, let $A$ be an associated split Cartan subgroup of
$H$, and for each $i
\in \{1, \ldots, t\}$ write $A_i = H_i \cap A$, so that $A_i$ is a split
Cartan subgroup of $H_i$. We write $\Phi =
\Phi_1 \cup \cdots \cup \Phi_t$
where $\Phi_i$ is the root system corresponding to $(H_i, A_i)$. 
Also fix measurable sections $\sigma_i: H \to H_i$.

Let the notations be as in \S \ref{section: verifying the axioms
1}. There is a partition 
$$\Delta = \Delta_1 \cup \cdots \cup \Delta_t$$
of the simple roots of $\goth{a}$ such that for each $j \in \{1,
\ldots, t\}$, 
$\Delta_j$ is a set of simple roots for $\goth{a}_j$. 
Assume by reordering that $\alpha_1 \in \Delta_1$, and let $j \in
\{2, \ldots, t\}$. Since $\Delta$ is indexed by $\{1, \ldots, r\}$, we
think of $\Delta_j$ as a subset of $\{1, \ldots, r\}$. For each $M>0$
and $T>0$, and each norm $\| \cdot 
\|$ on $\Mat_d(\R)$, define
$$\goth{a}^+(M,T, \| \cdot \|)= \left \{Y = \sum_{i=1}^r s_i 
\beta_i \in \goth{a}^+(T, \| \cdot \|)  : \forall i \in \Delta_j,\,
s_i \leq M \right \}. 
$$
Then $\goth{a}_j \cap \goth{a}^+(M,T, \| \cdot \|) $ is compact for
each $M>0$ and
hence 
$$L=L(M) = \bigcup_{T} \sigma_j \left (K \exp(\goth{a}^+(M,T,
\| \cdot \|)) K \right ) \subset H_j
$$
is precompact. 

Now let
$g_1=g_2=e,$ let  $C$ be as in Theorem \ref{thm: volume asymptotics}, let
$0<\delta < C$, and let 
$\goth{b}(M, T, \| \cdot \|)$ be as in \equ{eq: defn of goth b}. Since
the complement of $\goth{b}(M, T, \| \cdot \|)$ is contained in
$\goth{a}^+(M,T, \| \cdot \|)$, we see that for large $M$ \equ{eq:
small outside 
region} contradicts \equ{eq: balanced}. 
\end{proof}

\section{Riemannian skew balls and volume computations}
\label{section: verifying the axioms 2}

Suppose $G$ is a connected semisimple Lie group and $H$ is its
connected semisimple Lie subgroup. The main result of this section is
a computation of the asymptotics of 
the volume growth for certain `skew balls' in $H$,
with respect to a symmetric space
distance function.

Before formulating the precise result we introduce some notation. 
We are given a semisimple Lie group $G$, a semisimple Lie subgroup $H$
and a maximal compact subgroup
$K$. It follows from a theorem of Mostow \cite{mos} that for a
conjugate $H'$ of $H$, we can choose a maximal
compact subgroup $L$ of 
$H'$ respectively so that $ L \subset K$ and split Cartan
subgroups $D$ and $A$ of
$G$ and $H'$ (associated to $K$ and $L$) respectively so that $A
\subset D$. 
Let $\mathfrak{g}$, $\mathfrak{h}$, $\mathfrak{d}$, $\mathfrak{a}$ denote the
corresponding Lie algebras.

Let $X = K \backslash G$ be the right symmetric space of $G$, and let
$P: G \to X,\, d(\cdot , \cdot)$ be as in \S \ref{subsection:
symmetric spaces}. For the case of the left symmetric space $G/K$, see Remark \ref{r_left} below.
We will write $\bar{g} = P(g)$ to simplify
notation. Our goal will be to determine the asymptotics of
the volumes of all skew balls $H_T[g_1, g_2]$, as in \equ{eq: def skew
balls}. Since applying a conjugation in $G$ permutes the set of skew
balls, with no loss of generality we will replace $H$ with $H'$. 
The Riemannian metric on $X$ induces a scalar product
$(\cdot,\cdot)$ on $\mathfrak{d}$. We denote the corresponding norm by
$\| \cdot \|$,
and write
$$\goth{d}^0 = \{Y \in \goth{d} : \|Y\|=1 \}, \ \ \ \ \ \goth{a}^0 =
\goth{a} \cap \goth{d}^0.
$$
We recall the following facts about the geometry of $X$, see
e.g. \cite{Ballmann, Eberlein} for more details. 

\begin{prop}
\label{prop: geometry of X}
\begin{enumerate}
\item $X$ is a complete Riemannian manifold of nonpositive curvature.
\item The map $\goth{d} \to X, \ Y \mapsto P(\exp(Y)) $
is an isometry. In particular, the submanifold $P(D)$ is a totally
geodesic subset and for any $Y \in \goth{d}^0$,
the path 
$$t \mapsto  \gamma_Y(t)=\bar{e}\, \exp(tY)$$
is a unit speed geodesic. There is a continuous {\em
Busemann 
function} $\beta: X \times \goth{d}^0 \to \R$ such that 
$$\lim_{t \to \infty} d(\gamma_{Y}(t), x) - t = \beta(x,
Y)$$
and the convergence is uniform over compact subsets of $X\times \goth{d}^0$.

\item For any $g\in G$, $x\in X$, and $Y\in\goth{d}^0$,
$$\lim_{t \to \infty} d(\bar{g}\, \exp(tY), x)-d(\bar{g}\, \exp(tY),\bar{e}) = \beta(x,Y)$$
and the convergence is uniform over compact subsets of $G\times X\times \goth{d}^0$.
\end{enumerate}
\end{prop}
\begin{proof}
Parts (1) and (2) are well-known. 
To prove (3), note that $$d\left(\bar{e}\, \exp(tY), \bar{g}\, \exp(tY)\right) =
d(\bar{e},\bar{g}),$$ that is, $\{\bar{g} \, \exp(tY): t \ge 0\}$ stays
within bounded distance of the geodesic $\{\bar{e} \, \exp(tY): t \ge
0\}$. Now the result can be proved using \cite[Proposition 2.5]{Ballmann}.
\ignore{
To prove (3), we denote by $\sigma_t = \sigma^g_t$ 
the unit speed geodesic from $\bar{e}$ to $\bar{g}\, \exp(tY)$. Let
$\ell_t=d(\bar{e},\bar{g}\, \exp(tY))$ be the length of this geodesic. 
By the triangle inequality, for every $t \geq 0$, 
\begin{equation}
\label{eq: triangle}
|\ell_t-t|\le d\left(\bar{e}\, \exp(tY), \bar{g}\, \exp(tY)\right) =
d(\bar{e},\bar{g}).
\end{equation}
It is a feature of nonpositively curved geodesic spaces that the
function $s\mapsto d(\sigma_t(s),\gamma_Y(s))$ is convex, hence for $0\le s\le t$,
\[
\begin{split}
d(\sigma_t(s),\gamma_Y(s)) & \le \frac{t-s}{t}
d(\sigma_t(0),\gamma_Y(0))+\frac{s}{t}
d(\sigma_t(t),\gamma_Y(t)) \\
& =\frac{s}{t} d(\sigma_t(t),\gamma_Y(t)). 
\end{split}
\]
Also, by \equ{eq: triangle},
\[
\begin{split}
d(\sigma_t(t),\gamma_Y(t)) & \le
d(\sigma_t(\ell_t),\gamma_Y(t))+d(\bar{e}, \bar{g})\\
&=d(\bar{g} \, \exp(tY), \bar{e} \, \exp(tY))+d(\bar{e}, \bar{g}) \\
& =2d(\bar{e}, \bar{g}).
\end{split}
\]
Thus, for fixed $s$, $d(\sigma_t(s),\gamma_Y(s))\to 0$ as
$t\to\infty$, i.e. the geodesic $\sigma_t$ converges to the geodesic
$\gamma_Y$, and the convergence is uniform for $g$ and $s$ in a
compact subset of $G \times \R$. 
Hence, by \cite[Proposition 2.5]{Ballmann},
\begin{align*}
&\lim_{t \to \infty} d(\bar{g}\, \exp(tY), x)-d(\bar{g} \, \exp(tY), \bar{e})\\
=&\lim_{t \to \infty} d(\gamma_Y(t), x)-t=\beta(x,Y).
\end{align*}
This implies (3).} 
\end{proof}

Let $\goth{a}^+$ be a positive Weyl chamber in $\goth{a}$ that
corresponds to a system of positive roots $\Phi^+\subset \Phi$. 
For $T>0$, put $\goth{a}^+_T=\{a\in\goth{a}^+:\|a\|\le T\}$.
We will require the following standard fact concerning the functional
$\rho$ defined in \equ{eq: defn of
rho}.

\begin{lemma}
\label{lem: Weyl chamber wall}
There is a unique $Y_{\max} \in \goth{a}^+_1 \cap 
\interior \, \goth{a}^+ $ such that 
$\rho(Y_{\max}) = \max_{Y \in \goth{a}^+_1} \rho(Y).$
\end{lemma}

With this notation we have:

\begin{thm}
\label{thm: symmetric volume asymptotics}
Let $G$ be a semisimple Lie group and $H$ a connected semisimple
subgroup.  Then for any
$g_1, \, g_2 \in G$, 
\begin{equation}
\label{eq: volume growth}
\lambda\left(
H_T[g_1, g_2]
\right )
\sim C(g_1,g_2) T^{(d-1)/2}e^{\delta 
T},
\end{equation}
where  $d=\dim A, \,  \delta = 2\rho(Y_{\max}),$ and
\begin{align*}
 C(g_1, g_2) = &\left(\int_L \exp\left(-\delta \beta\left (\bar{g_1}
\ell,-Y_{\max}\right) \right)d\ell\right)\\ 
  &\times\left(\int_L \exp\left(-\delta \beta\left (\cl{g_2^{-1}}
\ell^{-1},Y_{\max}\right) \right)d\ell\right).
 \end{align*}
The convergence in \equ{eq: volume growth} is uniform for $g_1, g_2$
in compact subsets of $G$.
\end{thm}

We will need another standard result about integration of exponential
functions on balls. 
 Lebesgue measure on $\R^r$
is denoted by $dY$. We call $S \subset \R^r$ a {\em convex cone} 
$S$ is convex and for any $\mathbf{s} \in S$, the
ray $\{t \mathbf{s} 
: t>0 \}$ is contained in $S$. 

\begin{lemma}\label{lem_n}
Let
$\lambda$ be a linear functional on $\R^r$, and let 
$\delta = \max_{ Y \in \goth{a}^+_1}\lambda(Y).
$
Assume that $S \subset \goth{a}^+$ is an open convex cone such that  
$\delta = \max_{Y \in \goth{a}^+_1 \cap S} \lambda(Y). 
$
Then there is a constant $C$ such that
$$
\int_{\goth{a}^+_T \cap S} e^{\lambda(Y)}\, dY \sim
 \int_{\goth{a}^+_T}e^{\lambda(Y)}\, dY \sim C T^{(d-1)/2} e^{\delta T} .
$$
\end{lemma}

\begin{proof}[Proof of Theorem \ref{thm: symmetric volume asymptotics}
assuming Lemmas \ref{lem: Weyl chamber wall} and  \ref{lem_n}] 
\ignore{  Let $\Phi\subset \mathfrak{a}^*$  
be the restricted root system. Fix an ordering on $\Phi$, and let 
$\Phi^+$ be the set of positive roots and $\goth{a}^+$ the
corresponding positive Weyl chamber. }
Define $\xi$ by \equ{eq: defn xi}. Let
\begin{eqnarray*}
\mathfrak{a}^+_T[g_1,g_2]&=&\{Y\in \mathfrak{a}^+: d(\bar{g_1}\,
\exp(Y) g_2,
\bar{e})\le T\}\\ 
\mathfrak{a}^+_T&=&\mathfrak{a}^+_T[e,e] = \{Y \in \goth{a}^+ :
\|Y\| \leq T \}.
\end{eqnarray*}
First,
we determine the 
asymptotics of  
$$
\psi_{g_1,g_2}(T) \df \int_{\mathfrak{a}^+_T[g_1,g_2]}\xi(Y)dY
$$
as $T\to\infty$ with uniform convergence for $g_1$ and $g_2$ in a
fixed compact set $E\subset G$. To do this we will replace $\xi$ by
an exponential function.

By the triangle inequality, there exists $C=C(E)>0$ such that for
$g_1,g_2\in E$,  
\begin{equation}\label{eq_inc}
\mathfrak{a}^+_{T-C}[e,e]\subset\mathfrak{a}^+_T[g_1,g_2] \subset
\mathfrak{a}^+_{T+C}[e,e]. 
\end{equation}
 Expanding \equ{eq: defn xi} we obtain that there are
$\lambda_1, \ldots, \lambda_k \in \goth{a}^*, \, a_1, \ldots, a_k \in
\R$ such that 
\begin{equation*}
\xi(Y) = \frac{e^{2 \rho(Y)}}{2^m} + \sum_{i=1}^k a_i
e^{\lambda_i(Y)}, 
\end{equation*}
where $m=\sum_{\alpha \in \Phi} m_{\alpha}$ and for all $Y \in
\interior\, \goth{a}^+$, 
$$
2 \rho(Y) > \max_{i} \lambda_i(Y).
$$

Therefore
\begin{equation*}\label{eq_other}
\psi_{g_1,g_2}(T)=\frac{1}{2^{m}}\int_{\mathfrak{a}^+_T[g_1,g_2]}
e^{2\rho(Y)} \, dY + \sum_{i=1}^k  a_i \int_{\mathfrak{a}^+_T[g_1,g_2]}
e^{\lambda_i(Y)} \, dY.
\end{equation*}

Let $C'$ be a constant such that for all $T>0$ and all $g_1, g_2 \in
E$, the Lebesgue measure of
$\goth{a}^+_T[g_1, g_2]$ is at most $C' \, T^r$. Let $Y_{\max}$ be as in
Lemma \ref{lem: Weyl chamber wall}. Using the compactness of
$\goth{a}^+_1$, we find that there is $\eta>0$ such that
for each $i \in \{1, \ldots, k\}$, 
$$\max_{Y \in \goth{a}^+_1} \lambda_i(Y) \leq  
\delta - \eta.$$

This implies that 
\begin{equation}
\begin{split}
\frac{\psi_{g_1, g_2}(T) - 2^{-m} \int_{\goth{a}^+_T[g_1,g_2]}
e^{2\rho(Y)} \, 
dY}{e^{\delta T}} & = \frac{\sum_{i=1}^k a_i \int_{\goth{a}^+_T[g_1,
g_2]} e^{\lambda_i(Y)}\, dY}{e^{\delta T}} \\
& \leq \frac{ \sum_{i=1}^k C'T^r \,|a_i| \max_{Y \in \goth{a}^+_{T+C}}
e^{\lambda_i(Y)}}{e^{\delta T}} \\
& \leq \left (e^{\max_i \lambda_i(C)} \, C' T^r \sum_{i=1}^k |a_i|
\right) \, e^{-\eta T} \\
&
\mathop{\longrightarrow}_{T \to 
\infty} 0.  
\end{split}
\end{equation}

 Hence, 
in order to derive the asymptotics of $\psi_{g_1,g_2}(T)$,
it suffices to find $D_{g_1,g_2}>0$ such that 
\begin{equation}\label{eq_dom}
\int_{\mathfrak{a}_T^+[g_1,g_2]} e^{2\rho(Y)} \, dY \sim
D_{g_1,g_2}\, T^{(d-1)/2} e^{\delta T} 
\end{equation}
as $T\to\infty$ uniformly on $g_1,g_2\in E$.

\ignore{Let $S\subset \mathfrak{a}^+$ be a cone that contains $v_\rho$ in its
interior. Define  
$$
\mathfrak{a}_T^S[g_1,g_2]=\mathfrak{a}_T^+[g_1,g_2]\cap S.
$$
Note that 
$$
\max\{2\rho(a): a\in {\mathfrak{a}_1^+}\}>\max\{2\rho (a): a\in{\mathfrak{a}_1^+}-S\}\stackrel{def}{=}\delta''.
$$
By (\ref{eq_inc}),
$$
\int_{{\mathfrak{a}_T^+[g_1,g_2]}-S} e^{2\rho(a)} da\le
\int_{{\mathfrak{a}_{T+C}^+[e,e]}-S} e^{2\rho(a)} da\le
e^{\delta''(T+C)}\cdot\hbox{Vol}({\mathfrak{a}_{T+C}^+})=o(e^{\delta T})
$$
as $T\to\infty$. This implies that 
\begin{equation}\label{eq_check}
\psi_{g_1,g_2}(T)\sim\int_{\mathfrak{a}_T^S[g_1,g_2]} e^{2\rho(a)} da
\stackrel{def}{=}\psi^S_{g_1,g_2}(T) 
\end{equation}
as $T\to\infty$ provided that we show that
\begin{equation}\label{eq_dom0}
\psi^S_{g_1,g_2}(T)\ge C\cdot T^{(d-1)/2} e^{\delta T}
\end{equation}
as $T\to\infty$ for some $C>0$.
}

We will use polar coordinates on $\mathfrak{a}^+$. Thus we will
represent each $Y \in \goth{a}^+$ as $t \omega$, where $t = \|Y\|$ and
$\omega \in \goth{a}^1 = \goth{a}^0 \cap \goth{a}^+.$

Let
\[
\begin{split}
r_T=r_T(\omega,g_1,g_2) & =\inf \left \{r>0: r\omega \notin
\mathfrak{a}_T^+[g_1,g_2] \right \} \\
& = \inf \left \{r>0: d \left (\overline{g_1}\, \exp(r
\omega), \overline{g_2^{-1}}\right ) > T \right\},\\
R_T=R_T(\omega,g_1,g_2) & =\sup \left \{r>0: r\omega \in
\mathfrak{a}_T^+[g_1,g_2] \right \} \\
& = \sup \left \{r>0: d \left (\overline{g_1} \, \exp(r
\omega), \overline{g_2^{-1}}\right ) \le T \right\}.
\end{split}
\]
Then
\begin{equation}
\label{eq_incl}
\{r \omega : \omega\in \goth{a}^1, \, 0\le
r<r_T\}\subset \mathfrak{a}^+_T[g_1,g_2] 
\subset\{r \omega : \omega\in \goth{a}^1, \,
0\le r\le R_T\}. 
\end{equation}

By (\ref{eq_inc}), $r_T(\omega,g_1,g_2)\to\infty$ as $T\to\infty$
uniformly for $g_1,g_2\in E$ and $\omega\in \mathfrak{a}^1$. 

Let $s_T>0$ be such that $d \left (\overline{g_1}\, \exp(s_T\omega),
\overline{g_2^{-1}}\right )=T$ and $s_T\to\infty$. 
It follows from continuity that this condition holds for $s_T=r_T$ as well as
for $s_T=R_T$. 
By Proposition \ref{prop: geometry of X}(3),

\begin{equation}\label{eq_st}
\begin{split}
&\lim_{T\to\infty} (T-s_T) \\
&=\lim_{T\to\infty} \left[  d \left
(\overline{g_1} \exp(s_T\omega), \overline{g_2^{-1}}\right )-d \left
(\overline{g_1} \exp(s_T\omega), \overline{e}\right)
+ d \left (\overline{g_1} \exp(s_T\omega),
\overline{e}\right)-s_T \right] \\
& =\beta\left(\overline{g_2^{-1}},\omega\right)+\lim_{T\to\infty}
\left( d \left (\bar{e} \exp(-s_T\omega), 
\overline{g_1}\right)-s_T\right) \\
& = 
\beta\left(\overline{g_2^{-1}},\omega\right)+\beta\left(\overline{g_1},-\omega\right).
\end{split}
\end{equation}

In particular, this shows that
\begin{equation}\label{eq_rR}
R_T(\omega,g_1,g_2)-r_T(\omega,g_1,g_2)\to 0\quad\hbox{as}\quad T\to\infty
\end{equation}
with uniform convergence for $\omega\in\goth{a}^1$ and $g_1,g_2\in E$.
By (\ref{eq_incl}),
\begin{align*}
\int_{\goth{a}^1} \Lambda(r_T(\omega,g_1,g_2))\,d\omega\le
\int_{\mathfrak{a}_T^+[g_1,g_2]} e^{2\rho(Y)} \, dY\\ 
\le \int_{\goth{a}^1} \Lambda(R_T(\omega,g_1,g_2))\,d\omega,
\end{align*}
where
$$
\Lambda(r)=\int_0^r e^{2s\rho(\omega)} s^{d-1}ds,
$$
and $d\omega$ is a volume form on $\goth{a}^1$ such that $dY =t^{d-1}\, dtd\omega.$
One can check that for every $\vre>0$, there exists $\delta>0$ such
that $\Lambda(r+\delta)\le (1+\vre)\Lambda(r)$ 
for all sufficiently large $r$. Thus, it follows from (\ref{eq_rR}) that
$$
\int_{\goth{a}^1} \Lambda(R_T(\omega,g_1,g_2))\,
d\omega\, \sim \, \int_{\goth{a}^1} \Lambda(r_T(\omega,g_1,g_2))\,d\omega 
$$
as $T\to\infty$, uniformly on $g_1,g_2\in E$,
and to prove (\ref{eq_dom}), it suffices to derive asymptotics of one
of these integrals. 

Let 
$$
c_{g_1,g_2}=\beta\left(\overline{g_2^{-1}},Y_{\max}\right) +
\beta\left(\overline{g_1},-Y_{\max}\right).
$$
Let $\vre>0$. 
Since the Busemann function is continuous, we may choose an open cone
$S \subset \goth{a}^+$ sufficiently close to the ray through
$Y_{\max}$, so that 
for $\omega\in S \cap \goth{a}^1$ and $g_1,g_2\in E$,
$$
\left|\beta\left(\overline{g_2^{-1}},\omega\right)+
\beta\left(\overline{g_1},- \omega \right)-c_{g_1,g_2} \right|<\vre.
$$ 
Then by (\ref{eq_st}) for all sufficiently large $T$, $g_1,g_2\in E$, and $\omega\in
S\cap \mathfrak{a}^1$, 
$$
\left|T-c_{g_1,g_2}- r_T(\omega,g_1,g_2) \right | < 2\vre.
$$
This implies that
$$
S \cap \mathfrak{a}^+_{T-c_{g_1,g_2} - 2\vre}  \subset 
\{r \omega : \omega\in \goth{a}^0\cap S, \, 0\le r<r_T \}
\subset
S \cap \mathfrak{a}^+_{T-c_{g_1, g_2}+2\vre},
$$
hence using Lemma \ref{lem_n}
\begin{align*}
&\limsup_{T\to\infty} \frac{\int_{\goth{a}^1}
\Lambda(r_T(\omega,g_1,g_2))\,d\omega}{ T^{(d-1)/2} \, e^{\delta T}} 
 = \limsup_{T\to\infty} \frac{\int_{\mathfrak{a}^+_{T-c_{g_1,g_2}+2\vre}} e^{2
\rho(Y)} \, dY }{ T^{(d-1)/2} \, e^{\delta T}} \\ 
&\le  \limsup_{T\to\infty}
\frac{(T-c_{g_1,g_2}+2\vre)^{(d-1)/2}e^{\delta
(T-c_{g_1,g_2}+2\vre)}}{T^{(d-1)/2}e^{\delta T}} 
= e^{\delta (-c_{g_1,g_2}+2\vre)}.
\end{align*}
Since $c_{g_1,g_2}$ is bounded, convergence is uniform for $g_1,g_2\in E$.

Similarly, one has
$$
\liminf_{T\to\infty} \frac{\int_{\goth{a}^1}
\Lambda(r_T(\omega,g_1,g_2))\,d\omega}{T^{(d-1)/2}e^{\delta T}}\ge
e^{\delta (-c_{g_1,g_2} -2\vre)}, 
$$
and since $\vre$ was arbitrary, we have
$$
\int_{\goth{a}^1} \Lambda(r_T(\omega,g_1,g_2))\,d\omega\sim e^{-\delta
c_{g_1, g_2}}\cdot T^{(d-1)/2}e^{\delta T} 
$$
as $T\to\infty$. This proves \equ{eq_dom} with 
$$D_{g_1,
g_2} = e^{-\delta c_{g_1, g_2}} = \exp \left(-\delta
\left(\beta\left(\overline{g_2^{-1}},Y_{\max}\right) + 
\beta\left(\overline{g_1},-Y_{\max}\right) \right) \right).$$ 
Thus we have proved that 
$$\psi_{g_1, g_2} (T) \sim D_{g_1,
g_2} T^{(d-1)/2} e^{\delta T},
$$
as $T \to \infty$, with uniform convergence for $g_1,g_2\in E$.

We now obtain \equ{eq: volume growth} via \equ{eq_l_H}.
\end{proof}

\begin{remark}\label{r_left} {\rm
We use notations from Section \ref{subsection: symmetric spaces}.
To treat the case of the left symmetric space $G/K$, we observe that 
the transformation $g\mapsto g^{-1}$ maps Riemannian balls $\{g: d(P(e),P(g))<T\}$ to the balls $\{g: d'(P'(e),P'(g))<T\}$.
Hence, denoting $$H'_T[g_1,g_2]=\{h\in H: d'(P'(e),P'(g_1hg_2))<T\},$$ we have
$$
H'_T[g_1,g_2]=H_T[g_2^{-1},g_1^{-1}]^{-1}.
$$
This implies Theorem \ref{thm: symmetric volume asymptotics} for the balls $H'_T[g_1,g_2]$.
}\end{remark}

\begin{proof}[Proof of Lemma \ref{lem: Weyl chamber wall}]
For any $\chi \in \goth{a}^*$, let $v_{\chi} \in \goth{a}$ so that for
all $Y \in \goth{a}$,  
$$(v_{\chi}, Y) = \chi(Y).
$$
Since $\goth{a}_1$ is strictly convex, $\max_{Y \in \goth{a}_1}
\rho(Y)$ is attained at a unique point which we denote by
$Y_{\max}$. It is a standard application of Lagrange multipliers that
$Y_{\max} = \frac{v_{\rho}}{\|v_{\rho}\|},$
so it remains to show that $v_{\rho} \in \interior \, \goth{a}^+,$
that is, for all $\alpha \in \Phi$,
\begin{equation}
\label{eq: positivity}
(v_{\alpha}, v_{\rho})>0.
\end{equation}

The inner product $(\cdot, \cdot)$ is invariant under the action of
the Weyl group, that is, for any $\alpha, \beta, \gamma \in \Phi$, 
$$\left( v_{r_{\alpha}(\beta)}, v_{r_{\alpha}(\gamma)} \right) = \left(
v_{\beta}, v_{\gamma} \right),
$$
where $r_{\alpha} : \goth{a}^* \to \goth{a}^*$ denotes the reflection
in the root $\alpha$.
Let 
$$\Psi  = \Phi^+ \sm \spa (\alpha), \ \ \ \rho' = \frac12 \sum_{\beta \in
\Psi} m_{\beta} \beta,$$
so that 
$$\rho = \rho' + \, \frac12  \sum_{\beta
\in \Phi^+ \sm \Psi} m_\beta \beta .
$$

It is a standard fact about root systems (see
e.g. \cite[Prop. 1.1.2.5]{GWarner})
that $r_{\alpha}(\rho') = \rho'$. This implies that
$$(v_{\alpha}, v_{\rho'}) = (v_{r_{\alpha}(\alpha)},
v_{r_{\alpha}(\rho')}) = (v_{-\alpha}, v_{\rho'}) = -(v_{\alpha},
v_{\rho'}),
$$
so that $(v_{\alpha}, v_{\rho'})=0$. 
Moreover any $\beta \in \Phi^+ \sm \Psi$ is a positive multiple of
$\alpha$ and hence satisfies $(v_{\alpha}, v_{\beta})>0$. 
Therefore 
$$\left(v_{\alpha}, v_{\rho} \right) = \left(v_{\alpha}, v_{\rho'}
\right) + \frac12 \sum_{\beta \in \Phi^+ \sm \Psi} m_{\beta}
\left(v_{\alpha}, 
v_{\beta} \right) >0,
$$
and we have \equ{eq: positivity}.
\end{proof}

\begin{proof}[Proof of Lemma \ref{lem_n}]
Let 
$$g_T(z)= \left \{ \begin{array}{cc} e^{-z} z^{(r-1)/2} \left(2-\frac{z}{T}
\right)^{(r-1)/2}  & \ \ \ \ \ \ \ \  \mathrm{when \ } 0 \leq z \leq 2T \\ 0 
& \ \ \mathrm{when \ } z > 2T \end{array} \right. \, ,$$
and let 
$$g(z) = 2^{(r-1)/2} \, e^{-z} z^{(r-1)/2}.$$
By Lebesgue's dominated convergence theorem, 

\begin{equation}
\label{eq: gamma function}
\int_0^{\infty} g_T(z) \, dz \longrightarrow_{T \to \infty}
\int_0^{\infty} g(z) \, dz = C. 
\end{equation}

By rescaling, we may assume that 
 $\delta =1$. 
Let $\hbox{Vol}_{r-1}$ denote Lebesgue measure on the subspace $\ker
\lambda$. By translation it induces a measure on each affine subspace
parallel to $\ker \lambda$. For each $x \in [-T, T]$, 
$$B_{T,x} = \left\{Y\in  B(0,T): \lambda(Y)=x \right\}
$$
is a ball of radius $\sqrt{T^2 - x^2}$ in a translate of $\ker
\lambda$ and hence  
$$\hbox{Vol}_{r-1} (B_{T,x}) = C' \, \left(T^2 - x^2 \right)^{(r-1)/2},
$$
where $C'$ is a constant depending only on $r$.

Decomposing the integration into slices parallel to $\ker \lambda$, we
have 
\[
\begin{split}
\frac{\int_{B(0,T)} e^{\lambda(Y)} \, dY}{e^{T} \, T^{(r-1)/2}} & =
 \int_{-T}^T \frac{e^x \, \hbox{Vol}_{r-1}(B_{T,x})}{e^{T} \,
T^{(r-1)/2}} \, dx \\
& = C' \, \int_{-T}^T \frac{e^x \, \left(T^2 - x^2
\right)^{(r-1)/2}}{e^T \, T^{(r-1)/2}} \, dx \\
& = C' \, \int_{0}^{2T} e^{-z}
\left(\frac{z (2T-z)}{T}\right) ^{(r-1)/2} \, dz \\
& = C' \, \int_{0}^{\infty} g_T(z) \, dz.
\end{split}
\]

Applying \equ{eq: gamma function} we obtain 
\begin{equation}
\label{eq: first step}
\int_{B(0,T)}
e^{\lambda(Y)} \, dY \sim  C \, T^{(r-1)/2} \, e^T.
\end{equation}

Since $B(0,1)$ is strictly convex, the maximum $\max_{Y \in B(0,1)}
\lambda(Y)$ is attained at a unique point $Y_0$. By the hypothesis
$Y_0 \in S \subset \interior \, \goth{a}^+.$ 
There is $\vre>0$ such that for all $Y \in B(0,1) \sm S,$
$\lambda(Y) < \delta - \vre.$ 
Furthermore there is a cube $\mathcal{C} \subset \R^r$ of 
side length 2, such that 
$$B(0,1) \sm S \subset \mathcal{C}
$$
and 
$$\max_{Y \in \mathcal{C}} \lambda(Y) \leq \delta - \vre.
$$

For $T>0$ let 
$$T \, \mathcal{C} = \left \{t c : c \in \mathcal{C}, t \in [0,T]
\right \}.$$ 
We have 
\[
\begin{split}
\int_{B(0,T) \sm S} e^{\lambda(Y)} \, dY & \leq \int_{T \,
\mathcal{C}} e^{\lambda(Y)} \, dY \\
& \leq ( 2T)^r \max_{Y \in T \,
\mathcal{C}} e^{\lambda(Y)} \\
& \leq C'' T^r e^{(\delta - \vre)T}.
\end{split}
\]
Comparing with \equ{eq: first step} we obtain 
$$\frac{
\int_{B(0,T) \sm S} e^{\lambda(Y)} \, dY
}{
\int_{B(0,T)} e^{\lambda(Y)} \, dY
} \longrightarrow_{T \to \infty} 
0, 
$$
and so $\int_{B(0,T) \cap S} e^{\lambda(Y)} \, dY \sim \int_{B(0,T) }
e^{\lambda(Y)} \, dY
$. The assertion
follows.  
\end{proof}

\section{Ratner theory and linearization}

Our goal in this section is to prove
Theorem \ref{thm: using Ratner}. We will derive it from the following
result:

\begin{thm}
\label{thm: using Ratner, abstract}
Suppose the general setup holds, where $H$ is connected, semisimple,
and 
balanced, and assume also that UC and D1 hold. Suppose $g_0 \in G$
satisfies $\cl{H\pi(g_0)} = G/\Gamma$. Then for every $F \in
C_c(G/\Gamma)$ and any $g_1, g_2 \in G$ we have
\begin{equation}
\label{eq: using Ratner}
\frac{1}{\lambda \left(H_T[g_1, g_2] \right)} \int_{H_T[g_1, g_2]} F(h
\pi(g_0)) \, d\lambda(h) \mathop{\longrightarrow}_{T \to \infty} \int_{G/\Gamma}
F \, dm'.
\end{equation}
\end{thm}

The proof of Theorem \ref{thm: using Ratner, abstract} relies on the
fundamental 
results of Ratner on the dynamics of unipotent flows on homogeneous
spaces, and subsequent work of Dani, Margulis, Mozes, Shah, and others. Specifically, we use a
result of Nimish Shah. 

\ignore{
\begin{thm}\label{th_semi_equi}
Let $G$ be a Lie group, realized as a group of matrices via a faithful
representation $G \to \GL(V)$. Let $H$ be a closed connected
semisimple subgroup of $G$, and suppose that $H$ is balanced. Let
$\Gamma$ be a lattice in $G$.
Then for every $x\in G/\Gamma$ such that $Hx$ is
dense in $G/\Gamma$ and for every $f\in C_c(G/\Gamma)$,
$$
\frac{1}{\lambda(H_T)}\int_{H_T} f(hx) \, d\lambda(h)\mathop{\longrightarrow}_{T
\to \infty} \int_{G/\Gamma} f
\, dm'.
$$
\end{thm}

{\comment need to add a proof for the case that the averaging on $H$
is w.r.t. a symmetric space distance function.}
}
We say that a sequence
$\{h_n\}\subset H$ is {\em strongly 
divergent} if its projection 
on every simple factor is divergent. In the case of a matrix group,
and using the notation of \S 2.7, $\{h_n\}$ is strongly divergent if and only if
$\{\sigma_i(h_n)\} \subset H_i$ has no convergent subsequence for every
simple factor $H_i$ of $H$.

\begin{thm}[Shah]\label{th_shah}
Let $H$ be a connected semisimple Lie subgroup of a Lie group $G$.
Let $A$ be a split Cartan subgroup of $H$, $A^+$ a closed Weyl chamber 
and $K$ a 
maximal compact subgroup of $H$. Let $\Gamma$ be a lattice in $G$.
Let $\mu$ be a finite Borel measure on $K$ that is absolutely continuous
with respect to haar measure on $K$. Suppose that for $x\in G/\Gamma$, $Hx$ is
dense in $G/\Gamma$. Then for every $f\in C_c(G/\Gamma)$ and every
strongly divergent 
sequence $\{a_n\}\subset A^+$,
$$
\int_K f(a_nkx)\, d\mu(k)\to \mu(K) \, \int_{G/\Gamma}
f \, dm'\;\;\textrm{as}\;\; n\to\infty. 
$$
\end{thm}

\begin{remark}
{\rm
Theorem \ref{th_shah} is proved, but not formulated explicitly, in
\cite{Nimish: proc ind} (see the proof of Corollary 1.2).
}
\end{remark}

\begin{proof}[Proof of Theorem \ref{thm: using Ratner, abstract}]
Let $F \in C_c(G/\Gamma)$, and assume without loss of generality that
$F\ge 0$. Let $g_0, g_1, g_2$ be as in the statement of the
theorem. In proving \equ{eq: using Ratner}, to conserve ink
and computer memory, we 
will omit $g_1$ and $g_2$ 
from the notation. Thus $H_T$ stands for $H_T[g_1, g_2]$, $D(x)$ stands
for $D(g_1xg_2)$ 
and so on. 

Fix $\vre>0$, let $\vre_0>0$ such that 
$(1+\vre_0)^2 < 1+\vre,
$
 and by D1, let $\delta>0$ such that for all large enough $T$,
$$\lambda \left(H_{(1+\delta)^2 T} \right) \leq (1+\vre_0) \lambda
\left( H_T \right).$$
By UC, there exists a symmetric neighborhood
$\mathcal{O}$ of 
identity in $K$ such that for all $g \in G$ and $k \in \mathcal{O}$, 
\begin{equation}\label{eq_norm_ineq1}
D(gk) \le (1+\delta)D(g).
\end{equation}
Since $K$ is compact, there exists a finite cover 
$$
K=\bigcup_{i=1}^N \ell_i\mathcal{O}, \ \ \ \ \ell_1, \ldots, \ell_N
\in K.
$$
Denote by $\kappa$ the probability haar measure on $K$. Using a
partition of unity, there are measures $\mu_i, \, i=1, \ldots, N$ on
$K$, absolutely continuous with respect to $\kappa$, such that 
\begin{equation}\label{eq_mu_i_N}
\supp \, \mu_i \subset \ell_i \mathcal{O}, \ \ i=1, \ldots, N \ \ \ \ \ \ \
\mathrm{and} \ \
\sum_{i=1}^N \mu_i =\kappa. 
\end{equation}

Let $\nu$ denote the
measure on $\goth{a}^+$
which is equal to $\xi (Y)\, dY$, where $dY$ is Lebesgue measure. We
decompose $\lambda$ as in \equ{eq_l_H}, so that 
$$d\lambda(h) = d\kappa(k_1) \, d\nu(Y)\, d\kappa(k_2)$$ for $k_i \in K,
\, Y \in \goth{a}^+, \, h=k_1\exp(Y) k_2.$

Let $H=H_1 \cdots H_t$ be the decomposition of $H$ into an almost
direct product of simple factors and let $\sigma_j: H \to H_j$ be
measurable sections. 
Taking $x=\pi(g_0)$ in Theorem \ref{th_shah}, we find that there are
compact subsets $C_j \subset H_j, 
\, j=1, \ldots, t$
such that for any $h \in \til H$
\begin{equation}\label{eq_shah_eq}
\left| \int_K F(hk\pi(g_0))d\mu_i(k) - \mu_i(K)\int_{G/\Gamma} F \, dm'
\right| < \vre_0,
\end{equation}
where 
$$\til H =
\{h \in H: \sigma_j(h) \notin C_j, \, j=1, \ldots, t\}.
$$

For each $k_1, k_2 \in K$ and $T>0$ write
$$\goth{b}(k_1, k_2, T) = \{Y \in \goth{a}^+: \forall j, \,
\sigma_j(\exp(Y)) \notin C_j, \, D(k_1 \exp(Y)k_2)<T \}.
$$

By enlarging each $C_j$ we may assume that \equ{eq_shah_eq} holds
for any 
$h=k_1 \exp(Y)$, with $k_1 \in K, \, Y \in \goth{a}$ such that
$\sigma_j(\exp(Y)) \notin C_j$ for $j=1, \ldots, t.$ Further, we may
assume that  
$$\til H_T = \bigcup_{k_1, k_2 \in K} k_1 \exp \left(\goth{b}(k_1,
k_2, T) \right) k_2.
$$

It follows from \equ{eq_norm_ineq1} and \equ{eq_mu_i_N} that if $k_2
\in \supp \, \mu_i$ 
then for all $T>0$,
\begin{equation}
\label{eq: new one}
\goth{b}(k_1, k_2, T) \subset \goth{b}(k_1,\ell_i, (1+\delta)T)
\subset \goth{b}(k_1, k_2, (1+\delta)^2 T).  
\end{equation}

Writing $x=\pi(g_0)$, we have
\[
\begin{split}
& \int_{\til H_T} F(hx) \, d \lambda (h) \\
 =&  \int_{K} d\kappa(k_1) \,
\int_K d\kappa(k_2) \, \int_{\goth{b}(k_1, k_2, T)} F(k_1 \exp(Y) k_2
x) \, d\nu(Y) \\
\stackrel{\equ{eq_mu_i_N}}{=} & \sum_{i=1}^N \int_K d \kappa(k_1) \,
\int_{\ell_i \mathcal{O}} 
d\mu_i(k_2) \, \int_{\goth{b}(k_1, k_2, T)} F(k_1 \exp(Y) k_2
x) \, d\nu(Y) \\
\stackrel{\equ{eq: new one}}{\leq}  & \sum_{i=1}^N \int_K d \kappa(k_1) \,
\int_{\goth{b}(k_1, \ell_i, (1+\delta)T)}  d\nu(Y) \,
\int_{K}  F(k_1 \exp(Y) k_2 x)\,  d\mu_i(k_2)  \\   
\stackrel{\equ{eq_shah_eq}}{\leq} & (1+\vre_0) \int_{G/\Gamma} F \,
dm' \, \sum_{i=1}^N \int_K d\kappa(k_1) \int_{\goth{b}(k_1, \ell_i,
(1+\delta)T)} \mu_i(K) \, d\nu(Y) \\
\stackrel{\equ{eq: new one}}{\leq} & (1+\vre_0) \int_{G/\Gamma} F \,
dm' \sum_{i=1}^N \int_K d\kappa(k_1) \, \int_K d\mu_i(k_2) \,
\int_{\goth{b}(k_1, k_2, (1+\delta)^2T)}  \, d\nu(Y) \\
\leq & (1+\vre_0) \, \lambda \left(\til H_{(1+\delta)^2 T} \right) \,
\int_{G/\Gamma} F \, 
dm'  \\
\leq & (1+\vre_0)^2 \, \lambda \left(H_T \right) \, \int_{G/\Gamma} F \, dm' .
\end{split}
\]

Thus 
$$\frac{1}{\lambda \left( H_T \right)} \int_{\til H_T} F(hx) \,
d\lambda(h) < (1+\vre) \int_{G/\Gamma} F \, dm'.
$$
Since $H$ is balanced and $\vre$ was arbitrary we have
$$\limsup_{T \to \infty} \frac{1}{\lambda \left( H_T \right)}
\int_{H_T} F(hx) \, 
d\lambda(h) \leq \int_{G/\Gamma} Fdm'.
$$
The proof of the opposite inequality for $\liminf$ is similar, and we
have proved \equ{eq: using Ratner}.
\end{proof}

\begin{proof}[Proof of Theorem \ref{thm: using Ratner}]
Clearly $\left( ** \right)$ follows from \equ{eq: using Ratner} if we
replace $h$ by $h^{-1}$. To justify this, let 
$$\til D (g) = D(g^{-1})$$
and apply Theorem \ref{thm: using Ratner, abstract} to $G, H, \til D$
instead of $G, H, D$. Note that D1 for $\til D$ holds since $H$ is
unimodular, and UC for $\til D$ can be proved as in Propositions
\ref{prop: matrix distance functions have UC} and \ref{prop: right
invariant metric have UC}. 
\end{proof}

\section{Applications}
In this section we prove Corollaries \ref{th_oppenheim}, \ref{cor: dense projections},
\ref{cor_equi_appl2} and \ref{cor: application torus}.
We use the same notation as in \S 1.3. 

\begin{proof}[Proof of Corollary \ref{th_oppenheim}]
In terms of the identification \equ{eq: identify} of $\mathcal{F}_d$
with $G$, the 
map  
$$
\mathcal{F}_d\to \Mat_d(\R), \ \ f\mapsto \bar Q(f)
$$
is given by
\begin{equation}\label{eq_map_id}
G\to \Mat_d(\R) , \ \ g\mapsto {}^t g A_Q g,
\end{equation}
where $A_Q$ is the matrix of the quadratic form $Q$ with respect to
the standard basis $\mathbf{e}$. 
The set  
$$
\bar{Q}(\mathcal{F}_d)=\left\{{}^t g A_Q g:g\in G\right\}
$$
consists of symmetric matrices that have the same determinant and the
same signature as $A_Q$, 
\ignore{In particular, it follows that $\bar Q(\mathcal{F}_d)$ is a subset of
an algebraic set, and it contains an open set of this algebraic set
(such an open set can be defined by conditions on principal minors).
Then $\bar Q(\mathcal{F}_d)$ contains a nonsingular point of this algebraic set. Since $\bar Q(\mathcal{F}_d)$
is a single $G$-orbit, it is a manifold.
The map (\ref{eq_map_id}) induces a homeomorphism }
and we have an identification
$$
\rho: \HG \to \bar Q(\mathcal{F}_d), \  \ \rho(\tau(g)) = {}^t g A_Q
g=\bar{Q}(g\, {\bf e}),  
$$
where $H$ is the orthogonal group of $Q$. It can be shown that
$\bar{Q}(\mathcal{F}_d)$ is an algebraic variety and that 
$\rho$ is an isomorphism of algebraic varieties.

Given a norm $\| \cdot \|$ on $\R^d$, we define a norm on
$\Mat_d(\R)$, which we also denote by $\| \cdot \|$, by
$$
\|g\|=\max_{j=1,\ldots d} \|g_{*j}\|,\quad g=(g_{ij})\in\Mat_d(\R).
$$
Then for $g\in G$, $\|g\,{\bf e}\|<T$ if and only if $\|g\|<T$. 
Let $\Gamma = G(\Z)$  and $A\subset \bar Q(\mathcal{F}_d)$ as in the statement of the corollary.
We have
\begin{align*}
&\#\left\{f\in\mathcal{F}_d(\Z):\|f\|<T,\;\bar Q(f)\in A\right\}\\
= & \#\left\{\gamma\in \Gamma:\|\gamma\|<T,\; \bar Q(\gamma\, {\bf
e})\in A\right\}\\ 
= & \#\left\{\gamma\in \Gamma:\|\gamma\|<T,\; \tau(\gamma) 
\in\rho^{-1}(A)\right\}. 
\end{align*}

Let $H^o$ be
the connected component of the identity in $H$.
The group $H^o$ is a connected semisimple Lie group, and 
it is simple unless $Q$ is of signature $(2,2)$. If $H$ is simple,
then it is balanced. 
If $H$ is of signature $(2,2)$, then it can be shown by a direct
computation that it is balanced.
\ignore{is conjugate to
$\hbox{SO}(2,2)^o$, and it follows 
from Example \ref{so22} that it is balanced. }

Let $P:H^o\backslash G\to H\backslash G$ be the projection map.
Since $P$ has finite fibers, the set $\til A =
P^{-1}(\rho^{-1}(A))$ is a relatively compact subset 
with boundary of measure zero with respect to the
smooth measure class on $H^o\backslash G$. 
It is a well-known consequence of Ratner's orbit-closure theorem (see
e.g. \cite{bp92}) that $H^o\pi(e)$ is 
dense in $G/\Gamma$. Thus we may apply Corollary \ref{cor: putting together,
duality}, and obtain that
\begin{align*}
&\#\left\{ \gamma\in\Gamma:\|\gamma\|<T,\; \gamma H^o \in
\til A \right\}\\ 
 \sim \, & \, m \left ( \left\{g\in G: \|g\|<T,\; g H^o\in \til A
\right \} \right)\\
= \, & m \, \left(\left\{g\in G: \|g\|<T,\; \bar Q(g\, {\bf e})\in A
\right \} \right)\\
= \, & m \, \left( \left\{f\in\mathcal{F}_d:\|f\|<T,\;\bar
Q(f)\in A  \right\} \right).
\end{align*}
This proves the first assertion, and shows that these quantities are
asymptotic to $\nu (A) \lambda(H_T)$, where $\nu$ is the measure
$\nu_{x_0}$ for $x_0=\tau(e)$ defined by \equ{eq_nu_g}.

For the second assertion, we need to compute the asymptotics of
$\lambda(H_T)$. We use the notations of \S 7.  
In a suitable basis the quadratic form $Q$ is given by
$$
Q(x_1,\ldots, x_{p+q})=x_1x_{p+q}+\cdots + x_p x_{q+1}+x_{p+1}^2+\ldots+x_q^2.
$$
Then the split Cartan subgroup of $\hbox{SO}(Q)$ can be chosen to be
$$
A=\hbox{diag}(e^{s_1},\ldots,e^{s_p},1,\ldots,1,e^{-s_p},\ldots, e^{-s_1}),
$$
and the weights of the representation map $Y = (s_1, \ldots, s_p, 0,
-s_p, \ldots, -s_1)$ to its coordinates $
s_1,\ldots,s_p,0,-s_1,\ldots, -s_p.
$
Choosing an order on roots so that the positive root spaces are upper
triangular, we find that the positive roots are
\begin{align*}
s_i-s_j,\; 1\le i< j\le p,\quad &\hbox{multiplicity}=1,\\
s_i,\; 1\le i\le p,\quad &\hbox{multiplicity}=q-p,\\
s_i+s_j,\; 1\le i<j\le p,\quad &\hbox{multiplicity}=1.
\end{align*}
The dominant weight is $\lambda_1(s)=s_1$ and
$$
2\rho(s)=\sum_{i=1}^p (p+q-2i)s_i.
$$

First, we consider the case when $p<q$.
Then the system of simple roots is
$$
\{s_1-s_2,\ldots,s_{p-1}-s_p, s_p\}.
$$
The dual basis of $\mathfrak{a}$ is $\tilde{\beta}_k, \, k=1,\ldots,p,
$ where 
$$\tilde{\beta}_k=(1,\ldots,1,0,\ldots,0) \ \ \ (k \mathrm{\ ones} ).
$$
and
$$
{\beta}_k=\frac{\tilde{\beta}_k}{2\rho(\tilde{\beta}_k)}=\left(\frac{1}{k(p+q-k-1)},\ldots,\frac{1}{k(p+q-k-1)},0,\ldots,0\right).
$$
We have $\lambda_1(\beta_k)=\frac{1}{k(p+q-k-1)}$. Note that $\lambda_1(\beta_k)$ is strictly decreasing for $0\le k\le p\le\frac{p+q-1}{2}$.
Thus, $m_1=\frac{1}{p(q-1)}$ and condition {\bf G} is satisfied, so by
Theorem \ref{thm: volume asymptotics}
$$
\lambda(H_T) \sim CT^{p(q-1)}.
$$

Now if $p=q$, it may be shown that condition {\bf G} is not
satisfied. Recall that in Remark \ref{remark: about constants}(3) we
mentioned a generalization 
of Theorem \ref{thm: volume asymptotics} in the case when condition {\bf
G} does not hold. Using this generalization 
we are able to show that in this case the asymptotics are 
$C(\log T)T^{p(p-1)}$.
We omit the details.
\ignore{
. The system of simple roots
$$
\{s_1-s_2,\ldots,s_{p-1}-s_p, s_{p-1}+s_p\}.
$$
The dual basis of $\mathfrak{a}$ is
\begin{align*}
\tilde{\beta}_k=&(1,\ldots,1,0,\ldots,0), \ \ \ (k \mathrm{\ ones} )
\quad k=1,\ldots,p-2,\\ 
\tilde{\beta}_{p-1}=&\left(\frac{1}{2},\ldots,\frac{1}{2},-\frac{1}{2}\right),\\
\tilde{\beta}_{p}=&\left(\frac{1}{2},\ldots,\frac{1}{2}\right).
\end{align*}
and
\begin{align*}
{\beta}_k=&\left(\frac{1}{k(2p-k-1)},\ldots,\frac{1}{k(2p-k-1)},0,\ldots,0\right),\; k=1,\ldots,p-2,\\
{\beta}_{p-1}=&\left(\frac{1}{p(p-1)},\ldots,\frac{1}{p(p-1)},-\frac{1}{p(p-1)}\right),\\
{\beta}_{p}=&\left(\frac{1}{p(p-1)},\ldots,\frac{1}{p(p-1)}\right).
\end{align*}
Hence $\lambda_1(\beta_k)$ achieves its minimum for $k=p-1,p$ and $m_1=\frac{1}{p(p-1)}$. Thus, in this case the
asymptotics is
$$
C(\log T)T^{p(p-1)}.
$$}
\end{proof}

\medskip

\begin{proof}[Proof of Corollary \ref{cor: dense projections}]
For both assertions we apply Corollary \ref{cor: putting together,
duality} to the matrix 
norm distance function corresponding to $\Psi$ and $\|\cdot \|$, with
$\HG = S$ and $x = P(s_0)$. Note 
that all conditions are satisfied, and $\nu_x$ in our case is equal to
$\alpha(s_0, s)\, d\nu(s).$ So it remains to calculate the
density $\alpha$.

We begin with the first case. Since
$V = V_H \oplus V_S$, for any $s \in S, \, h \in H$, $\Psi(sh) =
\Psi(s) + \Psi(h) - \mathrm{Id}$, where $\mathrm{Id}$ is the identity
matrix on $V$. Therefore, given a bounded $S_0 \subset S$ there is
$M>0$ such that for $h\in H$ and $s \in S_0$,
$$
D(h)-M\le D(sh)\le D(h)+M.
$$
This implies that for $s_1,s_2\in S$,
$$
\lambda(H_{T-M}[e,e])\le \lambda(H_T[s_1,s_2])\le \lambda(H_{T+M}[e,e]).
$$
Via Theorem \ref{thm: volume asymptotics1} we obtain that $\alpha
\equiv 1$, and the first assertion follows. 

For the second assertion, we proceed as before to calculate $\alpha$. 
Note that for $s\in S$ and $h\in H$ and using the $p$-norms on $V_H$ and $V_S$,
$$
\|\Psi(sh)\|=(\dim V_S\cdot \dim V_H)^{-1/p}\cdot\|\Psi_S(s)\|\cdot\|\Psi_H(h)\|.
$$
Therefore
\[
\begin{split}
H_T[s_0,s] & = \{h \in H: \|\Psi(s_0hs)\| <T\} \\
& = \{h \in H: \|\Psi(hs_0s)\| <T\}\\ 
& = \{h \in H: \|\Psi_H(h)\|
< c_1 T\},
\end{split}
\]
where 
$$c_1 = \frac{(\dim V_S\cdot \dim V_H)^{1/p}}{\|\Psi_S(s_0s)\|}.$$
It follows using Theorem \ref{thm: volume asymptotics1} that 
\[
\begin{split}
\alpha(s_0,s) & =\lim_{T \to \infty} \frac{\lambda \left(H_T[s_0,s]
\right)}{\lambda \left(H_T\right)}\\ 
& = \lim_{T \to \infty} \frac{C \, (\log c_1T)^{\ell} (c_1T)^m }{C \,
(\log T)^{\ell} \, T^m}\\  
& = c\|\Psi_S(s_0s)\|^{-m}
\end{split}
\]
with 
$$c=(\dim V_S\cdot \dim V_H)^{m/p}  .$$
This proves the second statement.
\end{proof}

\medskip
\begin{proof}[Proof of Corollary
\ref{cor_equi_appl2}]
Let $G=L \times L$, let $H$ be the diagonal embedding of $L$ in $G$,
and let   
$$D(g)= D(\ell_1, \ell_2)= \max \left\{1, \|\ell_1\|, \| \ell_2\| \right\}.$$
Our choices of $G, \, H, \, D$ show that the general
setup holds, and that $G, \, H, \, D$ are standard. Let $\mu =
\frac{1}{m(L/\Lambda)}\, m$ be the $L$-invariant probability measure on
$L/\Lambda$. Since $\Lambda \pi(g_0)$ is dense in $L/\Delta$, $H
\pi'(y_0)$ is dense in $G/(\Delta \times \Lambda)$, where $y_0 =
(g_0,e)$ and $\pi': G \to G/(\Delta \times \Lambda)$ is the quotient
map. By Theorem   
\ref{thm: using Ratner}, we find that $\left(** \right)$
holds. Since the $H$-action on $L/(\Delta \times \Lambda)$ is
isomorphic with the $H$-action on $Y=\Delta \backslash L \times
L/\Lambda$, defined via \equ{eq: G action, 
induction} (with $\Delta$, $\Lambda$, $H$ replacing $H$, $\Gamma$, $G$
respectively), it follows that $\left(* \right)$ holds for the $H$-action on $Y$. Using   
Proposition \ref{prop: 
verifying axioms 1} we find that I1, I2 and UC hold, and using
Proposition \ref{prop: induction} we obtain the desired result. 
\end{proof}
\ignore{
\begin{prop}\label{prop: implies cor}
Suppose that the general setup holds for $L$ and a distance function
$D: L\to\R_+$, and assume that $D$ is balanced. 
Let $\pi(g)\in L/\Delta$ be such that
$\overline{\Lambda\pi(g)}=L/\Delta$. Then for every $\varphi\in
C_c(L/\Lambda)$, 
$$
\frac{1}{m(L_T)}\sum_{\lambda\in \Lambda\cap L_T}
\varphi(\lambda^{-1}\pi(g))\mathop{\longrightarrow}_{T\to\infty} \frac{1}{m(L/\Lambda)}
\int_{L/\Lambda}\varphi\, dm. 
$$
\end{prop}

\begin{proof}
Let $X=L/\Delta$ and $\Lambda$ acts on $X$ by 
$$
x\cdot\lambda=\lambda^{-1}x,\;\; x\in X,\; \lambda\in\Lambda.
$$
Let $\tilde{Y}=X\times L$ and $Y=\tilde{Y}/\Lambda$ be defined as in Section \ref{section: induction}.
It is easy to check that the map 
$$
(x,l)\mapsto (lx, l),\;\; x\in X,\; l\in L,
$$
defines a $L$-equivariant homeomorphism of $Y$ and $L/\Delta\times L/\Lambda$.
Thus, the corollary follows from .
\end{proof}
 }

\medskip

We now derive Corollary \ref{cor: application torus} from a more
general result.
Suppose that the general setup holds for a balanced semisimple Lie
group $G$ and a 
distance function $D$. Let $G$
act smoothly on a Lie group $H$ by automorphisms.
Let $\Gamma$ and $\Lambda$ be lattices in $G$ and $H$ respectively (in
particular $G$ and $H$ are unimodular).
Denote by $m$ and $\mu$ Haar measures on $G$ and $H$ respectively,
normalized by the requirement that
$m(G/\Gamma)=\mu(H/\Lambda)=1$.
We assume that $\Gamma\cdot \Lambda \subset \Lambda$. Thus, $\Gamma$
acts on $H/\Lambda$, preserving $\mu$. 

\begin{thm}\label{cor_auto_equi}
For every $x_0\in H/\Lambda$ such that $\overline{\Gamma\cdot x_0}=H/\Lambda$ and
every $\varphi\in C_c(H/\Lambda)$,
\begin{equation}\label{eq_cor_apl1}
\frac{1}{\# \Gamma_T}\sum_{\gamma\in \Gamma_T} \varphi(\gamma^{-1}
x_0)\mathop{\longrightarrow}_{T \to \infty}  
\int_{H/\Lambda}\varphi\, d\mu.
\end{equation}
\end{thm}

\begin{proof}
By Proposition \ref{prop: verifying axioms 1}, I1, I2 and UC are satisfied. 
Let $X=H/\Lambda$ and $\Gamma$ act on $X$ on the right by $x\cdot
\gamma=\gamma^{-1}\cdot x$. 
Define $\tilde{Y}=G\times X$ and $Y=Y/\Gamma$ as in \S \ref{section: induction}.
The map $G\ltimes H\to G\times X$ induces a $G$-equivariant homeomorphism of 
$(G\ltimes H)/(\Gamma\ltimes \Lambda)$ and $Y$. Let $\pi': \tilde{Y} \to
Y$ be the quotient map. 
Since $\cl{G \pi' (e,x_0)} = \tilde{Y}$, 
the corollary follows from Theorem \ref{thm: using Ratner} and Proposition \ref{prop: induction}.
\end{proof}

\ignore{
\begin{remark} {\rm
Although we are interested in the action of a lattice $\Gamma$ on $H/\Lambda$ by automorphisms,
we have made an extra assumption that this action comes from the action of the ambient group $G$.
This assumption is important in the proof. In the case when $G$ is connected semisimple
Lie group that acts smoothly on a Lie group $H$, and $\Gamma$ is
an irreducible lattice in $G$ that satisfy some mild conditions, every action of $\Gamma$ 
by automorphism is induced by an action of $G$. This follows from Margulis Superrigidity Theorem.
}
\end{remark}
}

\begin{remark}
{\rm
Arguing as in the proof of Theorem \ref{thm: using Ratner}, it is
possible in  Corollary \ref{cor_equi_appl2} (respectively, Theorem
\ref{cor_auto_equi}) to replace $\lambda^{-1}$ with $\lambda$ (resp., 
$\gamma ^{-1}$ with $\gamma$). 
}
\end{remark}

\section{Examples}
\label{section: examples}
In this section we collect some examples showing that our hypotheses
are not automatically satisfied, and that they are important for the
validity of our results.

\subsection{Condition S}

Let 
$$S^1 \df \{z \in \C : |z|=1\}, \ \ \, G \df \prod_1^{\infty} S^1,\ \
\, H \df \prod_{1}^{\infty} \{\pm 
1\},$$ 
equipped with the Tychonov topology, and let $\tau: G \to \HG$ be the
natural map. Then $G$ is compact and $H$
is a compact subgroup. Any section $\sigma: \HG \to G$ induces a
section $\{\pm 1\} \backslash S^1 \to S^1$ in each factor. Let $\mathcal{U}$ be
an open subset of 
$\HG$. By the definition of the product topology on $G$,
$\tau^{-1}(\mathcal{U})$ contains a subset of the form $V_1 \times
\cdots \times V_r \times \prod_{r+1}^{\infty} S^1$, where $V_1,
\ldots, V_r$ are open, and in particular
contains a copy of $S^1$ in one factor.  Since there is no continuous
section $\{\pm 1\} \backslash S^1 \to S^1$, there is no continuous section defined on 
$\mathcal{U}$.

\subsection{Condition D2}
Let $G \df \SL(3, \R)$ and let 
$$H= \left \{
\left(\begin{matrix}e^t \cos t & -e^t \sin t & x \\
e^t \sin t & e^t \cos t & y \\
0 & 0 & e^{-2t} \end{matrix} \right) : x,y,t \in \R \right \}.
$$
In coordinates $t,x,y$, the left haar measure on $H$ is given by
$d\lambda = e^{2t}\, dx \, dy \, dt.$

\ignore{
For $c>1$ to be specified below, let $D$ be a matrix norm 
distance function given by 
$D(g) = \|g\|,$
where
$$\left \| \left (a_{ij} \right)\right\| =
\max \left\{c|a_{11}|, \sqrt{a_{13}^2 + a_{23}^2}, |a_{ij}|, (i,j)
\neq (1,1), (1,3), (2,3)  \right\}.  
$$
}

For a constant $c>1$ to be specified below, let $D$ be a matrix norm 
distance function given by 
$D(g) = \|g\|,$
where
$$\left \| \left (a_{ij} \right)\right\| \df 
\max \left\{ \sqrt{ca_{11}^2 + a_{12}^2},\, \sqrt{ca_{22}^2 +
a_{21}^2}, \, 
\sqrt{a_{13}^2+a_{23}^2} , \, |a_{31}|, \,  |a_{32}|, \, |a_{33}|\right\}.
$$

We first compute the volume $\lambda(H_T)$ along
two subsequences. In $(t,x,y)$ coordinates,
$H_{T}$ is given by the inequalities 
$$ \sqrt{x^2+y^2} < T, \ e^{-2t}<T,
$$
and 
\begin{equation}
\label{eq: spiral}
f_c(t)<T, \ \ \ \ \mathrm{where \ \ \ } f_c(t)=e^{t}\sqrt{c^2
\cos^2 t + \sin^2 t}.
\end{equation}

For $c>1$ small enough, $f_c(t)$ is a monotonically increasing
function of $t$. 
This means that \equ{eq: spiral} is satisfied along a ray of the form
$(-\infty, \tau)$,
where $T=f_c(\tau)$.

\begin{equation}
\label{eq: computation for one norm}
\begin{split}
\lambda(H_{T}) & = \int_{H_{T}} e^{2t} \, dt \, dx \, dy \\
& = \int_{-\log T/2}^{\tau}\,  \int_{x^2+y^2 < T^2} e^{2t} \, dt
\, dx \, dy \\ 
& = \pi T^2 \int_{-\log T/2}^{\tau} e^{2t}dt \\
& = \frac{\pi}{2}T^2 \left(e^{2\tau} - (2/T)^2  \right).
\end{split}
\end{equation}

Now let $\tau=\tau_n = 2 n \pi, $ for $n \in \N$. We
have for $T_n  = f_c(\tau_n) =  ce^{2n\pi}:$
\begin{equation}
\label{eq: one expression}
\lambda(H_{T_n}) = \frac{\pi}{2}T_n^2 \left( (T_n/c)^2 - (2/T_n)^2
\right ) \sim 
\frac{\pi}{2c^2} T_n^4 \ \ \ \ \mathrm{as} \ n \to \infty.
\end{equation}

Similarly, for $\tau_n \df (2n+1/2) \pi, \, S_n =
f_c(\tau_n)$  we have 
\begin{equation}
\label{eq: another expression}
\lambda(H_{S_n}) = \frac{\pi}{2}S_n^2 \left( S_n^2 -
(2/S_n)^2 \right ) \sim 
\frac{\pi}{2} S_n^4 \ \ \ \mathrm{as} \ n \to \infty.
\end{equation}

In particular the quantity $\lambda(H_T)$ is not asymptotic to a
function of the form $kT^4$ as $T \to \infty$ for any constant $k$.

Now taking 
$$g_1=g_2= \left(\begin{matrix} 0 & 1 & 0 \\
-1 & 0 & 0 \\
0 & 0 & 1 \end{matrix} \right),
$$
and computing the shape of $H_T[g_1^{-1}, g_2]$ in $(t,x,y)$
coordinates we find the role of sin and cos switched in \equ{eq:
spiral}. This implies that the roles of $T_n$ and $S_n$ are
reversed, that is,  
 
$$\lambda(H_{T_n}[g_1^{-1},g_2]) \sim \frac{\pi}{2} T_n^4, \ \ \ \  \ \ \ \
\lambda(H_{S_n}[g_1^{-1},g_2]) \sim  
\frac{\pi}{2c^2} S_n^4.$$

In particular,
$$\frac{\lambda(H_{T_n}[g_1^{-1},g_2])}{\lambda(H_{T_n})} \to_{n \to
\infty} c^2, \ \ \ \ \ \  
\frac{\lambda(H_{S_n}[g_1^{-1},g_2])}{\lambda(H_{S_n})} \to_{n
\to \infty} \frac{1}{c^2},
$$
and D2 is not satisfied.

\ignore{xxxxxxxxxxxxxxxxxxxxxxxxxxxxx leftovers

Define a norm on $\Mat_d(\R)$ by 
$$\|(a_{ij})\|^2 = ca_{11}^2 + \sum_{(i,j) \neq (1,1)} a_{ij}^2,
$$
where $c>1$ will be specified below. Let $D$ be
the corresponding distance function.  We claim that there is a
positive continuous function $c(t)$ and $d>1$ such that 
$$\lambda(H_T) \sim c(T) T^4 \ \ \ \mathrm{and \ } c(t) = c(dt).$$

Let 
$$R_{\theta}  = \left(\begin{matrix} \cos \theta & -\sin \theta &
0 \\
\sin \theta & \cos \theta & 0 \\
0 & 0 & 1 \end{matrix} \right),$$
let $\| \cdot \|_{\theta}$ be the norm on $\Mat_3(\R)$ defined by 
$$\|A \|_{\theta} = \|R_{-\theta} A \|,$$
and let $E \in \Mat_3(\R)$ be the matrix with 1 in the (1,1) and (2,2)
entries and 0 in the other entries.

We claim that for any $x , y \in \R$ and $\theta \in [0, 2\pi)$ there
is a constant $C$ such that 
for $ T \gg 0$,
$$\left| \sup \{t \in \theta + 2\pi \Z : h(t,x,y) \in H_T\} -
\frac{\log T}{\log\|E\|_{\theta}} \right| < Ce^{-t}.$$ 

}
\medskip

To show that condition D2 is indeed necessary for the validity of our
results, we have the following:

\begin{thm}
Let $H, \, G, \, \| \cdot \|, \, \{S_n\}, \, \{T_n\}$ be as in the
above example. Then  there are two equivalent but different measures 
$\nu_1, \, \nu_2$ on $\HG$ such that for any
lattice $\Gamma$ in $G$, any $x_0 \in \HG$ such that $\HG 
= \cl{x_0 \, \Gamma},$ and any bounded $A \subset \HG$
with $\nu_i(\partial A)=0$ we have 
\begin{equation}
\label{eq: one asymptotics}
\frac{N_{S_n}(A, x_0)}{\lambda(H_{S_n})}
\mathop{\longrightarrow}_{n\rightarrow\infty} 
\nu_1(A)
\end{equation}
and
\begin{equation}
\label{eq: another asymptotics}
\frac{N_{T_n}(A, x_0)}{\lambda(H_{T_n})}
\mathop{\longrightarrow}_{n\rightarrow\infty} 
\nu_2(A).
\end{equation}
\end{thm}

\begin{proof}[Sketch of proof]
Suppose
$\Gamma$ is given. The subgroup of $H$ given by the requirement $t=0$
is the unipotent radical of a parabolic subgroup of $G$. Using the methods of 
\cite{Nimish: proc ind} one can show that condition
$\left(** \right)$ holds for any $g_0 \in G$ for
which  $\cl{x_0 \, \Gamma} = \HG$ (where $x_0 = \tau(g_0)$) . 

We now claim that  D1 is satisfied. Let $\mathcal{N}$ be
any bounded collection of norms on $\Mat_3(\R)$. For $| \cdot | \in
\mathcal{N}$, let 
$$H_T^{|\cdot |} = \{h \in H: |h|<T\}.$$
  We will show that there is a function $f(|\cdot |, T)$,
which depends continuously on both its arguments, such that 
\begin{equation}
\label{eq: will show asympt}
\lambda\left(H_T^{|\cdot |}\right) \sim f(|\cdot |,  T) \, T^4
\end{equation}
and such that for all $| \cdot | \in \mathcal{N}$ and all $T>0$, 
\begin{equation}
\label{eq: cyclic}
f(|\cdot |, aT)=f(|\cdot |, T), \ \ \ \mathrm{where \ } a=e^{2\pi}.
\end{equation}
It is easily seen that this implies D1.

For 
$$\omega = (x_0, y_0, s_0)  \in S^2 = \{(x_0, y_0, s_0):
x_0^2+y_0^2+s_0^2=1\}$$
and for $T, \, r>0$, write
$$h(\omega, r) = \left(\begin{matrix} rs_0 \cos \log (rs_0) & -rs_0
\sin \log (rs_0)  & rx_0 \\
rs_0 \sin \log (rs_0) & rs_0 \cos \log (rs_0) & ry_0 \\
0 & 0  & (rs_0)^{-2} \end{matrix} \right) 
$$
and 
$$\tau(\omega, T, | \cdot |) = \{r \geq 0: |h(\omega, r) | < T\}.$$
Then, letting $d \omega$ denote the standard volume form on $S^2$, we
have  
\begin{equation*}
\lambda\left(H^{|\cdot |}_T\right) = \int_{S^2} \, \int_{\tau(\omega,
T, |\cdot |)} s_0r^3 \, 
dr \, 
d\omega. 
\end{equation*}
Also write
\[
\begin{split}
E_{\theta} = \left(\begin{matrix} \cos \theta & -\sin \theta & 0 \\
\sin \theta & \cos \theta & 0 \\ 0 & 0 & 1 \end{matrix} \right)
, \ \ & E_1 =  \left(\begin{matrix} 0 & 0 & 1 \\
0 & 0 & 0 \\ 0 & 0 & 0 \end{matrix} \right)
, \\
E_2 =  \left(\begin{matrix} 0 & 0 & 0 \\
0 & 0 & 1 \\ 0 & 0 & 0 \end{matrix} \right)
, \ \ &  E_3  =  \left(\begin{matrix} 0 & 0 & 0 \\
0 & 0 & 0 \\ 0 & 0 & 1 \end{matrix} \right) .
\end{split}
\]
Then we have, for $\omega = (x_0, y_0, s_0) \in S^2$,  
$$\tau(\omega, T, |\cdot |) = \left\{r \geq 0: \left |r s_0 E_{\log(rs_0)}
+ r x_0 E_1 +ry_0E_2 + \frac{1}{(rs_0)^{2}}E_3 \right| < T
\right\}.$$
Define
\[
\begin{split}
\til \tau(\omega, T, |\cdot |) & = \left\{r \geq 0: \left |r s_0
E_{\log(rs_0)} +r x_0 E_1 +ry_0E_2 \right| < T
\right\} \\ 
& =  \left\{r \geq 0:  r  < T / |E(r, \omega) | 
\right\},
\end{split}
\]
where 
$$ E(r, \omega) = s_0 E_{\log(rs_0)}
+x_0 E_1 +y_0E_2.$$
By the computation \equ{eq: computation for one norm}, and since the
collection  
$\mathcal{N}$ is bounded, there are constants $c_1, c_2$ such that for
all $| \cdot | \in \mathcal{N}$:
\begin{equation}
\label{eq: grows like T^4}
0<c_1 \leq \liminf_{T \to \infty} \frac{\lambda(H^{|\cdot
|}_T)}{T^4} \leq 
\limsup_{T \to \infty} \frac{\lambda(H^{|\cdot |}_T)}{T^4} \leq
c_2<\infty.
\end{equation}
 
Given $\vre>0$, there is a constant $C$ such that for each $\omega$
with $s_0 \geq \vre$, and each $| \cdot  
| \in \mathcal{N}$, the symmetric difference
of the sets $\til 
\tau(\omega, T, |\cdot |)$ and $ \tau(\omega, T, |\cdot |)$ is
contained in an interval of length $C/T$. 
Using this and \equ{eq: grows like T^4} it is not hard to show that 
\begin{equation}
\label{eq: similar growth}
\lambda\left(H_T^{|\cdot |}\right)
 \sim \int_{S^2}\, \int_{\til \tau(\omega, T, |\cdot |)} r^3 \, dr \, d\omega. 
\end{equation}

Since $E(r, \omega)=E(ar, \omega)$ we have, for all $T>0$ and $\omega
\in S^2$,  
$$\til \tau(\omega, aT) = a\til\tau(\omega, T),$$
therefore by a change of variables, for any $\omega \in S^2$:
\begin{equation}
\label{eq: invariance}
\int_{\til \tau(\omega, aT)} r^3 \, dr \, d\omega = 
a^4 \int_{\til \tau(\omega, T)} r^3 \, dr \, d\omega,
\end{equation}
hence
$$f( |\cdot |, T) = \frac{\int_{S^2} \, \int_{\til \tau(\omega, T,
|\cdot |)} r^3 \, 
dr \, 
d\omega}{T^4}$$
depends continuously on both its
parameters and satisfies \equ{eq: cyclic}. 
Now \equ{eq: will show asympt} follows from  \equ{eq: similar growth}
and \equ{eq: invariance}.  

\ignore{
\[
\begin{split}
\lim_{n \to \infty} \frac{\lambda(H_{a^nT_0})}{(a^nT_0)^4} & = \lim_{n
\to \infty} 
\frac{\int_{S^2}\, \int_{\til \tau(\omega, a^nT_0)} r^3 \, dr \, d\omega
}{(a^nT_0)^4} \\
& = \frac{\int_{S^2}\, \int_{\til \tau(\omega, T_0)} r^3 \, dr \, d\omega
}{T_0^4} \df f(T_0).
\end{split}
\] 

Let $f: \R_+ \to \R_+$ be defined by the conditions
$$ f(t) =  \frac{\int_{S^2}\, \int_{\til \tau(\omega, t)} r^3 \, dr \,
d\omega }{t^4} 
$$
when $t \in [1,a]$ and $f(at)=f(t)$. Then $f$ depends continuously on
$t$ and on the norm $\| \cdot \|$, and it
follows from the  above computation that 
$$\lambda(H_T) \sim f(T) T^4.$$

In particular it follows that condition D1 is satisfied.

 Since $f_c(t)$ is a
monotonically increasing function of $t$, it has an inverse
$f_c^{-1}(T),$ and from \equ{eq: spiral} we have 
$$T = f_c(\tau) \ \ \Longleftrightarrow \ \ \tau = \log T + d(T),$$
where
$$d(T) = \sqrt{c^2 \cos^2 (f_c^{-1}(T)) + \sin^2
(f_c^{-1}(T))} \in [1,c]. 
$$
It is easily checked using \equ{} that 
$$\lambda(H_T) \sim d(T) T^4.$$
}

Now we may apply Theorem \ref{thm: duality} to obtain, for 
any $\varphi \in C_c(\HG)$,
$$S_{\varphi, x_0} (T) \sim \til S_{\varphi, x_0}(T). $$

\ignore{Note that in the current setup condition D1 does not hold, and so
Theorem \ref{thm: duality} cannot be applied directly. However,
examining the proof of the theorem one sees that is enough to assume
the following 
weakened version of hypothesis D1:

{\em 
For any bounded $B \subset G$
and any $\vre>0$ there are $T_0$ and $\delta$ such that for any
$T>T_0$ there is $ \delta_0 \geq \delta$ such that 
$$\left(H_{(1+\delta_0)T} [g_1,g_2] \right) \leq (1+\vre)\lambda \left(H_T
[g_1,g_2] \right).
$$
}
Moreover, this weakened hypothesis can be verified in our example.

Now let $\{S_n\}$ and $\{T_n\}$ be as in the above example. Using a
standard compactness argument, one can, by passing to subsequences
$\{s_n \} \subset \{S_n\}$ and $\{t_n\} \subset \{T_n\}$ ensure that
the limits 
$$\alpha_1(g_1, g_2) = \lim_{n \to \infty}
\frac{\lambda(H_{s_n}[g_1^{-1},g_2])}{\lambda(H_{s_n})}
, \ \alpha_2(g_1, g_2) =
\lim_{n \to \infty}
\frac{\lambda(H_{t_n}[g_1^{-1},g_2])}{\lambda(H_{t_n})}
$$ exist and are positive for all $g_1, g_2$. As was shown above, there are
$g_1, g_2$ for which $\alpha_1(g_1, g_2) \neq \alpha_2(g_1, g_2)$ and in
fact it is not hard to see that $\nu_1 \neq \nu_2,$ where $\nu_i$ are
defined  via \equ{eq_nu_g}, using $\alpha = \alpha_i$. 
}
Fix a sequence $t_n = e^{2\pi n} t_0$, for any $t_0>0$. The ratio 
$$\alpha_{t_0} (g_1, g_2) = \lim_{n \to \infty} \frac{\lambda \left
(H_{t_n}[g_1^{-1}, g_2] \right)}{\lambda \left(H_{t_n} \right) }$$ 
exists, is positive, and
depends continuously on $g_1, g_2$ by \equ{eq: will show asympt}. Let
$\nu_{t_0}$ be 
defined by \equ{eq_nu_g}, but using $\alpha_{t_0}$ in place of $\alpha$. 
Now repeating the arguments of \S \ref{section: describing the 
limit},   but taking limits as $n \to \infty$, we obtain the
conclusion of Corollary \ref{eq: indicator sets}. In particular,
taking in turn 
$t_n=S_n, \, t_n=T_n$ we obtain \equ{eq: one
asymptotics} and \equ{eq: another asymptotics}. Note that the two
limiting measures in this case are different by \equ{eq: one
expression} and \equ{eq: another expression}. 
\end{proof}

\subsection{Non-balanced semisimple groups}
\label{prop_balan_crit}
We now construct a non-balan\-ced semisimple group
and show that Theorem
\ref{thm: using Ratner} fails for this group, 
i.e. find an action of this group on a homogeneous space which does
not satisfy 
$\left(**
\right)$. 

Let 
$$H_1 = H_2 = \SL(2,\R), \ \ \ \ H=H_1 \times H_2.$$

Abusing notation, we consider $H_1$ and $H_2$ as subgroups of $H$.

Denote by $A_i$ a Cartan subgroup of $H_i$, $\goth{a}_i$ its
(one-dimensional) Lie
algebra, $Y_i$ a generator of $\goth{a_i}$. Then
$A=A_1A_2$ is a Cartan subgroup of $H$. We will write $(s_1, s_2)= s_1
Y_1+s_2Y_2 \in \goth{a}$. 
With respect to these coordinates, a root system of $(H, A)$ is
$\{\pm\alpha_1,\pm\alpha_2\}$, 
where
$$
\alpha_i(s_1, s_2)=2s_i
$$
and
$$
\goth{a}^+=\{(s_1, s_2): s_1 \geq 0, s_2 \geq 0\}.
$$

In the notation of \S 7 we have 
$$\rho = \frac12 \left(\alpha_1+\alpha_2 \right), \ \ \ \ \ \beta_i =
\frac{Y_i}{2}, \ i=1,2.
$$


We consider $H$ as a subgroup of $G = \textrm{SL}(2 \ell,\mathbb{R})$
where $H$ is embedded in $G$ via the
tensor product  
of irreducible representations of $H_1$ and $H_2$ of dimensions
2 and $\ell>2$ respectively. 
Note that the set of weights of $\mathfrak a$ is
$$\left\{(s_1, s_2) \mapsto is_1+js_2: \,
i \in \{\pm 1\}, \, j \in  \{1- \ell, 3 -\ell,\ldots,\ell-1\}\right\}, 
$$
and a highest weight corresponding to the choice of $\goth{a}^+$ is
$$\lambda_1(s_1, s_2) = s_1+(\ell -1)s_2.$$

We have 
$$\lambda_1(\beta_1) = \frac12 < \frac{\ell -1}{2} =
\lambda_1(\beta_2),$$
so condition {\bf G} is satisfied, and by Proposition \ref{prop: G and
balanced}, $H$ is not balanced, that is, there is a bounded open $L
\subset H_2$ such that 
$$
c=\limsup_{T \to \infty} \frac{\lambda(H_T^L)}{\lambda(H_T)}>0, \ \ \ \
\mathrm{where \ } H_T^L = \{h=(h_1, h_2): \|h\|<T, \, h_2 \in L\}
$$
for some (any, see Proposition \ref{prop: contribution to volume})
norm $\| \cdot \|$ on $\Mat_{2 
\ell}(\R)$. 
We have:
\begin{claim}\label{claim: lattice exists}
There is a $\Q$-subgroup $M$ of $G$ containing $H_1$ and an element
$m \in M$ such that $m^{-1} H m$ is not contained in any
proper $\Q$-subgroup of $G$. In particular, setting $x = \pi(m)$
and $\Gamma = \SL(2\ell , \Z)$ we
have that $H_1x \subset M \pi(e) = \cl{M\pi(e)}$ and $\cl{Hx} =
G/\Gamma$. 
\end{claim}

Assuming the claim is true, let $X_0 = LM\pi(e) \subset
G/\Gamma$. This is the closure of a locally closed submanifold of
$G/\Gamma$, being the image of 
$M/ \Delta \times L \subset H$ under the proper map $(g \Delta, \ell) \mapsto \ell
gx$, where $\Delta = M \cap m \Gamma m^{-1}.$ Since $\dim M <  \dim G$, we
have 
$m'(X_0)=0$. Therefore, for any compact $K_0 \subset G/\Gamma$ we can
find $\varphi \in C_c(G/\Gamma)$ such 
that $\int_{G/\Gamma} \varphi \, dm' < c/2$ and $\varphi|_{X_0 \cap K_0} \equiv
1.$ It follows from the non-divergence results used in \cite{Nimish:
proc ind} that we can make  
$K_0$ large enough so that 
$$\limsup_{T \to \infty} \frac{\lambda \left\{h \in H_T: hx \notin K_0
\right\}}{\lambda(H_T)} < c/2.
$$
This yields
\[
\begin{split}
& \liminf_{T \to \infty} \frac{1}{\lambda(H_T)} \int_{H_T} \varphi(hx) \, d\lambda(h) \\
& \geq \liminf_{T \to \infty} \frac{\lambda \left\{h \in H_T: hx \in K_0
\cap X_0\right\}}{\lambda(H_T)} \\
& \geq
\liminf_{T \to \infty} \frac{\lambda(H_T^L)}{\lambda(H_T)}  -
\limsup_{T \to \infty}  \frac{\lambda \left\{h \in H_T: hx \notin K_0
\right\}}{\lambda(H_T)}\\
& \geq c -c/2 >
\int_{G/\Gamma} \varphi \, dm',
\end{split}
\]
and $\left( ** \right)$ fails. 
 
\medskip

It remains to prove Claim \ref{claim: lattice exists}. Let $\{u_1,
u_2\}$ and $\{v_1, \ldots, v_{\ell} \}$ be the standard bases of $\R^2$ and
$\R^{\ell}$ respectively, and let $\mathcal{B} =\{u_i \otimes v_j\}$, a basis of
$\R^2 \otimes \R^{\ell} = \R^{2\ell}.$ The $\Q$-structure on $G$ is
defined via $\mathcal{B}$, and the $H_1$-action (respectively, the
$H_2$) action on $\R^{2\ell}$ 
is induced by its action on the $u_i$'s (respectively, $v_j$'s). Now
let $M$ be the subgroup of $G$ leaving invariant each of the subspaces
$V_j = \spa \{u_1 \otimes v_j, u_2 \otimes v_j\}$. Clearly each $V_j$
is $H_1$-invariant and hence $H_1 \subset M$. It is also clear that
$M$ is defined over $\Q$ and hence $M\pi(e)$ is closed. It remains to
show that there exists $m \in M$ so that $m^{-1}Hm$ is not contained
in any proper $\Q$-subgroup of $G$. Suppose otherwise; since the
number of $\Q$-subgroups of $G$ is countable, this would imply that
there is a fixed proper $\Q$-subgroup $T \subset G$ such that for all
$m \in M, \, m^{-1}Hm \subset T$. However it is not difficult to show
(we omit the computation) that the set $\{m^{-1}hm : m \in M, h \in
H\}$ generates $G$, and this is a contradiction.

\subsection{A simple case revisited}
Our results also enable us to generalize Ledrappier's result,
discussed in \S 1.1, to general norms. Namely we have:

\begin{thm}
Let $\Gamma$ be a lattice in $\hbox{\rm SL}(2,\R)$ and let $\|\cdot\|$
be a norm on $\Mat_2(\R)$. 
Suppose that $v\in V=\R^2$ satisfies $\overline{v\cdot
\Gamma}=V$, and let $dw$ denote Lebesgue measure on $\R^2$. Then for
every $\varphi\in C_c(V)$,  
\begin{equation}\label{eq_ledd}
S_{ \varphi, v} (T) \sim \left(c_{\Gamma} \, \int_V \varphi(w)
\alpha_v(w)\, 
dw\right)\, T, 
\end{equation}
where $c_{\Gamma}>0$ is a constant depending on $\Gamma$ and 
$$
\alpha_v(w)=
\left\|\left( \begin{array}{cc}
0 & 0 \\
1 & 0
\end{array}
\right )\right\| \, \cdot \,
\left\|\left(
\begin{tabular}{cc}
$-v_2w_1$ & $-v_2w_2$\\
$v_1w_1$ & $v_1w_2$
\end{tabular}
\right)\right\|^{-1}.
$$
\end{thm}

\begin{proof}
Via the map 
$$
\tau: G\to V, \ \ \ \tau :g \mapsto (1,0)\cdot g,
$$
the space $V \sm \{0\}$ is identified with $H\backslash G$, where
$G=\hbox{SL}(2,\R)$ and 
$$H=\left\{u_t=\left( 
\begin{tabular}{cc}
$1$ & $0$\\
$t$ & $0$
\end{tabular}
\right):t\in\R\right\}.$$
Note that haar measures on $H$ and $G$ are only defined up to a
constant multiple; we equip $H$ with the haar measure $dt$, and choose
haar measure $\mu$ on $G$ so that $\nu_{\HG}$ is Lebesgue measure. This
induces a choice $\mu'$ of $G$-invariant measure on $G/\Gamma$. Since
our results were formulated for the choice making $\mu'$ a probability
measure, we set $c_{\Gamma} = \frac{1}{\mu'(G/\Gamma)},$ and $m' =
c_{\Gamma} \mu'$.


We check that the hypotheses of Theorems \ref{thm: duality} and
\ref{thm: 
identifying the limit} hold 
in this case. 
For $g_1,g_2\in G$,
$$
H_T[g_1,g_2]=\{u_t:\, \|a+tb\|<T\}
$$
for $a=g_1g_2$ and $b=g_1\left(
\begin{tabular}{cc}
$0$ & $0$\\
$1$ & $0$
\end{tabular}
\right)g_2$.
This implies 
$$
\left\{u_t:\, |t|<\frac{T- \|a\|}{\|b\|}\right\}\subset
H_T[g_1,g_2]\subset\left\{u_t:\, |t|<\frac{T+ \|a\|}{\|b\|}\right\}.
$$
Hence
\begin{equation}\label{eq_l_asy}
\lambda(H_T[g_1,g_2])\sim  \frac{2T}{\|b\|},
\end{equation}
with uniform convergence for $g_1, g_2$ in a compact subset of $G$. 
In particular, hypotheses D1 and D2 are satisfied. Hypothesis
$\left(**\right)$ 
follows from the equidistribution of the horocycle flow \cite{Dani
Smillie:
horocycle flow}.
Thus, Corollary \ref{cor: limit for indicator sets} applies,
and (\ref{eq_ledd}) holds for the function $\alpha_v(w)$ as in
Proposition 
\ref{prop: measure well defined}(iv).

To calculate $\alpha_v(w)$, define $V_0 = \left\{(x_1, x_2) \in V: x_1
\neq 0 \right\}$
and assume that $v \in V_0$ and $\supp \, \varphi \subset
V_0$. Consider a measurable section 
$
\sigma: V \to G$ whose restriction to $V_0$ is continuous and defined by
$$\sigma(x_1,x_2) =  \left(
\begin{tabular}{cc}
$x_1$ & $x_2$\\
$1$ & $\frac{x_2+1}{x_1}$
\end{tabular}
\right).
$$

By \equ{eq_dual_as} and (\ref{eq_l_asy}),
\begin{eqnarray*}
\alpha_v(w)&=&\lim_{T\to\infty}
\frac{\lambda(H_T[\sigma(v)^{-1},\sigma(w)])}{\lambda(H_T)}\\ 
&=&\left\|\left(\begin{tabular}{cc}
$0$ & $0$\\
$1$ & $0$
\end{tabular}
\right)\right\|\cdot \left\|\sigma(v)^{-1}\left(
\begin{tabular}{cc}
$0$ & $0$\\
$1$ & $0$
\end{tabular}
\right)\sigma(w)\right\|^{-1}
\end{eqnarray*}
for $v,w\in V_0$.
This implies the corollary. If $v \notin V_0$ or $\supp \, \varphi
\not \subset V_0$ we complete the proof by taking a different section
in the obvious way.
\end{proof}

\end{document}